\documentclass[10pt]{amsart}
\usepackage[left=1in,top=1in,right=1in,bottom=1in]{geometry}
\usepackage{amsfonts}
\usepackage{enumerate,subfigure,color}
\usepackage{amssymb,latexsym, amsmath,graphicx}%,epsf}
\usepackage{epstopdf}
\usepackage{framed}
\usepackage{multirow}

\usepackage[round,sort&compress]{natbib}

\usepackage{hyperref}

\usepackage{enumitem}
\usepackage[utf8]{inputenc}
\usepackage{tikz}
\usetikzlibrary{shapes.geometric, arrows}

\tikzstyle{title} =  [rectangle, minimum width=3cm, minimum height=1cm,text centered, text width=5cm, draw=black, fill=purple!20]

\tikzstyle{startstop} = [rectangle, rounded corners, minimum width=3cm, minimum height=1cm,text centered, text width=5cm, draw=black, fill=orange!10]
\tikzstyle{io} =  [rectangle, rounded corners, minimum width=3cm, minimum height=1cm,text centered, text width=5cm, draw=black, fill=blue!20]
\tikzstyle{process} = [rectangle, rounded corners, minimum width=3cm, minimum height=1cm, text width=5cm, text centered, draw=black, fill=yellow!30]
\tikzstyle{algend} = [rectangle, rounded corners, minimum width=3cm, minimum height=1cm,text centered, text width=4.5cm, draw=black, fill=green!10]

\tikzstyle{decision} = [diamond, minimum width=3cm, minimum height=1cm, text centered, text width=3cm,draw=black, fill=green!30]

\tikzstyle{arrow} = [thick,->,>=stealth]
\tikzstyle{startstop2} = [rectangle, rounded corners, minimum width=3cm, minimum height=1cm,text centered, text width=8cm, draw=black, fill=orange!10]
\tikzstyle{io2} =  [rectangle, rounded corners, minimum width=3cm, minimum height=1cm,text centered, text width=7.5cm, draw=black, fill=blue!20]
\tikzstyle{process2} = [rectangle, rounded corners, minimum width=3cm, minimum height=1cm, text width=6cm, text centered, draw=black, fill=yellow!30]
\tikzstyle{algend2} = [rectangle, rounded corners, minimum width=3cm, minimum height=1cm,text centered, text width=10cm, draw=black, fill=green!10]

\tikzstyle{algend3} = [rectangle, rounded corners, minimum width=3cm, minimum height=1cm,text centered, text width=12.5cm, draw=black, fill=green!10]

\newcommand{\dbar}{\bar{\partial}}

\newcommand{\T}{{\mathbf{t}}}

\newcommand{\R}{{\mathbb R}}

\newcommand{\C}{{\mathbb C}}

\DeclareMathOperator{\de}{\partial}

\DeclareMathOperator{\dez}{\de_z}
\DeclareMathOperator{\dbarz}{\dbar_z}
\DeclareMathOperator{\dbark}{\dbar_k}

\DeclareMathOperator{\by}{\times}
\DeclareMathOperator{\bndry}{\partial\Omega}

\DeclareMathOperator{\Sexp}{S^{\mbox{\tiny{\textbf{exp}}}}}
\DeclareMathOperator{\Psiexponetwo}{\Psi_{12}^{\mbox{\tiny{\textbf{exp}}}}}
\DeclareMathOperator{\Psiexptwoone}{\Psi_{21}^{\mbox{\tiny{\textbf{exp}}}}}
\DeclareMathOperator{\Sexponetwo}{S_{12}^{\mbox{\tiny{\textbf{exp}}}}}
\DeclareMathOperator{\Sexptwoone}{S_{21}^{\mbox{\tiny{\textbf{exp}}}}}

\DeclareMathOperator{\Qexponetwo}{Q_{12}^{\mbox{\tiny{\textbf{exp}}}}}
\DeclareMathOperator{\Qexptwoone}{Q_{21}^{\mbox{\tiny{\textbf{exp}}}}}
\DeclareMathOperator{\Mexp}{M^{\mbox{\tiny{\textbf{exp}}}}}
\DeclareMathOperator{\sigexp}{\sigma^{\mbox{\tiny{\textbf{exp}}}}}
\DeclareMathOperator{\muexp}{\mu^{\mbox{\tiny{\textbf{exp}}}}}

\DeclareMathOperator{\gamexp}{\gamma^{\mbox{\tiny{\textbf{exp}}}}}

\DeclareMathOperator{\Lambdaref}{\Lambda_{\mbox{\tiny{\textbf{ref}}}}}

\DeclareMathOperator{\tdiff}{\mathbf{t}^{\mbox{\tiny{\textbf{diff}}}}}
\DeclareMathOperator{\texp}{\mathbf{t}^{\mbox{\tiny{\textbf{exp}}}}}

\DeclareMathOperator{\Spsidiff}{S^{\mbox{$\Psi$\tiny{\textbf{diff}}}}}
\DeclareMathOperator{\Sdiffonetwo}{S_{12}^{\mbox{\tiny{\textbf{diff}}}}}
\DeclareMathOperator{\Sdifftwoone}{S_{21}^{\mbox{\tiny{\textbf{diff}}}}}
\DeclareMathOperator{\Psidiffonetwo}{\Psi_{12}^{\mbox{\tiny{\textbf{diff}}}}}
\DeclareMathOperator{\Psidifftwoone}{\Psi_{21}^{\mbox{\tiny{\textbf{diff}}}}}

\DeclareMathOperator{\gammadiff}{\gamma^{\mbox{\tiny{\textbf{diff}}}}}
\DeclareMathOperator{\sigmadiff}{\sigma^{\mbox{\tiny{\textbf{diff}}}}}
\DeclareMathOperator{\mudiff}{\mu^{\mbox{\tiny{\textbf{diff}}}}}

\usepackage{soul,color}
\soulregister\cite7
\soulregister\ref7
\soulregister\pageref7
\newcommand{\trev}[1]{{#1}}

\newcommand{\tRev}[1]{{#1}}
\begin{document}
%--------------------------------------------------------------------
\bibliographystyle{plainnat} % for Phys Meas.

%\title[Robust D-bar (Preliminary)]{Robust Computation of 2D EIT Absolute Images with D-bar Methods}
%\title{Robust Computation of 2D EIT Absolute Images with D-bar Methods}
\title[Robust Absolute EIT Imaging with D-bar Methods]{Robust computation in 2D absolute EIT (a-EIT) using D-bar methods with the ``exp'' approximation}
\author{S.~J. Hamilton, J.~L. Mueller, and T.~R.~Santos}

\date{June 2018}

\thanks{S.~J. Hamilton is with the Department of Mathematics, Statistics, and Computer Science; Marquette University, Milwaukee, WI 53233 USA,  email: \texttt{sarah.hamilton@marquette.edu}}

\thanks{J.~L. Mueller is with the Department of Mathematics and School of Biomedical Engineering, Colorado State University, Fort Collins, CO 80523 USA}

\thanks{T.~R~Santos is with the Department of Mechanical Engineering, University of S\~ao Paulo, Brazil}

%--------------------------------------------------------------------
\begin{abstract}
{\it Objective:} Absolute images have important applications in medical Electrical Impedance Tomography (EIT) imaging, but the traditional minimization and statistical based computations are very sensitive to modeling errors and noise. In this paper, it is demonstrated that D-bar reconstruction methods for absolute EIT are robust to such errors.  {\it Approach:} The effects of errors in domain shape and electrode placement on absolute images computed with 2-D D-bar reconstruction algorithms are studied on experimental data.  {\it Main Results:} It is demonstrated with tank data from several EIT systems that these methods are quite robust to such modeling errors, and furthermore the artefacts arising from such modeling errors are similar to those occurring in classic time-difference EIT imaging.  {\it Significance:} This study is promising for clinical applications where absolute EIT images are desirable, but previously thought impossible.
\end{abstract}
%--------------------------------------------------------------------
%\pacs{PUT NUMBERS IN}
%\keywords{magnetic moment, solar neutrinos, astrophysics}
%\submitto{\PM}
%--------------------------------------------------------------------
\maketitle % uncomment if separate title page is required
%--------------------------------------------------------------------
\section{Introduction}\label{sec:intro}
%--------------------------------------------------------------------
Absolute images, defined as images computed independently of a reference set of data or images, may have important applications in medical EIT for distinguishing between pulmonary abnormalities that appear the same in difference images, such as pneumothorax versus hyperinflation or atelectasis versus pulmonary edema, and for the classification of a breast lesion as a cyst or tumor, benign or malignant.  \tRev{For thoracic imaging, 2D images are particularly relevant for ARDS patients, for which a single cross-sectional image well represents the heterogeneous mechanical properties of dependent and non-dependent lung regions [\cite{ElDash2016}].}  However, due to the sensitivity of the severely ill-posed inverse conductivity problems to modeling errors such as errors in electrode locations, contact impedance, and domain shape, absolute images with good spatial resolution and few artefacts are notoriously difficult to compute.  Minimization-based methods, such as Gauss-Newton approaches require an accurate forward model that can predict the voltages on the electrodes from a candidate conductivity distribution with high precision.  A typical forward model computed with the finite element method (FEM) requires precise knowledge of boundary shape and electrode positions, as well as a very high number of elements to achieve accuracy, and so solving the forward problem at each iteration in the Gauss-Newton method has high computational cost.  When real-time imaging is desirable, such as in pulmonary applications, this computational burden and imprecise boundary knowledge pose significant challenges.  

D-bar methods are direct (non-iterative) and therefore do not require repeated high-accuracy simulations to compute an image.  Their real-time capabilities have recently been demonstrated in [\cite{Dodd2014}].  
Here, we demonstrate that high quality absolute images can be computed quickly from D-bar methods, and that they are robust in the presence of the intrinsic system noise, as well as errors in electrode placement and domain shape.  Absolute images are computed from experimental tank data using three D-bar methods, and the effects of errors in electrode location and domain shape are compared.  The D-bar methods studied here are the $\texp$ implementation of the D-bar method for real-valued conductivities based on the global uniqueness proof [\cite{Nachman1996}] with further developments and implementations in [\cite{Siltanen2000,Isaacson2004,Knudsen2009}], and two formulations of the D-bar method for complex-valued conductivities, \tRev{also with the exp-approximation} which are based on [\cite{Francini2000,Hamilton2012, Herrera2015}].

The first absolute images computed with a D-bar method were presented in [\cite{Isaacson2004}], where they were compared to absolute images computed using the NOSER [\cite{Cheney1990a}] algorithm.  The images were found to be nearly free of the high conductivity artefacts along the boundary caused by the presence of the electrodes, but the target positions (agar heart and lungs) were of lower spatial accuracy than those of the NOSER images.  However, the conductivity values computed by the D-bar method were more accurate in value and dynamic range.  In [\cite{Murphy2009}]  the effects of errors in input currents, output voltages, electrode placement, and domain shape modeling in the D-bar method for real-valued conductivities were studied on simulated data and the method was demonstrated to be quite robust.  \tRev{Additionally, the D-bar reconstruction methods for anisotropic conductivities, which uniquely recover $\sqrt{\det\sigma}$ up to a change of coordinates, are identical to their isotropic counterparts [\cite{Henkin2010}, \cite{Hamilton2014a}].  As incorrect domain modeling is a known source of anisotropy, leading to EIT data that would only arise from an anisotropic conductivity even if the true conductivity is isotropic, these results help explain the robustness observed for D-bar based reconstruction methods for EIT.}  This paper contains the first study of the effects of errors in electrode placement and domain shape modeling when using experimental data with D-bar methods and handles permittivity, as well as conductivity, imaging.

The paper is organized as follows.  Section \ref{sec:methods} provides overviews of the equations in the D-bar methods studied here with computational details of the implementation of the D-bar method for complex conductivities in [\cite{Herrera2015}] given in Section \ref{sec:CompDetails}.  Section~\ref{sec:results} presents absolute and time-difference EIT reconstructions using experimental data from three EIT systems.  Results are compared for various errors in electrode placement and domain shape modeling \trev{and a} discussion of the results \trev{included. Conclusions are} stated in Section~\ref{sec:conclusions}.

%--------------------------------------------------------------------
\section{Methods}\label{sec:methods}
%--------------------------------------------------------------------
Letting $\Omega\subset\R^2$ denote a simply connected domain with Lipschitz boundary, the electric potential inside $u$ inside $\Omega$ can be modeled by a generalized Laplace equation, also called the {\it admittivity equation}
\begin{equation}\label{eq:admitt}
\nabla\cdot\gamma(z)\nabla u(z) = 0,\quad z\in\Omega\subset\R^2,
\end{equation}
where $\gamma(z)=\sigma(z)+i\omega\epsilon(z)$ denotes the admittivity with conductivity $\sigma$, permittivity $\epsilon$, and frequency $\omega$.  The EIT problem is then to recover the coefficient $\gamma(z)$ for all $z\in\Omega$ from electrical measurements taken on the surface \trev{$\bndry$}.  These measurements take the form of Neumann-to-Dirichlet (ND) data $(g,f)$ where $g$ are the applied currents and $f$ the corresponding voltages: $\mathcal{R}_\gamma g=f$.  Knowledge of the ND map allows one to predict the resulting boundary voltages for any applied current.  For full boundary data, the Dirichlet-to-Neumann (DN) map is the inverse of the ND map: $\Lambda_\gamma = \left(\mathcal{R}_\gamma\right)^{-1}$.

D-bar methods are based on using a nonlinear Fourier transform of the DN data tailored to the EIT problem.  The admittivity is then recovered from the transformed data by solving a $\overline{\partial}_{k}$ (D-bar) equation in the transform variable $k$ for solutions known as  {\it Complex Geometrical Optics (CGO) solutions} to a related partial differential equation.  While there are several D-bar methods for EIT, they all have the same basic form:

\[ \begin{array}{c} \text{Current/Voltage}\\
\text{data}\end{array} \longrightarrow \begin{array}{c} \text{Scattering}\\ \text{data}\end{array}  \longrightarrow \begin{array}{c} \text{CGO}\\ \text{solutions}\end{array} \longrightarrow \begin{array}{c} \text{Admittivity.}\end{array}\]

Below we provide a brief overview of the  D-bar methods for EIT studied here.  We compare the results of the D-bar method for real-valued conductivities, which is based on transforming the admittivity equation \eqref{eq:admitt} to the Schr\"odinger equation in the global uniqueness proof [\cite{Nachman1996}], implemented and further developed in [\cite{Siltanen2000}, \cite{Isaacson2004}, \cite{Knudsen2009}],  with the results of two implementations of the D-bar method for complex conductivities [\cite{Hamilton2012}, \cite{Herrera2015}], in which the problem is transformed to a first order elliptic system [\cite{Francini2000}].  Computational details will be discussed in Section~\ref{sec:CompDetails}.  

Here and throughout we associate $\R^2$ and $\C$ via the mapping $z=(x,y)\mapsto x+iy$ and make use of the common $\dez$ and $\dbarz$ derivative operators 
\[\dez = \frac12\left(\partial_x -i\partial_y\right),\quad  \dbarz = \frac12\left(\partial_x +i\partial_y\right).\]

%--------------------------------------------------------------------
\subsection{The D-bar Method for Real-valued Conductivities}\label{sec:NachDbar}
%--------------------------------------------------------------------
Here we review the proof given in [\cite{Nachman1996}], first implemented in [\cite{Siltanen2000}] and established as a rigorous regularization strategy in [\cite{Knudsen2009}].  Since this method only applies to real-valued conductivities, consider $\gamma=\sigma$.  In this method, equation \eqref{eq:admitt} is first transformed into the Schr\"odinger equation $[-\Delta + q(z)]\tilde{u}(z) = 0$, for $z\in\Omega\subset\R^2$ via the change of variables $\tilde{u}=\sigma^{1/2}u$, $q(z)=\sigma^{-1/2}\Delta\sigma^{1/2}$.  Assume $\sigma$ is the constant $1$ in a neighborhood of the boundary, and extend  the conductivity  from $\Omega$ to the full plane as $\sigma(z)\equiv 1$ for $z\in\R^2\setminus \Omega$.  Note that the the assumption that $\sigma$ is constant near the boundary can be dropped, and the DN map extended through analytic continuation, so the method then applies on a slightly larger domain.  See [\cite{SiltanenTamminen2012}, \cite{Nachman1996}] for further details on this approach.

Nachman proved the existence of \trev{unique} CGO solutions $\psi_s(z,k)$ to the Schr\"odinger equation
\begin{equation}\label{eq:schro}
[-\Delta + q(z)]\psi_s(z,k) = 0,\qquad z\in\C,\;\;k\in\C\setminus0,
\end{equation}
where \trev{$e^{-ikz}\psi_s(z.k)-1 \in W^{1,p}(\R^2)$, $p>2$, i.e. $\psi_s\sim e^{ikz}$ for large $|k|$ or $|z|$.}  Here $k$ is a complex number, and $kz$ is complex multiplication.  The CGO solutions $\mu(z,k)\equiv e^{-ikz}\psi_s(z,k)$ satisfy a D-bar equation in the transform variable $k$:
\begin{equation}\label{eq:Nach_dbark}
\dbar_k\mu(z,k)=\frac{1}{4\pi\bar{k}}\T(k)e(z,-k)\overline{\mu(z,k)},
\end{equation}
where $e(z,k)\equiv e^{i(kz+\bar{k}\bar{z})}=e^{2i\Re(zk)}$ and $\T(k)$ denotes the scattering data, which is related to the DN map via
\begin{equation}\label{eq:Nach_scatBndry}
\T(k)=\int_{\bndry} e^{i\bar{k}\bar{z}}(\Lambda_\sigma-\Lambda_1)\psi_s(z,k)\;dz,
\end{equation}
where $\Lambda_1$ represents the DN map corresponding to a constant conductivity of 1.   The conductivity $\sigma(z)$ can be recovered from the CGOs $\mu(z,k)$ by:
\begin{equation}\label{eq:Nachcond}
\sigma(z)=\left(\lim_{|k|\to0}\mu(z,k)\right)^{1/2}.
\end{equation}

 A {\it Born approximation} $\texp$ to the fully nonlinear scattering data $\T(k)$ can be made by replacing $\psi_s$ in \eqref{eq:Nach_scatBndry} by its asymptotic condition $\psi_s\sim e^{ikz}$:
\begin{equation}\label{eq:Nach_texp}
\texp(k)=\int_{\bndry} e^{i\bar{k}\bar{z}}(\Lambda_\sigma-\Lambda_1)e^{ikz}\;dz.
\end{equation}
This approximation circumvents computing the trace of $\psi_s(z,k)$ on $\partial\Omega$ which requires using another boundary integral equation.  \tRev{Approximating $\psi$ by its asymptotic behavior is equivalent to making a leading order approximation from the asymptotic series for $\psi$.}  This assumption is most accurate in the case of low contrast in conductivity/permittivity.  

The reconstruction process is then:
\begin{enumerate}[label=\Roman*.]
\item Given $[\Lambda_\sigma,\Lambda_1]$, evaluate \eqref{eq:Nach_texp} for $\texp(k)$,  $0<|k|\leq R$ for some chosen radius $R$.
\item Solve the $\dbark$ equation \eqref{eq:Nach_dbark} for each $z\in\Omega$ and recover the conductivity $\sigexp(z)=\left(\muexp(z,0)\right)^{1/2}$.
\end{enumerate}

The method above produces {\it absolute}, also called {\it static}, EIT images.  To obtain time-difference EIT images relative to a reference data set $\Lambdaref$, we use the {\it differencing scattering transform} $\tdiff$ introduced in [\cite{Isaacson2006}] 
\begin{equation}\label{eq:Nach_tdiff}
\tdiff(k):=\texp^{\mbox{\tiny ,$\gamma$}}(k) - \texp^{\mbox{\tiny ,$\gamma_{\mbox{\tiny ref}}$}}(k)=\int_{\bndry} e^{i\bar{k}\bar{z}}(\Lambda_\sigma-\Lambdaref)e^{ikz}\;dz,
\end{equation}
in equation \eqref{eq:Nach_dbark} and compute $\sigmadiff\equiv \left[\mudiff(z,0)\right]^2-1$, since $\sigma=1$ near the boundary.  If $\sigma=\sigma_0\neq 1$ near $\bndry$, the problem can be scaled as in [\cite{Isaacson2004}].

%--------------------------------------------------------------------
\subsection{The D-bar Method for Complex-Valued Conductivities}\label{sec:FranDbar}
%--------------------------------------------------------------------
The methods in this section can reconstruct real or complex-valued admittivities $\gamma(z)=\sigma(z)+i\omega\epsilon(z)$.  We briefly review the methods developed in [\cite{Francini2000, Hamilton2012, Herrera2015}].  

First, equation \eqref{eq:admitt} is transformed to a first-order elliptic  system of equations
\begin{equation}\label{eq:DQ}
\left[-D+ Q(z)\right]\Psi(z,k)=0, \quad z,k\in\C,
\end{equation} 
where 
\[D=\left[\begin{array}{cc}
\bar\partial_z & 0\\
0 &\partial_z
\end{array}\right], \quad 
Q(z)=\left[\begin{array}{cc}
0 & -\frac{1}{2} \partial_z \log \gamma\\
-\frac{1}{2} \bar\partial_z \log \gamma & 0\\
\end{array}\right], \quad \text{and } \Psi(z,k)=\gamma^{1/2}(z)\left[\begin{array}{cc}
\dez u_1& \dez u_2\\
\dbarz u_1 & \dbarz u_2\\
\end{array}\right],
\]
where $u_1\sim \frac{e^{ikz}}{ik}$ and $u_2\sim \frac{e^{-ik\bar{z}}}{-ik}$ are CGO solutions to the admittivity equation \eqref{eq:admitt}, whose existence was established in [\cite{Hamilton2012}].

 Functions  defined by $M(z,k) \equiv \Psi(z,k)\left[\begin{array}{cc}
e^{-ikz} &0\\
0 & e^{ik\bar{z}}\\
\end{array}\right]$ are CGO solutions to a system of D-bar equations in the transform variable $k$
%\vspace{-1em}
\begin{equation}\label{eq:Fran_dbark}
\overline{\partial}_k M(z,k)=M(z,\bar{k})\left[\begin{array}{cc}
e(z,\bar{k}) & 0 \\
0 & e(z,-k)
\end{array}\right]S(k),
\end{equation}
where $S(k)$ is now a matrix of scattering data defined by
\begin{equation}\label{eq:Fran_scatDom}
S(k)=\left[\begin{array}{cc}
0 & \frac{i}{\pi}\int_{\R^2} e(z,-\bar{k})Q_{12}(z)M_{22}(z,k)dz\\
-\frac{i}{\pi}\int_{\R^2} e(z,k)Q_{21}(z)M_{11}(z,k)dz & 0\\
\end{array}\right].
\end{equation}
Since $\gamma\equiv 1$ in $\R^2\setminus \Omega$, the matrix potential $Q$ has compact support, and integration by parts in \eqref{eq:Fran_scatDom} results in 
\begin{eqnarray}
S_{21}(k) &=& -\frac{i}{2\pi}\int_{\bndry}e^{i\bar{k}\bar{z}}\Psi_{21}(z,k)\overline{\nu(z)}\;ds(z)\label{eq:Fran_PsiScat21} \\
S_{12}(k) &=& \frac{i}{2\pi}\int_{\bndry}e^{-i\bar{k}{z}}\Psi_{12}(z,k){\nu(z)}\;ds(z),   \label{eq:Fran_PsiScat21}
\end{eqnarray}
where $\nu$ denotes the unit outward normal vector to $\bndry$. 

The admittivity $\gamma(z)$ can be recovered from the CGO solutions at $k=0$ by first reconstructing $Q$
\begin{equation}\label{eq:MtoQ}
Q_{12}(z)=\frac{\dez\left[M_{11}(z,0)+M_{12}(z,0)\right]}{M_{22}(z,0)+M_{21}(z,0)},\qquad \text{and}\qquad Q_{21}(z)=\frac{\dbarz\left[M_{22}(z,0)+M_{21}(z,0)\right]}{M_{11}(z,0)+M_{12}(z,0)},
\end{equation}
and then undoing the change of variables by computing
\begin{equation}\label{eq:QtoGamma}
\gamma(z)=\exp\left\{-\frac{2}{\pi \bar{z}}\ast Q_{12}(z)\right\} \quad \text{or equivalently }\quad \gamma(z)=\exp\left\{-\frac{2}{\pi {z}}\ast Q_{21}(z)\right\},
\end{equation}
where $\ast$ denotes convolution in $z$ over $\R^2$.  

The two formulations of this matrix-based D-bar method differ in their connections of the scattering data to the  DN map $\Lambda_\gamma$.  
In the first approach [\cite{Hamilton2012}], traces of the CGO solutions $\Psi_{12}(z,k)$ and $\Psi_{21}(z,k)$ are computed for $z\in\bndry$ via
\begin{equation}\label{eq:Fran_Psibndry}
\begin{array}{lcl}
\Psi_{12}(z,k)&=& \int_{\bndry} \frac{e^{i\bar{k}(k-\zeta)}}{4\pi(z-\zeta)}\left(\Lambda_\gamma-\Lambda_1\right)u_2(\zeta,k)ds(\zeta)\\
\Psi_{21}(z,k)&=& \int_{\bndry} \overline{\left[\frac{e^{i{k}(k-\zeta)}}{4\pi(z-\zeta)}\right]}\left(\Lambda_\gamma-\Lambda_1\right)u_1(\zeta,k)ds(\zeta).
\end{array}
\end{equation}
 Replacing $u_1$ and $u_2$ by their asymptotic behaviors ($u_1\sim \frac{e^{ikz}}{ik}$ and $u_2\sim \frac{e^{-ik\bar{z}}}{-ik}$) yields the {\it Born approximations}  as in [\cite{Hamilton2017_PhysMeas}]:
\begin{equation}\label{eq:Fran_PsiExp}
\begin{array}{lcl}
\Psiexponetwo(z,k)&=& \int_{\bndry} \frac{e^{i\bar{k}(k-\zeta)}}{4\pi(z-\zeta)}\left(\Lambda_\gamma-\Lambda_1\right)\left(\frac{e^{-ik\bar{\zeta}}}{-ik}\right)ds(\zeta)\\
\Psiexptwoone(z,k)&=& \int_{\bndry} \overline{\left[\frac{e^{i{k}(k-\zeta)}}{4\pi(z-\zeta)}\right]}\left(\Lambda_\gamma-\Lambda_1\right)\left(\frac{e^{ik\zeta}}{ik}\right)ds(\zeta),
\end{array}
\end{equation}
and corresponding scattering data
\begin{equation}\label{eq:Fran_scatPsiExp}
\begin{array}{lcl}
S^{\mbox{\tiny{$\Psi$exp}}}_{12}(k)&=&\frac{i}{2\pi}\int_{\bndry}e^{-i\bar{k}{z}}\Psiexponetwo(z,k){\nu(z)}ds(z)\\
S^{\mbox{\tiny{$\Psi$exp}}}_{21}(k)&=&\frac{i}{2\pi}\int_{\bndry}e^{i\bar{k}\bar{z}}\Psiexptwoone(z,k)\bar{\nu}(z)ds(z).
\end{array}
\end{equation}

The  reconstruction process is then:

\vspace{0.5em}
\noindent
\underline{\bf Approach $1$:}
\begin{enumerate}[label=\Roman*.]
\item Given $[\Lambda_\gamma,\Lambda_1]$, evaluate the approximate CGOs $\left\{\Psiexponetwo,\; \Psiexptwoone \right\}$ from \eqref{eq:Fran_PsiExp} and compute the corresponding approximate scattering data for $\left\{S^{\mbox{\tiny{$\Psi$exp}}}_{12}(k),\;S^{\mbox{\tiny{$\Psi$exp}}}_{21}(k)\right\}$ from \eqref{eq:Fran_scatPsiExp}, for  $0<|k|\leq R$ for some chosen radius $R$.
\item Solve the $\dbark$ equation \eqref{eq:Fran_dbark}, using $S^{\mbox{\tiny{$\Psi$exp}}}(k)$, for each $z\in\Omega$.
\item Recover the admittivity $\gamma^{\mbox{\tiny{$\Psi$exp}}}$ \eqref{eq:QtoGamma} using the matrix potential $Q^{\mbox{\tiny{$\Psi$exp}}}$ \eqref{eq:MtoQ}.
\end{enumerate}

\vspace{1em}
The second approach is to re-write the equations for the scattering data directly in terms of the DN map $\Lambda_{\gamma}$ [\cite{Herrera2015}]:
\begin{equation}\label{eq:Fran_scatNata}
\begin{array}{lcl}
S_{12}(k) & =& \frac{i}{4\pi}\int_{\partial\Omega} e^{-i\bar{k}z}\left(\Lambda_{\gamma}+i\partial_{\tau}\right)u_2(z,k)\;ds(z)\\
S_{21}(k) & =& -\frac{i}{4\pi}\int_{\partial\Omega} e^{i\bar{k}\bar{z}}\left(\Lambda_{\gamma}-i\partial_{\tau}\right)u_1(z,k)\;ds(z),
\end{array}
\end{equation}
where $\partial\tau$ denotes the tangential derivative operator: $\partial_\tau f(z) = \nabla f(z) \cdot \tau$.  Again, using the asymptotic behavior of the CGO solutions $u_1$ and $u_2$ gives a  {\it Born approximation} to the scattering data \eqref{eq:Fran_scatNata}
\begin{equation}\label{eq:Fran_scatExp}
\begin{array}{lcl}
\Sexponetwo(k) & =& \frac{i}{4\pi}\int_{\partial\Omega} e^{-i\bar{k}z}\left(\Lambda_{\gamma}+i\partial_{\tau}\right)\left(\frac{e^{-ik\bar{z}}}{-ik}\right)ds(z)\\
\Sexptwoone(k) & =& -\frac{i}{4\pi}\int_{\partial\Omega} e^{i\bar{k}\bar{z}}\left(\Lambda_{\gamma}-i\partial_{\tau}\right)\left(\frac{e^{ikz}}{ik}\right)ds(z),
\end{array}
\end{equation}

The reconstruction process is then:

\vspace{0.5em}\noindent
 \underline{\bf Approach $2$:}
\begin{enumerate}[label=\Roman*.]
\item Given $\Lambda_\gamma$, compute the {\it Born} scattering data for $\left\{\Sexponetwo(k),\;\Sexptwoone(k)\right\}$ from \eqref{eq:Fran_scatExp},  for $0<|k|\leq R$ for some chosen radius $R$.
\item Solve the $\dbark$ equation \eqref{eq:Fran_dbark}, using $\Sexp(k)$, for each $z\in\Omega$.
\item Recover the admittivity $\gamma^{\mbox{\tiny{exp}}}$ \eqref{eq:QtoGamma} using the matrix potential $Q^{\mbox{\tiny{exp}}}$ \eqref{eq:MtoQ}.
\end{enumerate}

\vspace{1em}
Note that {\it `Approach~2'} does not require the DN map $\Lambda_1$ and thus no simulation of data is needed to obtain absolute images.  Aside from Calder\'on's original method [\cite{Calderon1980}, \cite{Bikowski2008}, \cite{Muller2017}], this is the only method for absolute imaging of EIT that does not need any simulated data.  This has the key benefits of not needing to know the locations of electrodes with high precision, contact impedances at the electrodes, or even near perfect knowledge of the boundary shape.  As we demonstrate below in Section~\ref{sec:results}, the method is remarkably stable against perturbations in electrodes locations as well as incorrect information about boundary shape.  In fact, all of the D-bar methods presented here perform quite well (see Figures~\ref{fig:ACT3_HnL_Nach}, \ref{fig:ACT3_HnL_FranSarah} and \ref{fig:ACT3_HnL_FranNata})\tRev{.   We remark that the action of the DN map $\Lambda_1$ on the complex exponential functions ($e^{ikz}$ for $\texp$, and $\frac{e^{ikz}}{ik}$ and $\frac{e^{-ik\bar{z}}}{-ik}$ for {\it `Approach 1'}) can be approximated by using the definition of the DN map: $\Lambda_\gamma u=\gamma\nabla u\cdot \nu$.  In [\cite{Isaacson2004}], this was done by analytically computing the action of $\Lambda_1$ on the basis of applied trigonometric current patterns.  Alternatively, one could use a continuum type approximation directly computing, e.g., $\Lambda_1 e^{ikz}=1\nabla e^{ikz}\cdot\nu = ik\nu e^{ikz}$.  However,} due to the flexibility of not needing any simulated data, we focus this work on {\it `Approach 2'}.
To obtain time-difference EIT images, define the {\it differencing scattering transform}:
\begin{itemize}
\item \underline{\bf Approach~1 - Difference Imaging:} [\cite{Hamilton2017_PhysMeas}]
\begin{equation}\label{eq:Fran_scatPsiExpDiff}
\begin{array}{lcl}
S^{\mbox{\tiny{$\Psi$diff}}}_{12}(k)&=&\frac{i}{2\pi}\int_{\bndry}e^{-i\bar{k}{z}}\Psidiffonetwo(z,k){\nu(z)}ds(z)\\
S^{\mbox{\tiny{$\Psi$diff}}}_{21}(k)&=&\frac{i}{2\pi}\int_{\bndry}e^{i\bar{k}\bar{z}}\Psidifftwoone(z,k)\bar{\nu(z)}ds(z),
\end{array}
\end{equation}
where
\begin{equation}\label{eq:Fran_PsiDiff}
\begin{array}{lcl}
\Psidiffonetwo(z,k)&=& \int_{\bndry} \frac{e^{i(\bar{k}(k-\zeta))}}{4\pi(z-\zeta)}\left(\Lambda_\gamma-\Lambdaref\right)\left(\frac{e^{-ik\bar{\zeta}}}{-ik}\right)ds(\zeta)\\
\Psidifftwoone(z,k)&=& \int_{\bndry} \overline{\left[\frac{e^{i({k}(k-\zeta))}}{4\pi(z-\zeta)}\right]}\left(\Lambda_\gamma-\Lambdaref\right)\left(\frac{e^{ik\zeta}}{ik}\right)ds(\zeta).
\end{array}
\end{equation}

\vspace{1em}
\item \underline{\bf Approach~2 - Difference Imaging:}  [\cite{Herrera2015}]
\begin{equation}\label{eq:Fran_scatdiff}
\begin{array}{lcl}
\Sdiffonetwo(k) & =& \frac{i}{4\pi}\int_{\partial\Omega} e^{-i\bar{k}z}\left(\Lambda_{\gamma}-\Lambdaref\right)\left(\frac{e^{-ik\bar{z}}}{-ik}\right)ds(z)\\
\Sdifftwoone(k) & =& -\frac{i}{4\pi}\int_{\partial\Omega} e^{i\bar{k}\bar{z}}\left(\Lambda_{\gamma}-\Lambdaref\right)\left(\frac{e^{ikz}}{ik}\right)ds(z),
\end{array}
\end{equation}

\end{itemize}
then, solve the $\dbark$ equation  \eqref{eq:Fran_dbark} using the differencing scattering data and recover $\gammadiff$ by subtracting $1$ from \eqref{eq:QtoGamma}, as 1 is the constant admittivity near the boundary $\bndry$.  As before, if $\gamma=\gamma_0\neq1$ near $\bndry$, the problem can be scaled.

Note that if the admittivity is purely real-valued, then simplifications to the first-order system method exist [\cite{Brown1997}, \cite{Knudsen2003}, \cite{Knudsen2002}].

%--------------------------------------------------------
\subsection{Computational Details}\label{sec:CompDetails}
%--------------------------------------------------------
Here we present a self-contained summary of the implementation of {\it `Approach~2'} for complex-valued conductivities.  For numerical implementations of the D-bar method for real-valued conductivities see [\cite{Mueller2012}, \cite{DeAngelo2010}] and for {\it `Approach~1'} for complex-valued conductivities see [\cite{Hamilton2012, Hamilton2013}].

We present the numerical solution method based on trigonometric current patterns, however any set of $L-1$ linearly independent measurements (for a system with $L$ electrodes) can be transformed via a change of basis to synthesize the measurements that would have occurred if the trigonometric current patterns presented here had been applied.  

To form the discrete approximation $\mathbf{R}_\gamma$ to the ND map $\mathcal{R}_\gamma$ we use an inner product of the measured currents and voltages on the electrodes as in [\cite{Isaacson2004}].  Let $\Phi^j_\ell$ denote the current on the $\ell$-th electrode arising from the $j$-th trigonometric current pattern:
\begin{equation}\label{eq:trigCPs}
\Phi^j_\ell=\Phi(\ell,j):=\begin{cases}
A\cos(j\theta_\ell), & 1\leq j\leq \frac{L}{2},\;\;1\leq \ell\leq L\\
A\sin((\frac{L}{2}-j)\theta_\ell) & \frac{L}{2}+1\leq j\leq L-1,\;\;1\leq \ell\leq L,
\end{cases}
\end{equation}
where $A$ is the amplitude of the applied current, and $V^j_\ell$ denotes the corresponding voltage measurement.  We form $\mathbf{R}_\gamma$ using the normalized currents $\phi^j=\frac{\Phi^j}{\|\Phi^j\|_2}$ and voltages $v^j_\ell$ satisfying $\sum_{\ell=1}^L v^j_\ell=0$ and $v^j=\frac{v^j}{\|\Phi^j\|_2}$:
\begin{equation}\label{eq:NDcomp}
\mathbf{R}_\gamma(m,n):=\sum_{\ell=1}^L \frac{\phi^m_\ell v^n_\ell}{|e_\ell|},\quad 1\leq m,n\leq L-1,\;\; 1\leq \ell\leq L,
\end{equation}
where $|e_\ell |$ is the area of the $\ell$-th electrode.  The discrete approximation $\mathbf{L}_\gamma$ to the DN map $\Lambda_\gamma$ is then $\mathbf{L}_\gamma = \left(\mathbf{R}_\gamma\right)^{-1}$.  While the method as described assumes \trev{$\gamma=\gamma_0 =1$} on $\bndry$ and the domain has maximal radius 1,  the DN matrix $\mathbf{L}_\gamma$ can be scaled by $1/{\gamma_0}$ and to a max radius 1 by multiplying the DN matrix $\mathbf{L}_\gamma$ by $r$, the radius of the smallest circle containing the domain, as in [\cite{Isaacson2004}].  Note that an estimate to $\gamma_0$ is the best constant-admittivity approximation to measured data, as in [\cite{Cheney1990a}].  The resulting conductivity or admittivity at the end of the algorithm is then rescaled by multiplying by $\sigma_0$ or $\gamma_0$ respectively.

The evaluation of the scattering data \eqref{eq:Fran_scatExp} requires computing the action of the DN map on the asymptotic behaviors of $u_1$ and $u_2$, namely $\mathbf{L}_\gamma \left(\frac{e^{ikz}}{ik}\right)$ and $\mathbf{L}_\gamma \left(\frac{e^{-ik\bar{z}}}{-ik}\right)$.  To accomplish this, we follow the method outlined in [\cite{DeAngelo2010}] and expand the functions $\frac{e^{ikz}}{ik}$ and $\frac{e^{-ik\bar{z}}}{-ik}$ in the basis of normalized current patterns $\{\phi^j\}_{j=1}^{L-1}$ as
%\begin{eqnarray*}
%\frac{e^{ikz_\ell}}{ik}&\approx& \sum_{\ell}^L a_j(k)\phi^j_\ell =:\; e^{\mbox{\tiny u1}}_\ell(k)\\
%\frac{e^{-ik\overline{z_\ell}}}{-ik}&\approx& \sum_{\ell}^L b_j(k)\phi^j_\ell =:\; e^{\mbox{\tiny u2}}_\ell(k),
%\end{eqnarray*}
\[
\frac{e^{ikz_\ell}}{ik}\approx \sum_{j=1}^{L-1} a_j(k)\phi^j_\ell =:\; e^{\mbox{\tiny u1}}_\ell(k), \qquad\text{and} \qquad
\frac{e^{-ik\overline{z_\ell}}}{-ik}\approx \sum_{j=1}^{L-1} b_j(k)\phi^j_\ell =:\; e^{\mbox{\tiny u2}}_\ell(k),
\]
where $z_\ell$ denotes the center of the $\ell$-th electrode.  Choosing a scattering radius $R>0$,  the integrals in \eqref{eq:Fran_scatExp} are computed using a Simpson's rule approximation
\begin{equation}\label{eq:Fran_scatExp_numerics}
\begin{array}{lcl}
\Sexponetwo(k) & \approx & \begin{cases}
\frac{i}{4\pi}\frac{P}{L}\left[e^{-i\bar{k}\overline{\mathbf{z}}}\right]^T\left[\phi\mathbf{L}_\gamma\mathbf{e}^{\mbox{\tiny u2}}(k) + \mathbf{d}_{12}(k)\right] & 0<|k|\leq R\\
0 & |k|>R
\end{cases}\\
&&\\
\Sexptwoone(k) & \approx & \begin{cases}
-\frac{i}{4\pi}\frac{P}{L}\left[e^{i\bar{k}\mathbf{z}}\right]^T\left[\phi\mathbf{L}_\gamma\mathbf{e}^{\mbox{\tiny u1}}(k) + \mathbf{d}_{21}(k)\right]& 0<|k|\leq R\\
0 & |k|>R,
\end{cases}\\
\end{array}
\end{equation}
where $P$ is the perimeter of the boundary $\bndry$, $T$ the matrix transpose, and $\mathbf{z}$ is the vector of the centers of the electrodes $z_\ell$.  The $\mathbf{d}_{12}$ and $\mathbf{d}_{21}$ terms represent the $32\by1$ vectors resulting from the action of the tangential derivative map $\partial_\tau$ acting on $\frac{e^{-ik\overline{\mathbf{z}}}}{-ik}$ and $\frac{e^{ik\mathbf{z}}}{ik}$, respectively.  These can be computed analytically as
\[\partial_\tau\left(\frac{e^{-ik\overline{\mathbf{z}}}}{-ik}\right)=\nabla\left(\frac{e^{-ik\overline{\mathbf{z}}}}{-ik}\right)\cdot\tau=e^{-ik\bar{z}}(1,-i)\cdot(\tau_1,\tau_2) = e^{-ik\bar{z}}\overline{\tau(z)}=-i\overline{\nu(z)}e^{-ik\bar{z}},\]
and $\partial_\tau\left(\frac{e^ik\mathbf{z}}{ik}\right)=i\nu(z)e^{ikz}$.  The unit normal vector $\nu(z)$ can be computed via numerical forward-differences using a parameterization $r(\theta)$ of the boundary.  The unit tangent vector to the boundary can then be computed as
\[\tau(\theta) = \frac{r(\theta+\varepsilon)-r(\theta)}{|r(\theta+\varepsilon)-r(\theta)|}=\tau_1(\theta)+i\tau_2(\theta),\]
for some small $\varepsilon>0$ and then the unit normal $\nu = \tau_2-i\tau_1$.  Note that this step is trivially parallelizable in $k$.  The scattering data for $k=0$ can be computed by interpolation.

Next, the $\dbark$ system \eqref{eq:Fran_dbark} decouples into two systems of two equations, each which can be solved using the integral form
\begin{equation}\label{eq:Fran_dbark_convs}
\begin{array}{c}
\left\{\begin{array}{lcl}
1 & =& M_{11}(z,k) - \frac{1}{\pi k}\ast\left[M_{12}(z,\bar{k})e(z,-k)\Sexptwoone(k)\right]\\
0 & =& M_{12}(z,k) - \frac{1}{\pi k}\ast\left[M_{11}(z,\bar{k})e(z,\bar{k})\Sexponetwo(k)\right]
\end{array}\right.\\
\\
\left\{\begin{array}{lcl}
1 & =& M_{22}(z,k) - \frac{1}{\pi k}\ast\left[M_{21}(z,\bar{k})e(z,\bar{k})\Sexponetwo(k)\right]\\
0 & =& M_{21}(z,k) - \frac{1}{\pi k}\ast\left[M_{22}(z,\bar{k})e(z,-k)\Sexptwoone(k)\right],
\end{array}\right.\\
\end{array}
\end{equation}
where $\ast$ denotes convolution in $k$ over the disc of radius $R$, since $\Sexp$ has compact support in $|k|\leq R$, and we have used the fundamental solution $\frac{1}{\pi k}$ for the $\dbark$ operator as well as the asymptotic condition $M(z,k)\sim \left[\begin{array}{cc} 1 & 0\\ 0 & 1\\ \end{array}\right]$ for large $|k|$ or $|z|$.  The solution can be obtained using a modified version of Vainikko's method [\cite{Vainikko2000}] for solving integral equations with weakly singular kernels.  The convolutions can be implemented using two-dimensional  {\sc FFT}s on a uniform $k$-grid of size $(2^{N}+1)\by(2^{N}+1)$ with uniform step size $h_k$ as
\[\frac{1}{\pi k}\ast f(k)=h_k^2\;IFFT\left[FFT\left(\frac{1}{\pi k}\right)\;FFT\left(f(k)\right)\right].\]
The system in \eqref{eq:Fran_dbark_convs} can be solved using matrix-free {\sc GMRES} for each desired $z\in\Omega$.  Note that this step is trivially parallelizable in $z$, and that after each solution of \eqref{eq:Fran_dbark_convs} only the $k=0$ entry is required to recover the admittivity $\gamexp(z)$.

The potentials $\Qexponetwo$ and $\Qexptwoone$ can be computed from \eqref{eq:MtoQ}, with $\Mexp(z,0)$, using a numerical derivative such as centered finite-differences, and the admittivity obtained via convolutions and two-dimensional {\sc FFT}s
\begin{eqnarray*}
\gamexp(z) &=& \exp\left\{-2h_z^2\;IFFT\left[FFT\left(\frac{1}{\pi \bar{z}}\right)\;FFT(\Qexponetwo)\right]\right\}\\
&=& \exp\left\{-2h_z^2\;IFFT\left[FFT\left(\frac{1}{\pi {z}}\right)\;FFT(\Qexptwoone)\right]\right\}
\end{eqnarray*}
where $h_z$ denotes the uniform step size in the $z$ grid, and again we have used the fundamental solutions for $\dez$ and $\dbarz$, respectively.  If the DN matrix $\mathbf{L}_\gamma$ was scaled by $\gamma_0\neq 1$, then undo the scaling now by multiplying by $\gamma_0$.

%--------------------------------------------------------------------
\subsection{\trev{Evaluation methods}}\label{sec:eval_mtds}
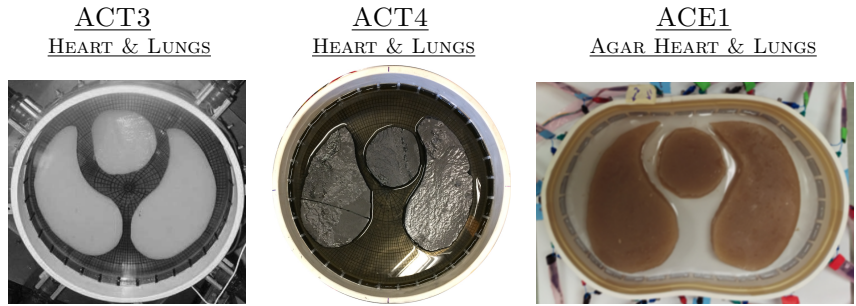
\begin{figure}[h!]
\begin{picture}(320,130)
\put(0,0){\includegraphics[width=90pt]{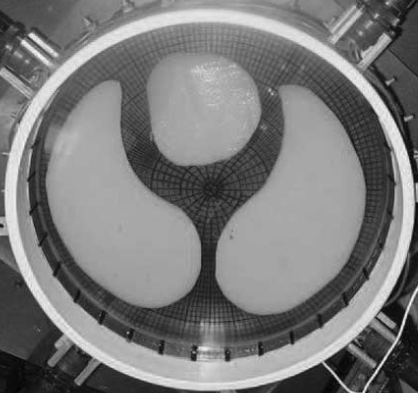}}
%\put(100,0){\includegraphics[width=90pt]{HLSA_plueral_saline_DeepDbar_DICOM_smaller.jpg}}
\put(100,0){\includegraphics[width=90pt]{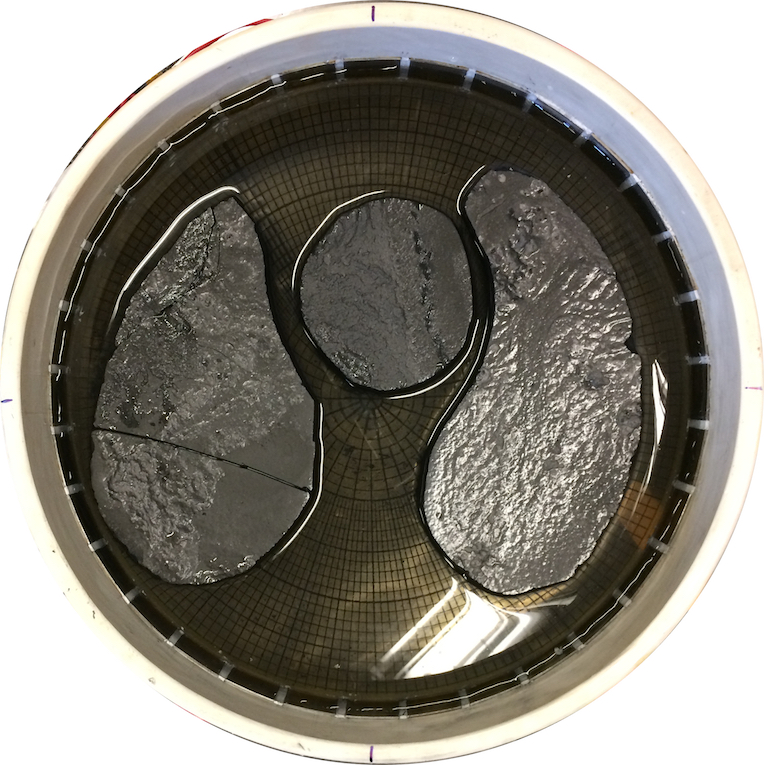}}

\put(200,0){\includegraphics[width=120pt]{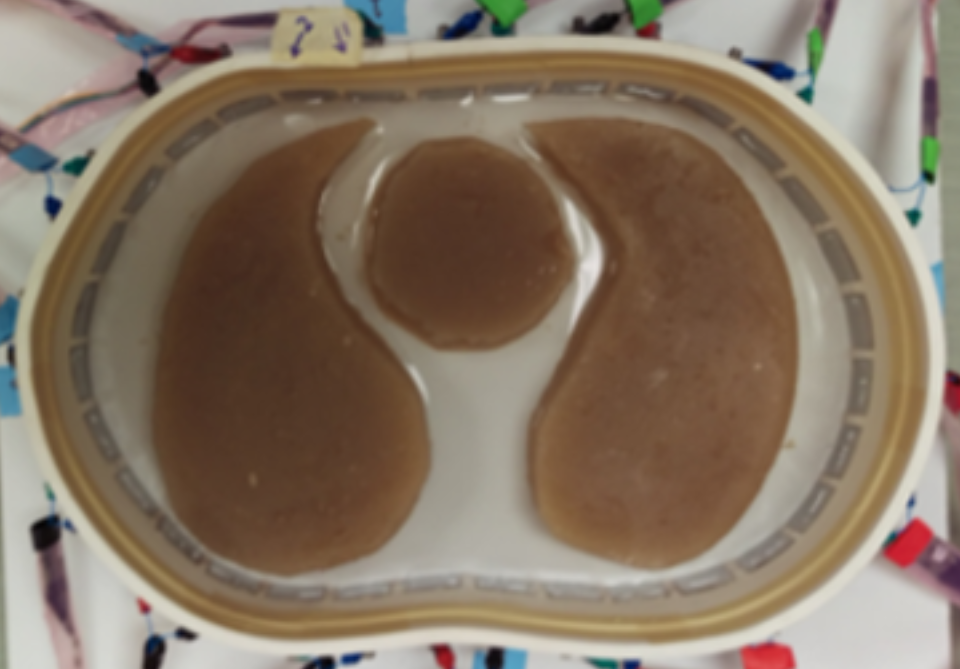}}
%\put(330,0){\includegraphics[width=120pt]{MelonTrimmedRot_smaller.png}}
%\put(300,0){\includegraphics[width=90pt]{HLSA_plastic_saline_DICOM.jpg}}
%\put(200,5){\includegraphics[width=70pt]{eH1Lpic_clipped.jpg}}

\put(25,105){\sc \underline{ACT3}}
\put(15,95){\sc \footnotesize \underline{Heart \& Lungs}}

\put(130,105){\sc \underline{ACT4}}
%\put(125,95){\sc \footnotesize \underline{HLSA PE}}
\put(115,95){\sc \footnotesize \underline{\trev{Heart \& Lungs}}}

\put(245,105){\sc \underline{ACE1}}
\put(220,95){\sc \footnotesize \underline{Agar Heart \& Lungs}}

%\put(375,105){\sc \underline{ACE1}}
%\put(325,95){\sc \footnotesize \underline{Cucumber Heart, Melon Lungs}}

%\put(330,105){\sc \underline{ACT 4}}
%\put(315,95){\sc \footnotesize \underline{HLSA-Plastic}}

\end{picture}
\caption{Experimental setups for test phantoms from the ACT3, ACT4, and ACE1 EIT systems.}
\label{fig:phantoms}
\end{figure}

%--------------------------------------------------------------------
The robustness of the D-bar methods are demonstrated on EIT data from the {\sc ACT3} and {\sc ACT4} EIT systems of Rensselaer Polytechnic Institute (RPI) [\cite{Cook1994,Liu2005}] and the {\sc ACE1} EIT system of Colorado State University (CSU) [\cite{MellenthinEMBS_2015}].  See Figure~\ref{fig:phantoms} for the experimental setups.

%--------------------------------------------------------
% Phantoms used (ACT3, ACT4, and ACE1)
%--------------------------------------------------------
%--------------------------------------------------------

The ACT3/ACT4 data were collected on a saline-filled tank of radius 15 cm with 32 electrodes of width 2.5 cm on the boundary.  The ACT3 data set is archival, with an agar heart (0.75~S/m) and two agar lungs (0.24~S/m) in a saline background (0.424~S/m) filled to a depth of 1.6 cm. Trigonometric current patterns of amplitude 0.2 mA were applied at a frequency of 28.8~kHz (see [\cite{Isaacson2004}]).  Results are presented in Section~\ref{sec:results_ACT3}.  The targets for the ACT4 data were made of agar with added graphite to simulate a chest phantom with a heart \trev{and two lungs} in a saline bath of 0.3~S/m and height 2.25~cm.   %Using the DICOM orientation, the lower portion of the right lung was removed and replaced with agar/graphite matching that of the heart.  
The admittivities of the targets, measured using the {sc SFP-7} bioimpedance meter sold by {\sc Impedimed}, \trev{were as follows: heart {$0.68 +0.05i$~S/m} and lungs {$0.057 + 0.011i$~S/m}.}  The ACT4 system applies voltages and measures currents rather than vice-versa.  In this experiment, trigonometric voltage patterns with maximum amplitude 0.5 V were applied at \trev{3} kHz.  ACT4 results are presented in Section~\ref{sec:results_ACT4}.

%The data from the ACE1 system were taken on a chest shaped tank of perimeter 1.026 m with 32 electrodes of width 2.54 cm with bipolar adjacent current patterns applied at 3.3 mA and 125 kHz.  Results from two experiments are studied here.  The first uses agar heart and lung targets (0.45~S/m and 0.09~S/m, respectively) placed in a saline bath of 0.2~S/m filled to a height of 2.04 cm.  The second uses a cucumber heart (0.03~S/m) and watermelon lungs (left: 0.032~S/m, right: 0.054~S/m) in a saline bath (0.187~S/m) filled to a height of 1.7~cm.  The susceptivity of the melon and cucumber are known to be positive but their specific values are unknown.   The conductivities of the cucumber and melon were measured using an Omega CDH221 conductivimeter.   Reconstructions from the ACE1 data are presented in Section~\ref{sec:results_ACE1}.

The data from the ACE1 system were taken on a chest shaped tank of perimeter 1.026 m with 32 electrodes of width 2.54 cm with bipolar adjacent current patterns applied at 3.3 mA and 125 kHz.  \trev{A}gar heart and lung targets (0.45~S/m and 0.09~S/m, respectively) \trev{were} placed in a saline bath of 0.2~S/m filled to a height of 2.04 cm.  Reconstructions from the ACE1 data are presented in Section~\ref{sec:results_ACE1}.

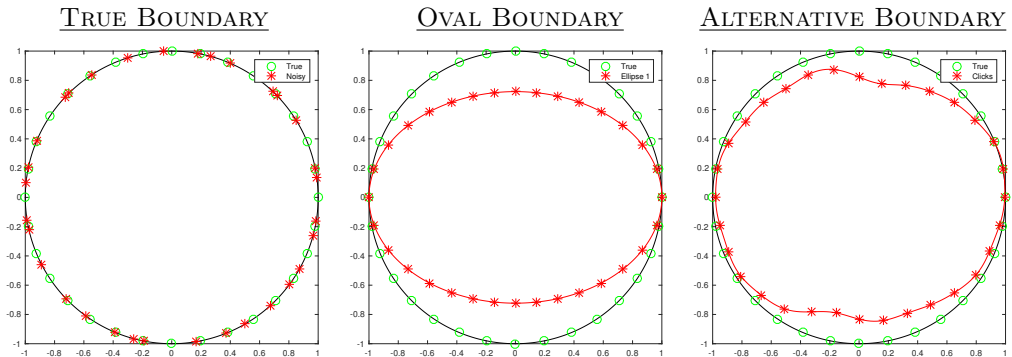
\begin{figure}[h!]
\begin{picture}(380,140)
\put(0,0){\includegraphics[width=120pt]{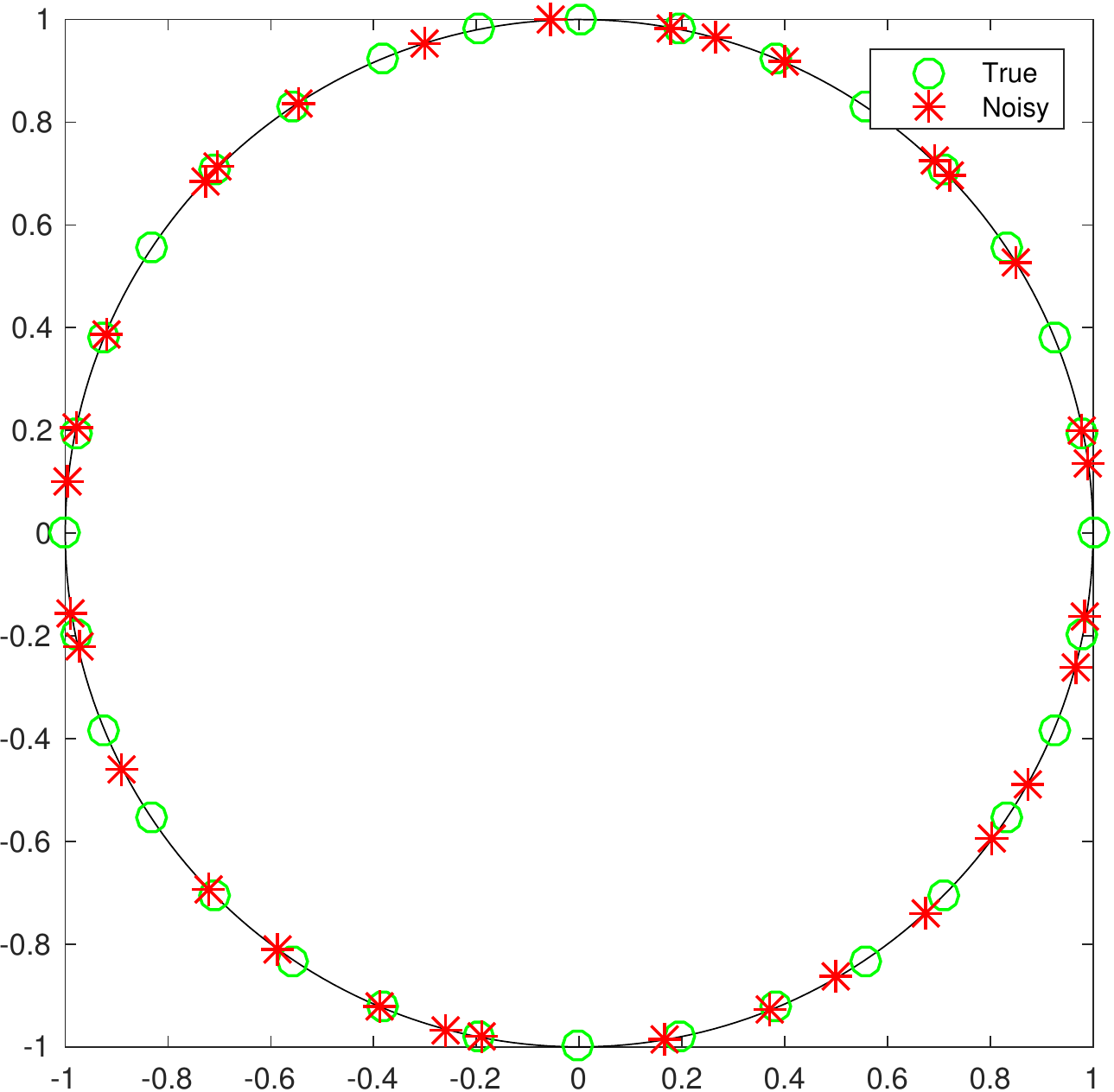}}
\put(130,0){\includegraphics[width=120pt]{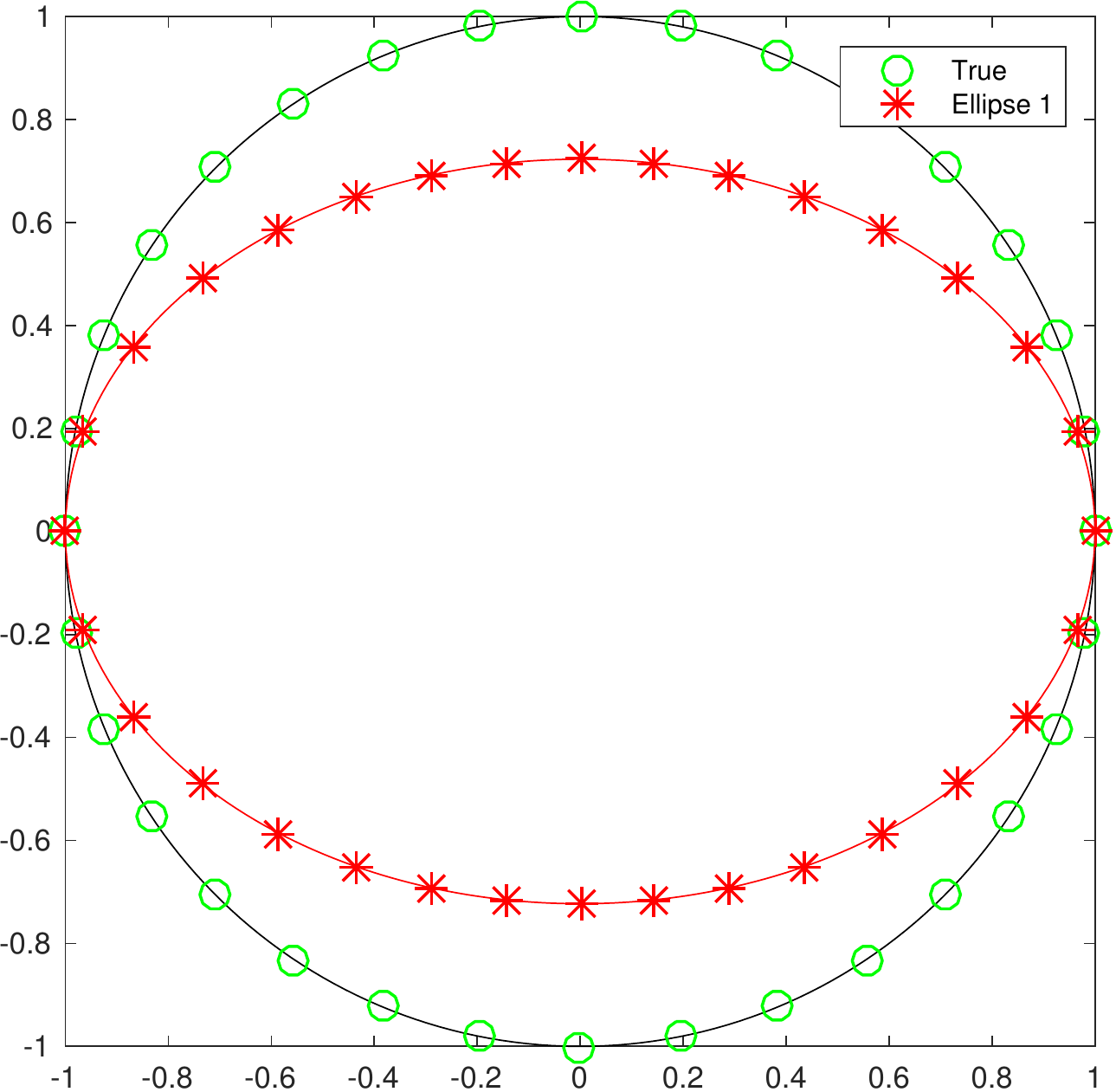}}
\put(260,0){\includegraphics[width=120pt]{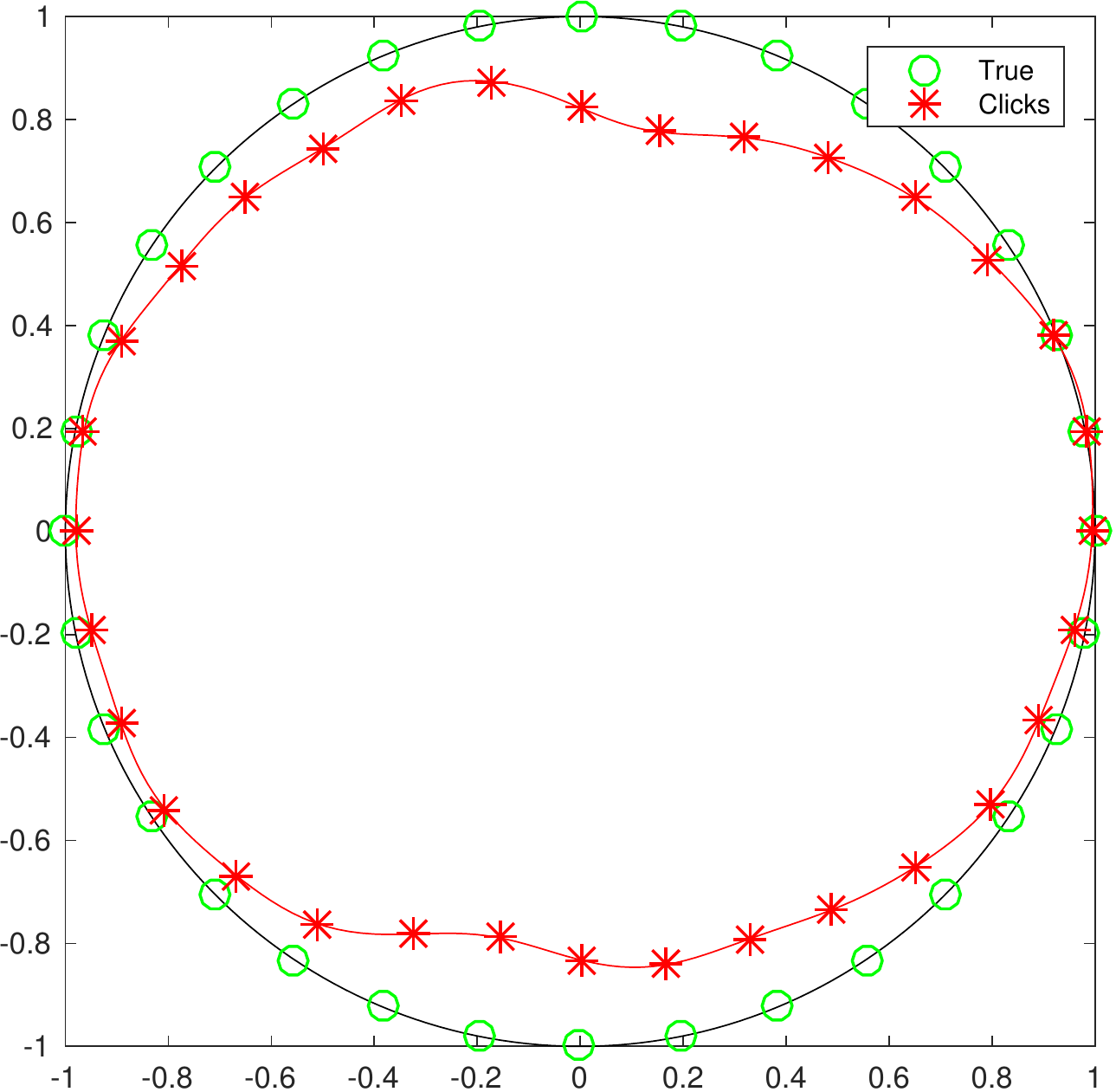}}

%\put(240,0){\includegraphics[width=110pt]{ACT4circ_Oval-eps-converted-to.pdf}}
%\put(240,0){\includegraphics[width=110pt]{ACT4circ_Oval_NoisyElecs-eps-converted-to.pdf}}

\put(20,125){\sc \underline{True Boundary}}
%\put(5,80){\it \footnotesize Heart \& Lungs}

\put(155,125){\sc\underline{Oval Boundary}}

\put(262,125){\sc\underline{Alternative Boundary}}
%\put(105,80){\it \footnotesize Heart \& Lungs}

%\put(260,115){\sc\underline{HLSA-Plastic}}
%\put(210,80){\it \footnotesize Heart \& Lung}
\end{picture}
\caption{Boundary shapes and electrode locations tested for the ACT3/ACT4 data.  Left: True boundary shape with incorrect electrode locations in red.  Middle: Oval (incorrect) boundary shape.  Right: An incorrect boundary created from a picture of the experiment, referred to as ``alternative boundary''. }
\label{fig:ACT4_circ_elecsNoisy}
\end{figure}

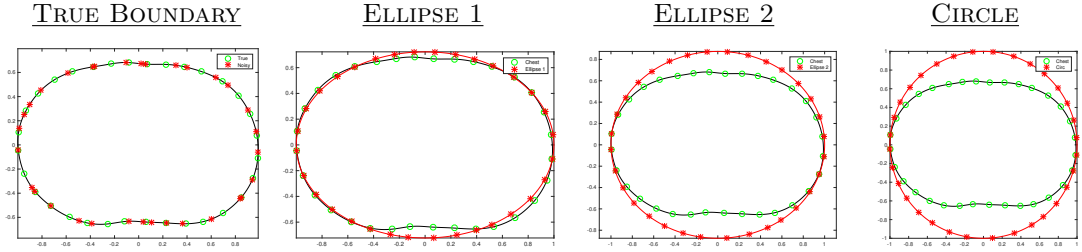
\begin{figure}[h!]
%%\begin{picture}(350,300)
%
%\put(0,155){\includegraphics[height=100pt]{ACE1chest_NoisyElecs-eps-converted-to.pdf}}
%\put(180,155){\includegraphics[height=100pt]{ACE1chest_Oval1-eps-converted-to.pdf}}
%
%\put(0,0){\includegraphics[height=120pt]{ACE1chest_Oval2-eps-converted-to.pdf}}
%\put(200,0){\includegraphics[height=120pt]{ACE1chest_Circ-eps-converted-to.pdf}}
%
%\put(40,285){\sc \underline{True Boundary}}
%\put(245,285){\sc\underline{Ellipse 2}}
%
%\put(60,130){\sc\underline{Ellipse 2}}
%
%\put(248,130){\sc\underline{Circle}}

%\end{picture}
\begin{picture}(405,110)

\put(0,0){\includegraphics[height=75pt]{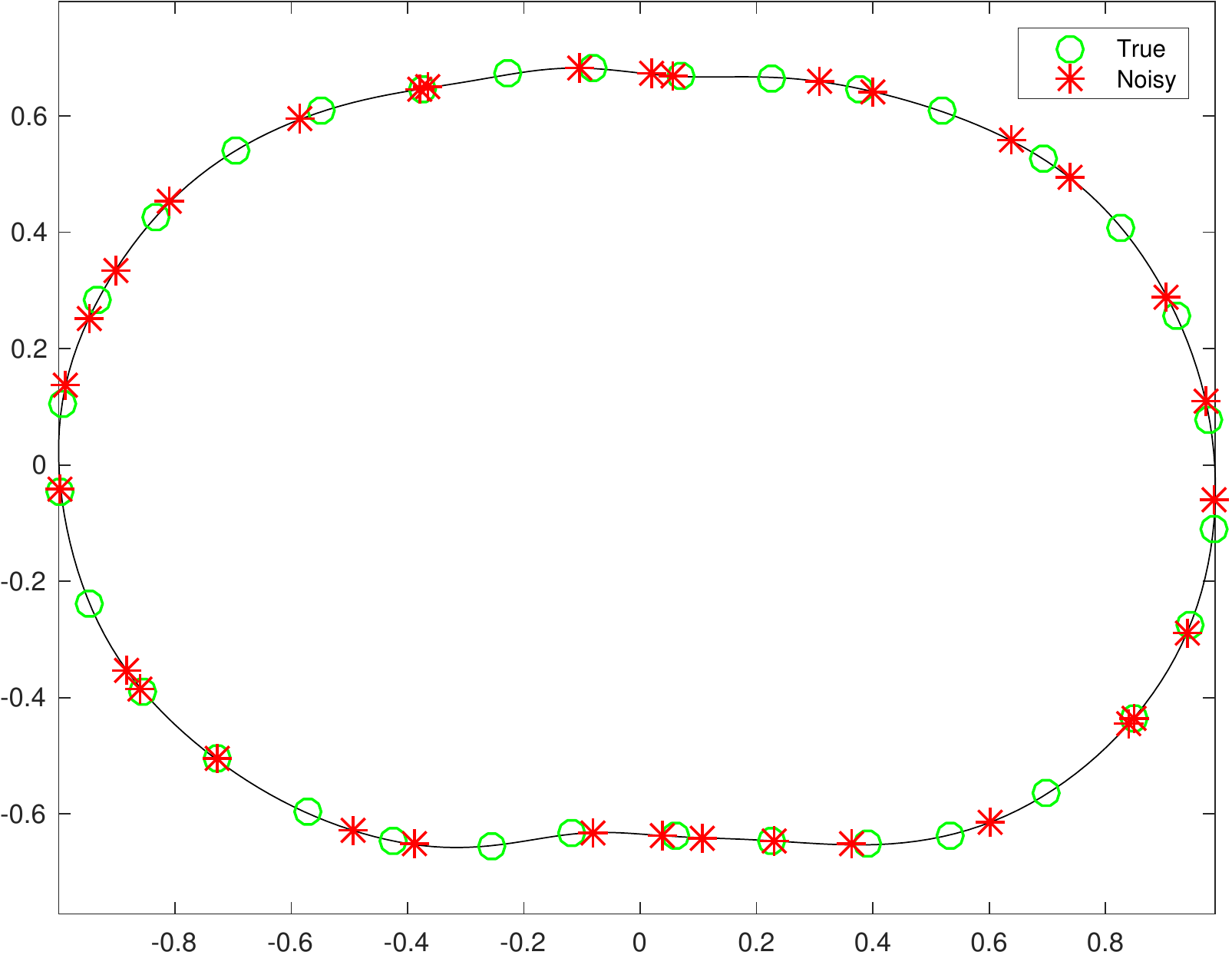}}
\put(105,0){\includegraphics[height=75pt]{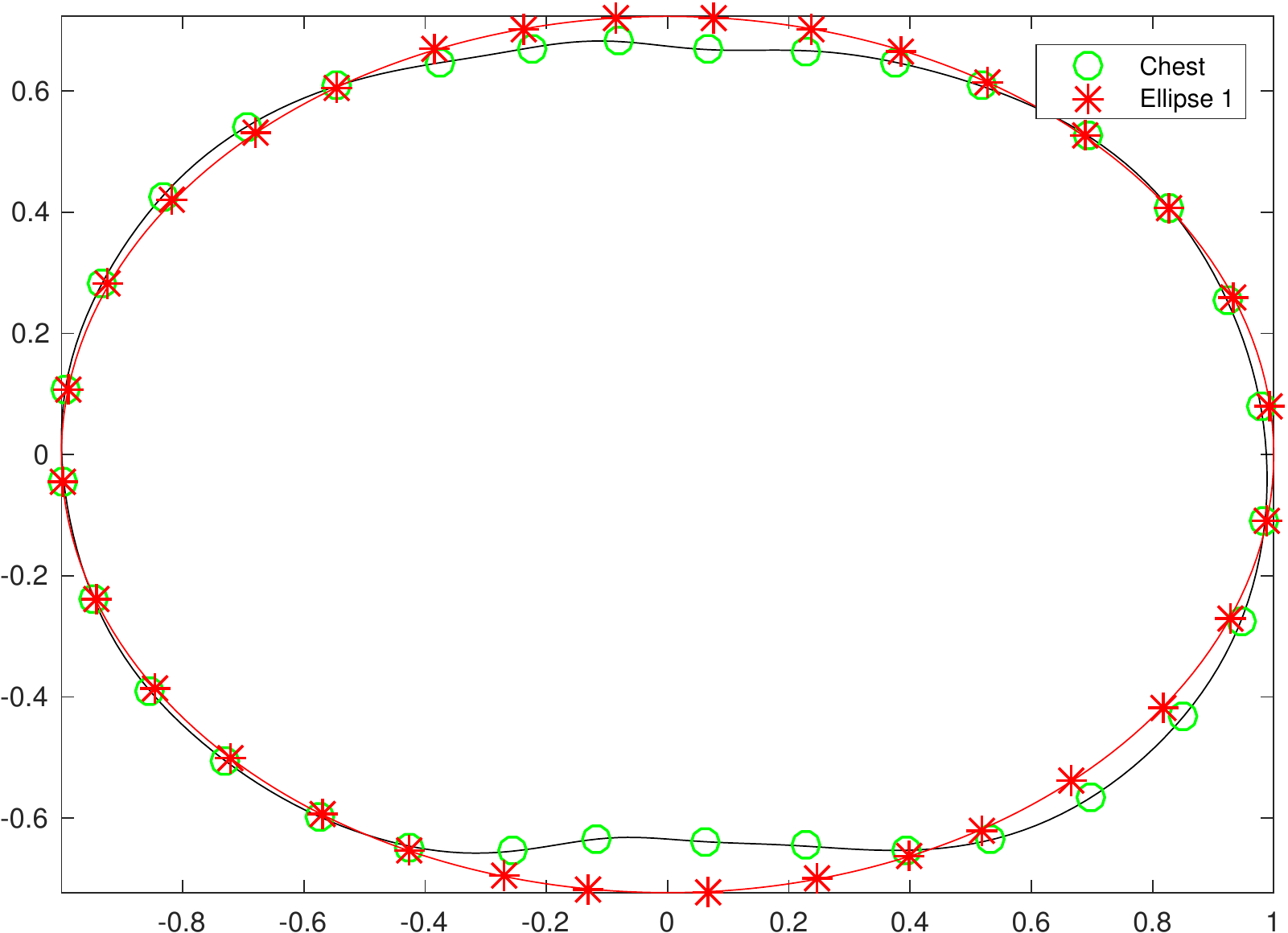}}
\put(220,0){\includegraphics[height=75pt]{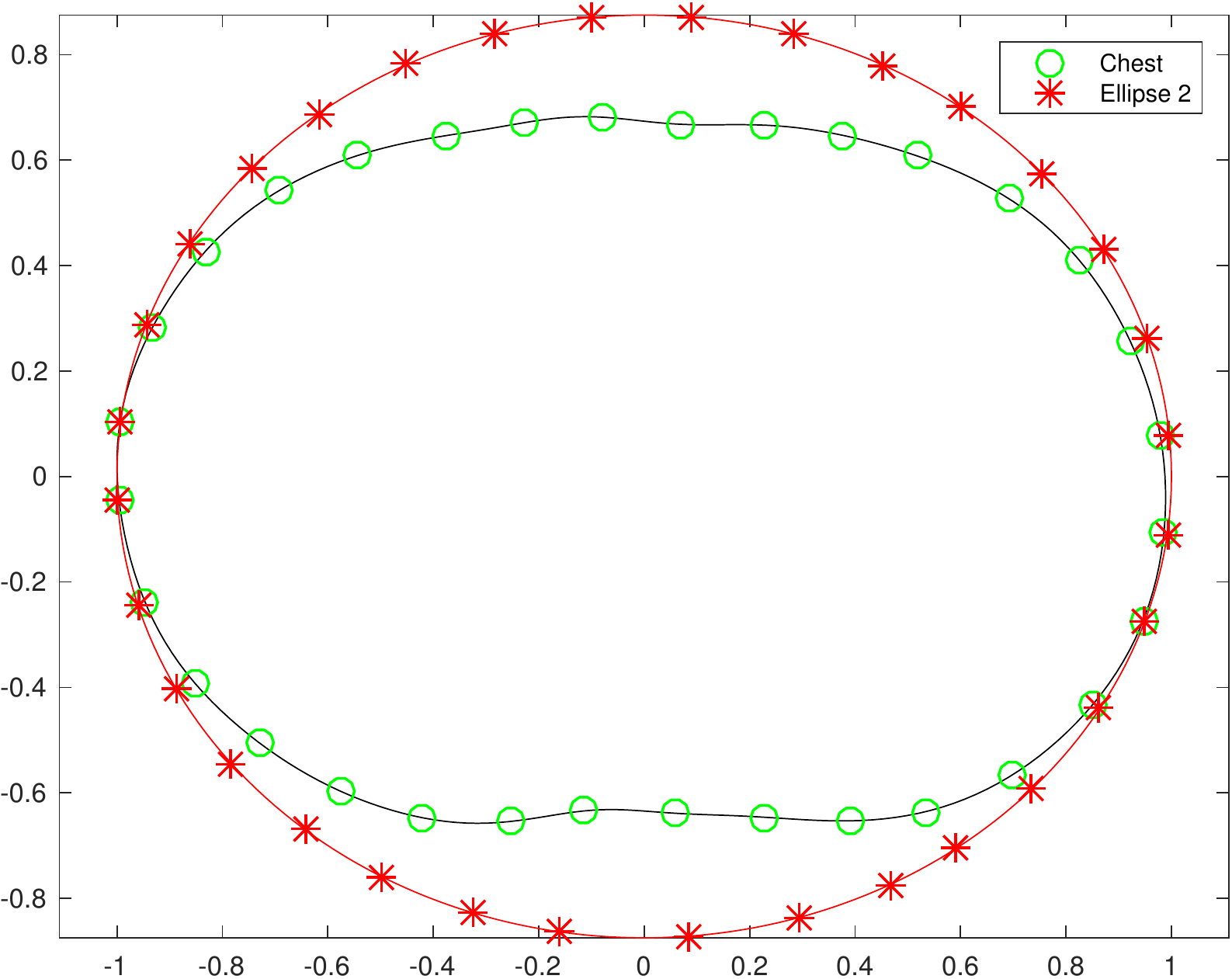}}
\put(330,0){\includegraphics[height=75pt]{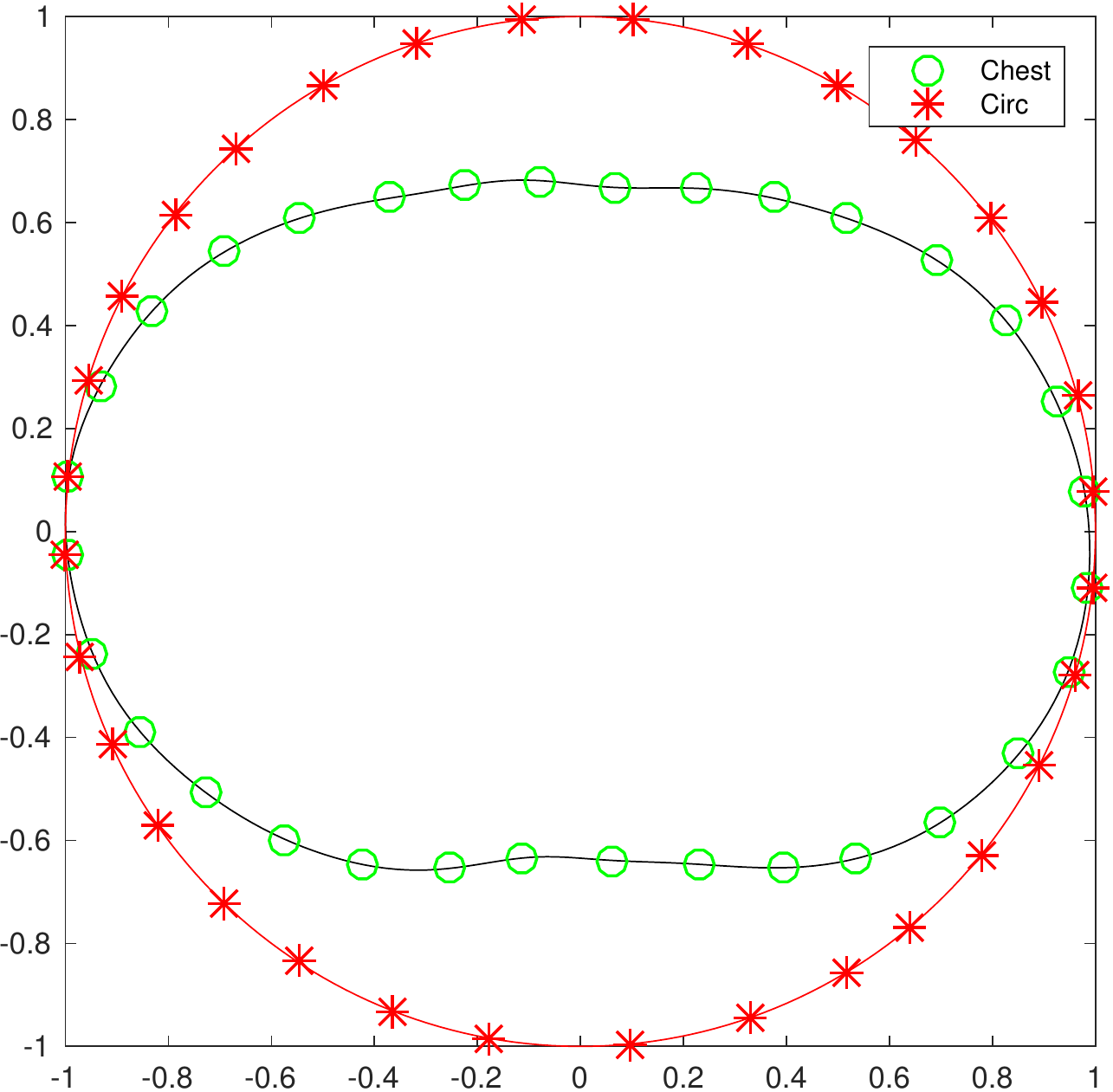}}

\put(10,85){\sc \underline{True Boundary}}
\put(135,85){\sc\underline{Ellipse 1}}
\put(245,85){\sc\underline{Ellipse 2}}
\put(350,85){\sc\underline{Circle}}

\end{picture}

\caption{Boundary shapes and electrode locations tested for the ACE1 phantoms. From \trev{left to right:} True boundary shape with correct and incorrect electrode angles, two ellipses of different eccentricities, and finally a circular boundary shape.   The green circles denote the `true' electrode locations whereas the red stars denote the noisy or incorrect electrode locations.}
\label{fig:CSU_elecSetups}
\end{figure}

Figure~\ref{fig:ACT4_circ_elecsNoisy} demonstrates the various boundary information scenarios that were tested for the ACT3 and ACT4 data.  The incorrect electrode locations are shown in the leftmost image in red with the markers representing the electrode centers, the ovular boundary is shown in the middle image, and another incorrect boundary is shown in the rightmost image.  For the ACE1 data, which has a chest-shaped boundary, we examine the alternative boundaries and electrode positions shown in Figure~\ref{fig:CSU_elecSetups} with incorrect electrode locations and boundaries again shown in red.

\trev{`True' boundary functions were formed as Fourier series approximations using coordinate data from imported photographs of the experiments.  The `alternative' boundary shown in Figure~\ref{fig:ACT4_circ_elecsNoisy} (right) was formed using imprecise clicks on the imported experiment photo (Fig.~\ref{fig:phantoms}, left).  The ovular boundaries were defined in terms of their major and minor principal semi-axes.  `Correct Angle' electrode locations were simulated by defining their centers $d\theta=\frac{2\pi}{L}$  apart.}

%--------------------------------------------------------
% Pics of the various incorrect data scenarios for a circular tank.
%--------------------------------------------------------
%--------------------------------------------------------

%%--------------------------------------------------------
%% ACE1 - Chest Shaped Tank - Incorrect Electrode Locations
%%--------------------------------------------------------
%\begin{figure}[h!]
%\begin{picture}(405,110)
%
%\put(0,0){\includegraphics[height=75pt]{ACE1chest_NoisyElecs-eps-converted-to.pdf}}
%\put(105,0){\includegraphics[height=75pt]{ACE1chest_Oval1-eps-converted-to.pdf}}
%\put(220,0){\includegraphics[height=75pt]{ACE1chest_Oval2-eps-converted-to.pdf}}
%\put(330,0){\includegraphics[height=75pt]{ACE1chest_Circ-eps-converted-to.pdf}}
%
%\put(10,85){\sc \underline{True Boundary}}
%\put(135,85){\sc\underline{Ellipse 1}}
%\put(245,85){\sc\underline{Ellipse 2}}
%\put(350,85){\sc\underline{Circle}}
%
%\end{picture}
%\caption{Boundary shapes and electrode locations tested for the ACE1 agar heart and lungs phantom. From left to right: True boundary shape with correct and incorrect electrode angles, two ellipses of different eccentricities, and finally a circular boundary shape.   The green circles denote the `true' electrode angles whereas the red stars denote the noisy electrode angles.}
%\label{fig:CSU_elecSetups}
%\end{figure}
%%--------------------------------------------------------
%--------------------------------------------------------
% ACE1 - Chest Shaped Tank - Incorrect Electrode Locations
%--------------------------------------------------------
%--------------------------------------------------------

For both the ACT4 and ACE1 data, the measured currents/voltages were used to synthesize the measurements that would have occurred if trigonometric current patterns of amplitude 1 mA had been applied.  Reconstructions for all examples were performed on a uniform $33\by33$ $k$-grid of stepsize $h_k=0.4706$ using a cutoff radius $R=4.0$ and nonuniform truncation threshold of 0.4.  The threshold was enforced such that $\Sexponetwo(k)=0$ if $|\Re(\Sexponetwo(k))|>0.4$ or $|\Im(\Sexponetwo(k))|>0.4$, and similarly for $\Sexptwoone$.  For consistency, the DN matrices $\mathbf{L}_\gamma$, $\mathbf{L}_\sigma$, and $\mathbf{L}_1$ were all scaled by $r$, the radius of the smallest circle enclosing the physical domain, effectively scaling the problem to be contained in the unit disc.  The spatial grid was represented by a uniform $64\by64$ $z$-grid on $[-1.05,1.05]^2$ with stepsize $h_z=0.0677$.  

For a point of comparison, we include reconstructions of the \trev{ACE1 data}, computed by the Gauss-Newton reconstruction algorithm.  The inverse problem was solved by minimizing the objective functional 
\begin{equation}\label{eq:objective_functional}
%r({\zeta}) = \parallel {\phi}_m - {\phi}_{c} \parallel_{2} + \alpha^2 \parallel ({\zeta}-{\zeta}\;^*){F^T F}({\zeta}-{\zeta}\;^*) \parallel_2,  
r({\zeta}) = \parallel ( {\phi}_m - {\phi}_{c} )\parallel_{2}^2 + \alpha^2 \parallel {F(\zeta}-{\zeta}\;^*) \parallel_{2}^2, 
\end{equation}
where $\phi_m$ is the vector of measurements, $\phi_c$ is the computed forward problem, $\zeta$ is the impeditivity distribution of the domain, $\alpha$ is the regularization parameter, $F$ is a Gaussian high-pass filter and ${\zeta}\;^*$ is a constant estimate of the impeditivity that also is used as the initial guess.   
The weight of the regularization was adjusted by visual examination, and in the reconstructions presented was $\alpha = 0.005$.  To reduce computational effort and time consumption, the number of elements in the finite element mesh is reduced by applying the approximation error theory [\cite{Kaipio2004a}], where a Bayesian modeling error approach is used to treat approximation and modeling errors.

 In this implementation, the forward problem was solved at each iteration with a finite element method with 12,000 elements.  The algorithm converged in five iterations, and took approximately 20 minutes to obtain each image. 
Before running the minimization algorithm, the approximation error vector is computed requiring additional 30 minutes of computation time.   

\trev{The maximum, minimum, and average reconstructed values of conductivity, and permittivity where applicable, in each of the organ regions were calculated and are provided in the tables in Section~\ref{sec:results} for each reconstruction method and dataset.  The organ boundaries were identified from the photos in Figure~\ref{fig:phantoms} and superimposed on each of the reconstructions.  No distortion of the superimposed organ shapes was imposed for the reconstructions on the incorrect boundary shapes.}  

The dynamic range of each conductivity reconstruction was computed by the formula 
\begin{equation} \label{DynRange}
\mbox{Dynamic Range} = \frac{\sigma^{recon}_{max} - \sigma^{recon}_{min}}{\sigma^{true}_{max} -\sigma^{true}_{min}} \times 100\%,
\end{equation}
where $\sigma^{recon}_{max}$ is the maximum value over all pixels in the reconstruction, and $\sigma^{true}_{max}$ is the maximum value over each of the targets as measured by the conductivimeter, and $\sigma^{recon}_{min}$ and $\sigma^{true}_{min}$ are defined analogously.  \trev{The dynamic range for the susceptivity is similarly defined.}  The {\it `Approach 2'} D-bar reconstructions computed here required approximately \trev{13} seconds of computation time using a four-core \trev{laptop} computer with \trev{2.9}~GHz Intel Core i7 processors \trev{and} a {\sc Matlab} implementation that has not been optimized for speed.

%--------------------------------------------------------------------
\section{Results and Discussion}\label{sec:results}
%--------------------------------------------------------------------
%--------------------------------------------------------
\subsection{Reconstructions from ACT3 Data}\label{sec:results_ACT3}
\begin{figure}[!h]

\begin{picture}(300,230)

% Difference images
\put(20,0){\includegraphics[width=80pt]{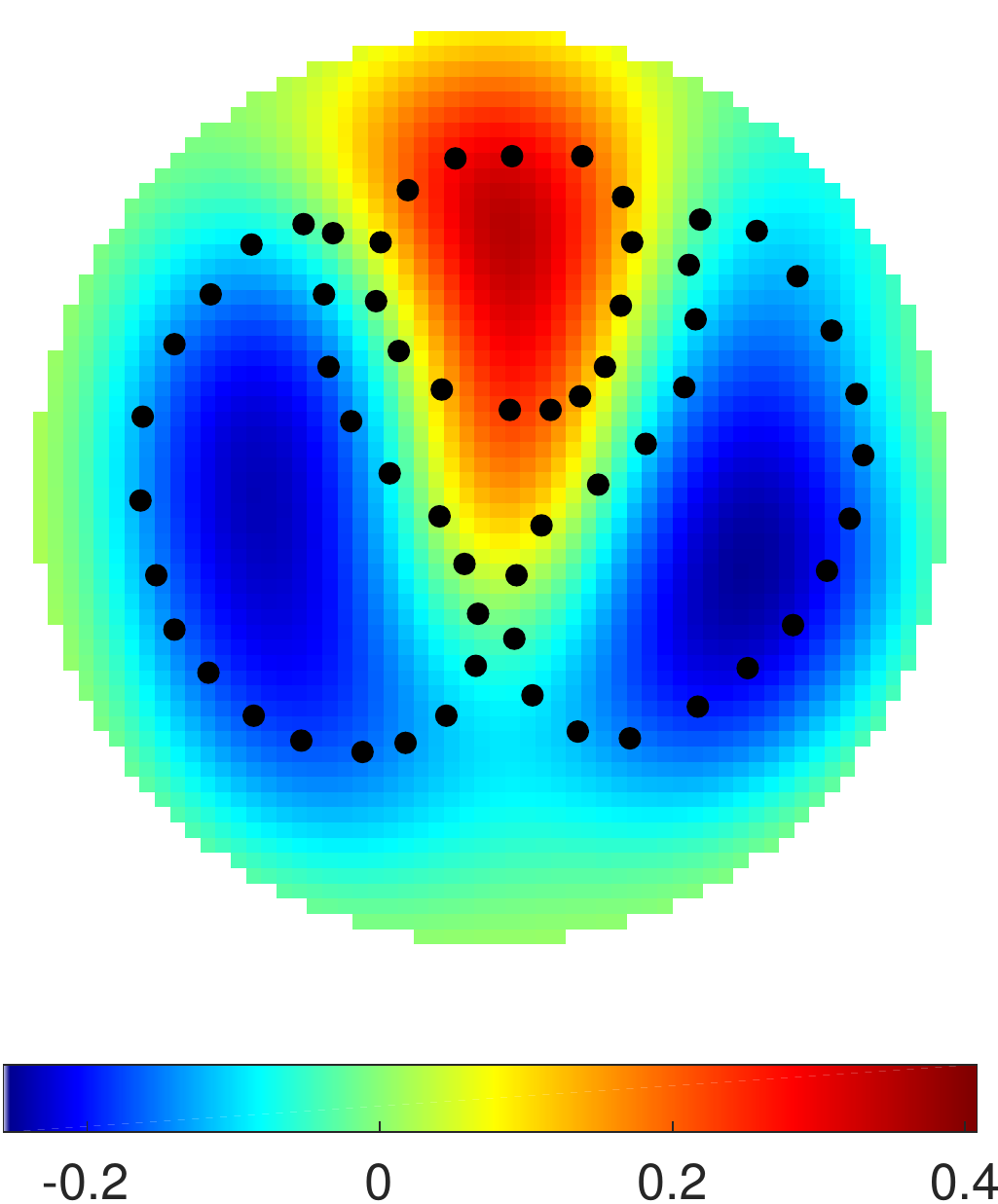}} % True angles and boundary image
\put(120,0){\includegraphics[width=80pt]{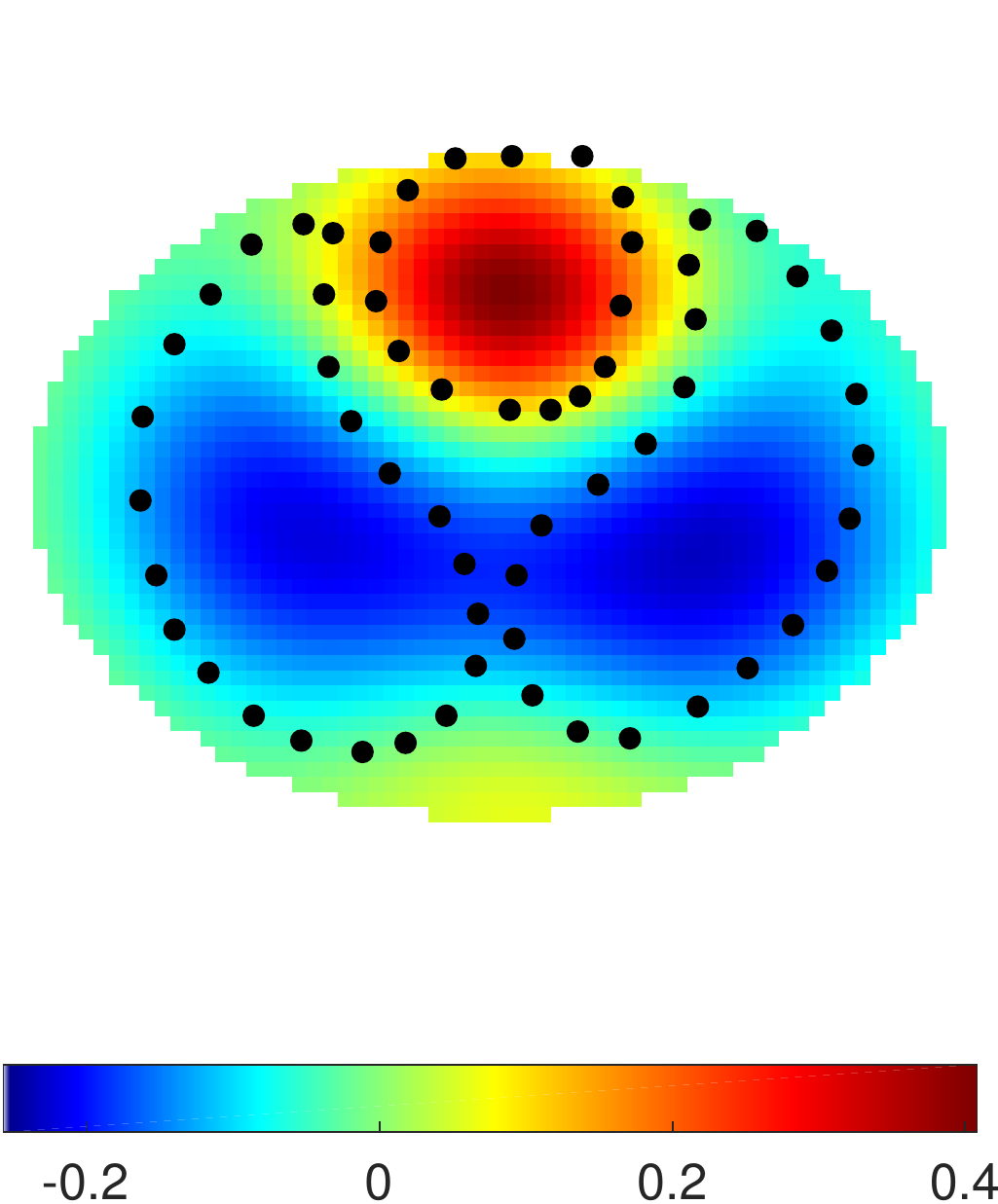}} % Incorrect electrode angles with correct boundary
\put(220,0){\includegraphics[width=80pt]{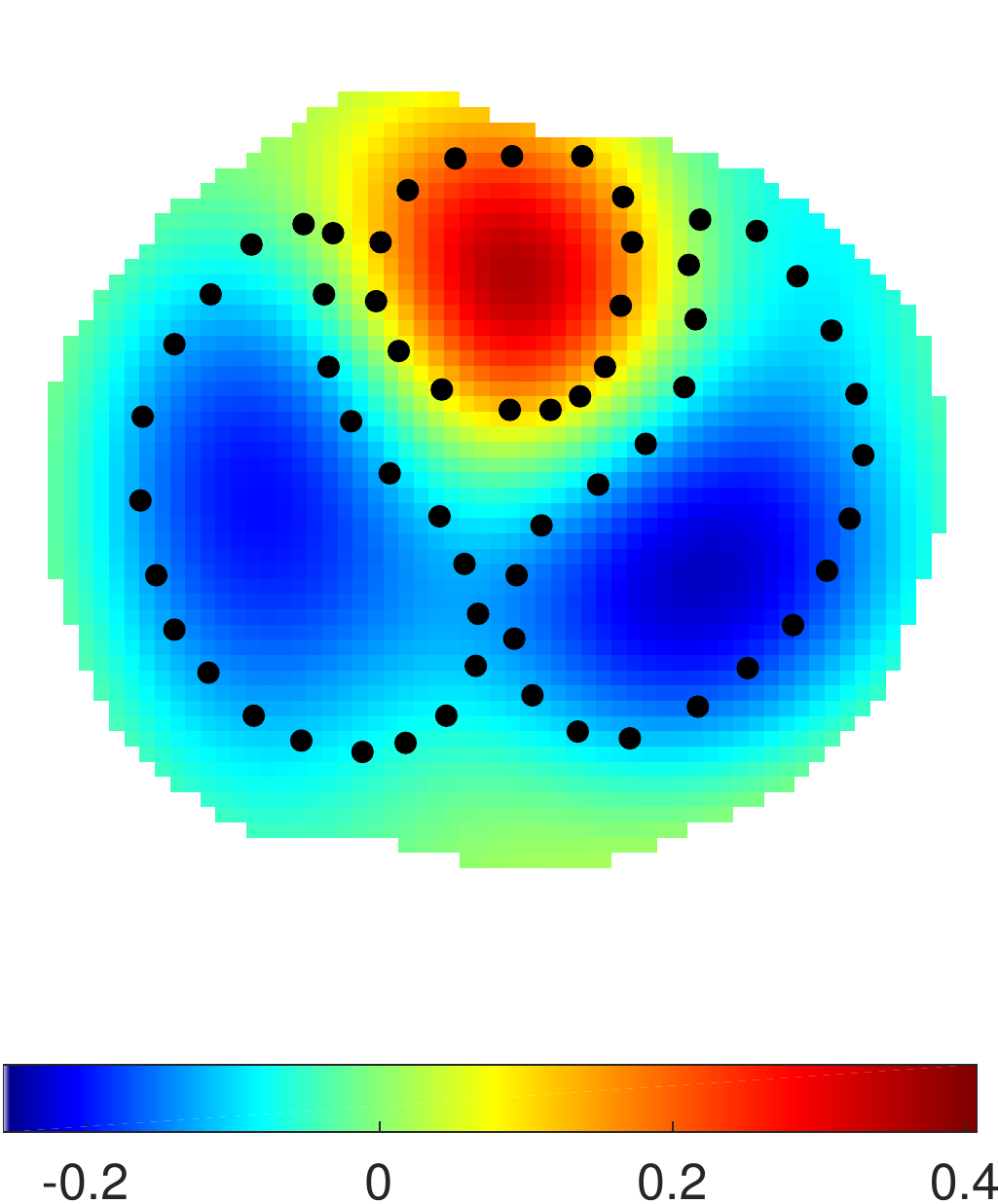}} % True angles and oval1 boundary
%\put(320,0){\includegraphics[width=80pt]{GAM_ACT3_HnL_tExp_Chest1_ElecsTrue_DIFF_PriorOverlay_R_jet-eps-converted-to.pdf}} % True angles and oval1 boundary

% Absolute images
\put(20,105){\includegraphics[width=80pt]{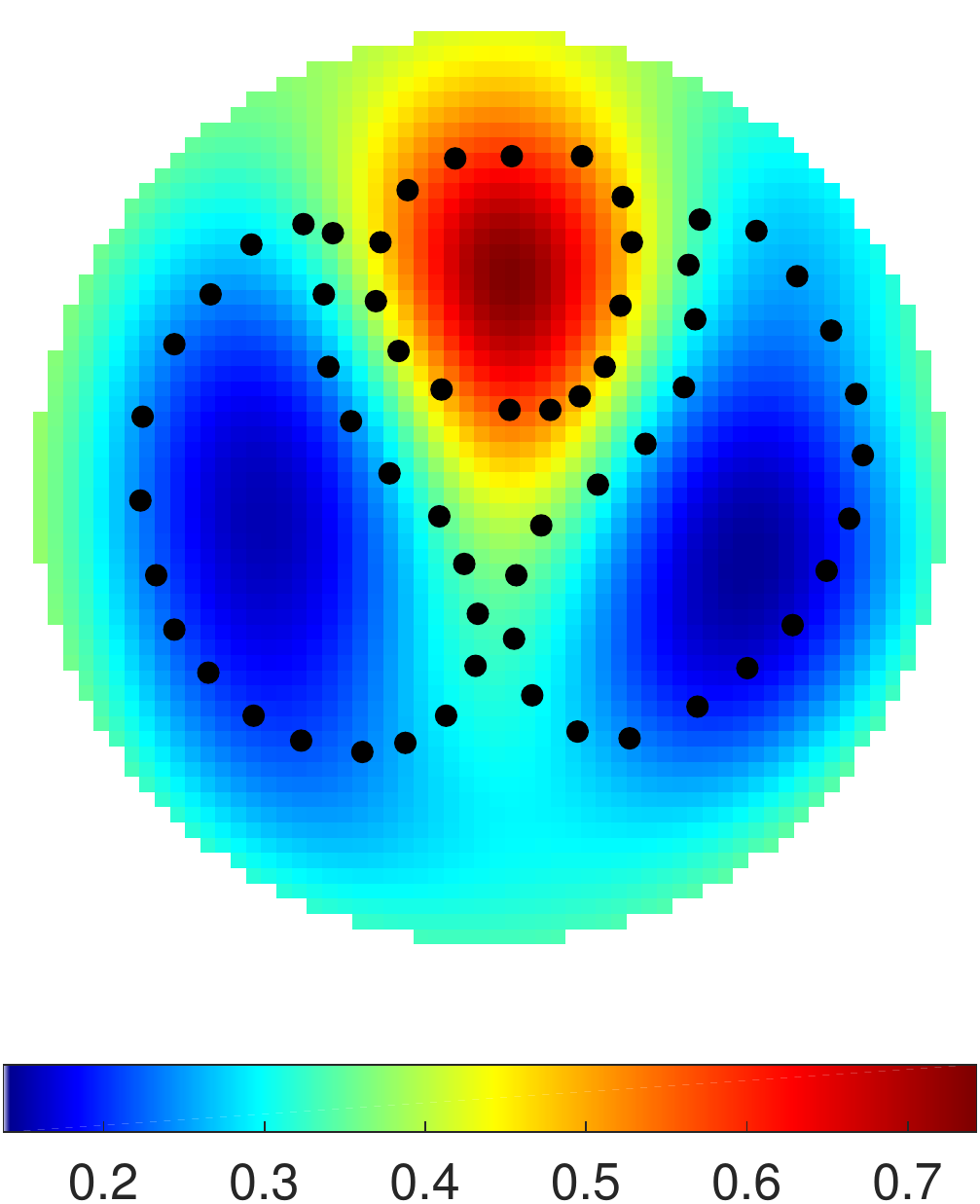}} % True angles and boundary image
\put(120,105){\includegraphics[width=80pt]{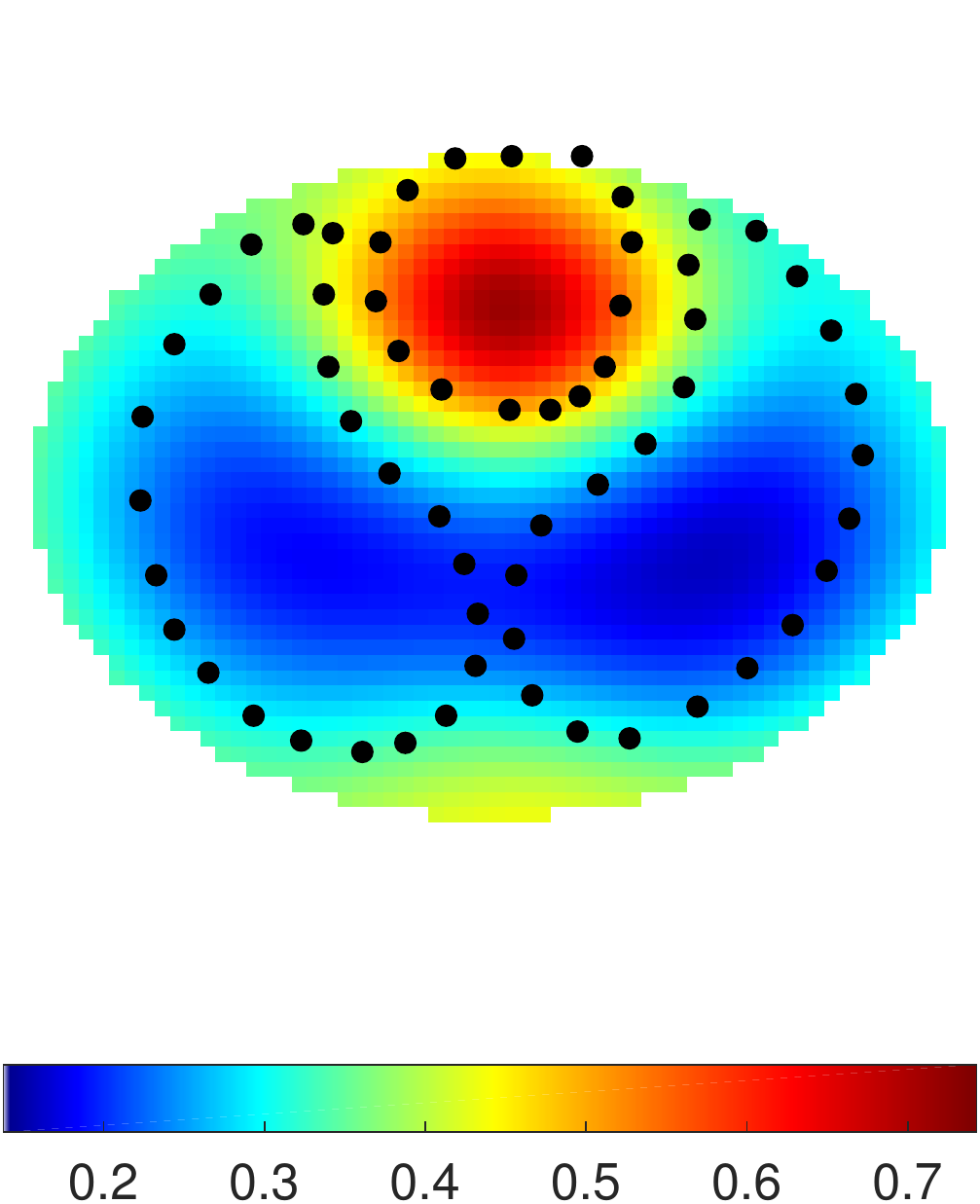}} % Incorrect electrode angles with correct boundary
\put(220,105){\includegraphics[width=80pt]{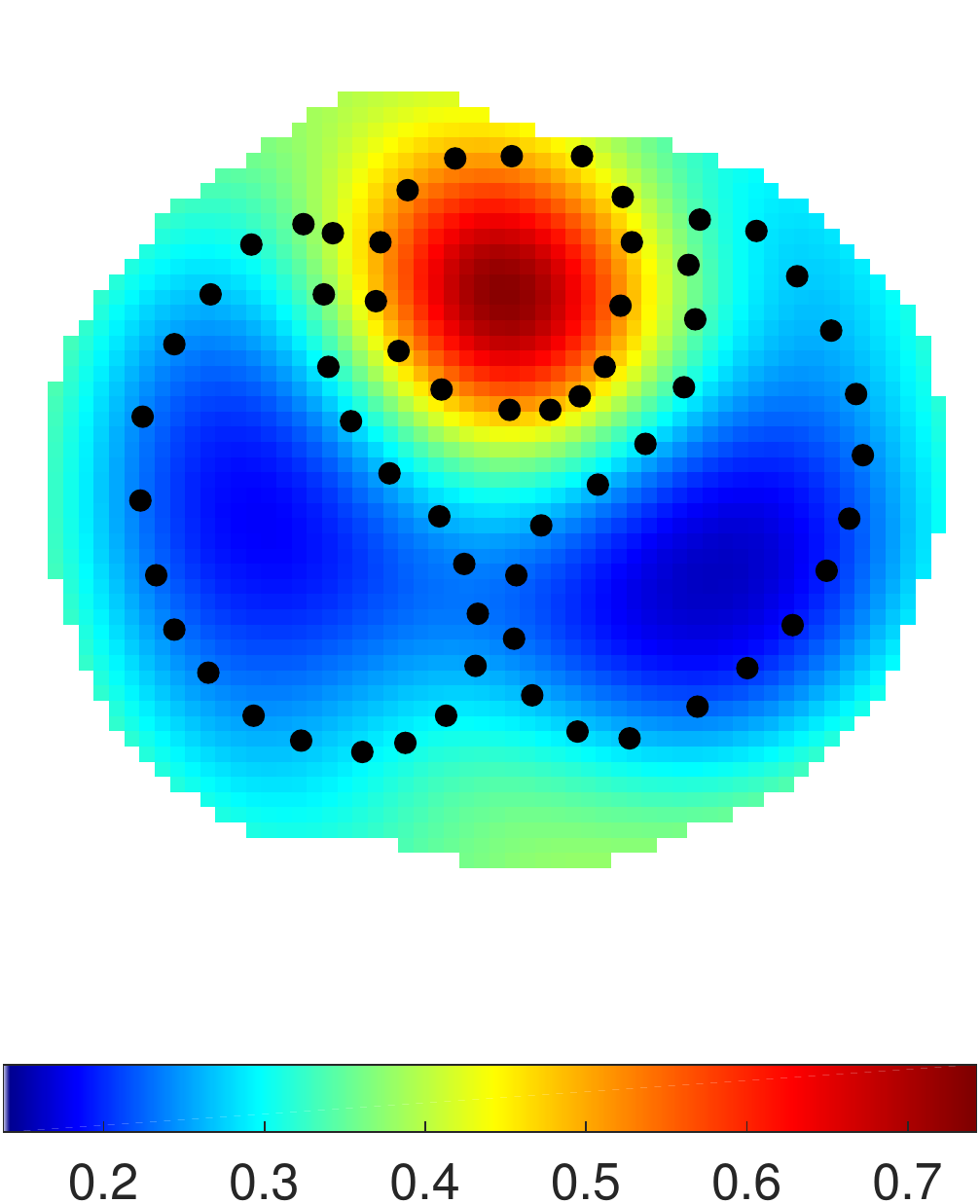}} % True angles and oval1 boundary
%\put(320,105){\includegraphics[width=80pt]{GAM_ACT3_HnL_tExp_Chest1_ElecsTrue_ABS_PriorOverlay_R_jet-eps-converted-to.pdf}} % True angles and clicks boundary

% side labels 
\put(0,30){\rotatebox{90}{\scriptsize{\sc Difference}}}
\put(8,42){\rotatebox{90}{\scriptsize{\sc Images}}}

\put(0,135){\rotatebox{90}{\scriptsize{\sc Absolute}}}
\put(8,145){\rotatebox{90}{\scriptsize{\sc Images}}}

% top labels:
\put(35,215){\footnotesize {\sc True Angles}}
\put(30,205){\footnotesize {\sc True Boundary}}

\put(135,215){\footnotesize {\sc True Angles}}
\put(125,205){\footnotesize {\sc Oval Boundary}}

\put(235,215){\footnotesize {\sc True Angles}}
\put(215,205){\footnotesize {\sc Alternative Boundary}}

\end{picture}
\caption{\label{fig:ACT3_HnL_Nach} Comparison of conductivity \trev{(S/m)} reconstructions from the ACT3 data (see Figure~\ref{fig:phantoms}) using  the $\texp$ approximation in Section \ref{sec:NachDbar} for knowledge of true vs. incorrect electrode angles as well as boundary shape.  Absolute images are in row 1\tRev{, plotted on the same color scale}.  Difference images are in row 2\tRev{, plotted on the same color scale.}}
\end{figure}

\begin{figure}[!h]

\begin{picture}(300,230)

% Difference images
\put(20,0){\includegraphics[width=80pt]{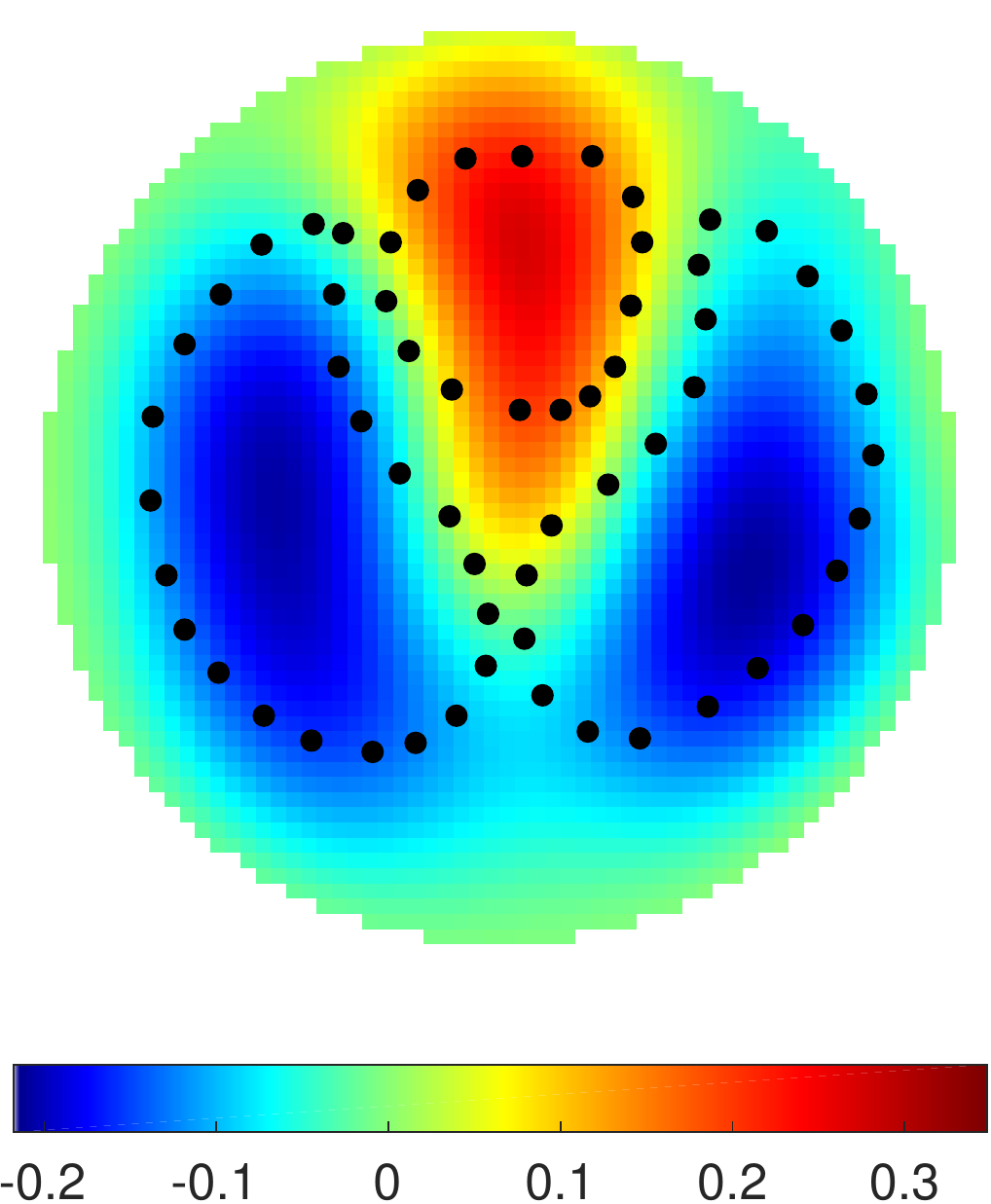}} % True angles and boundary image
\put(120,0){\includegraphics[width=80pt]{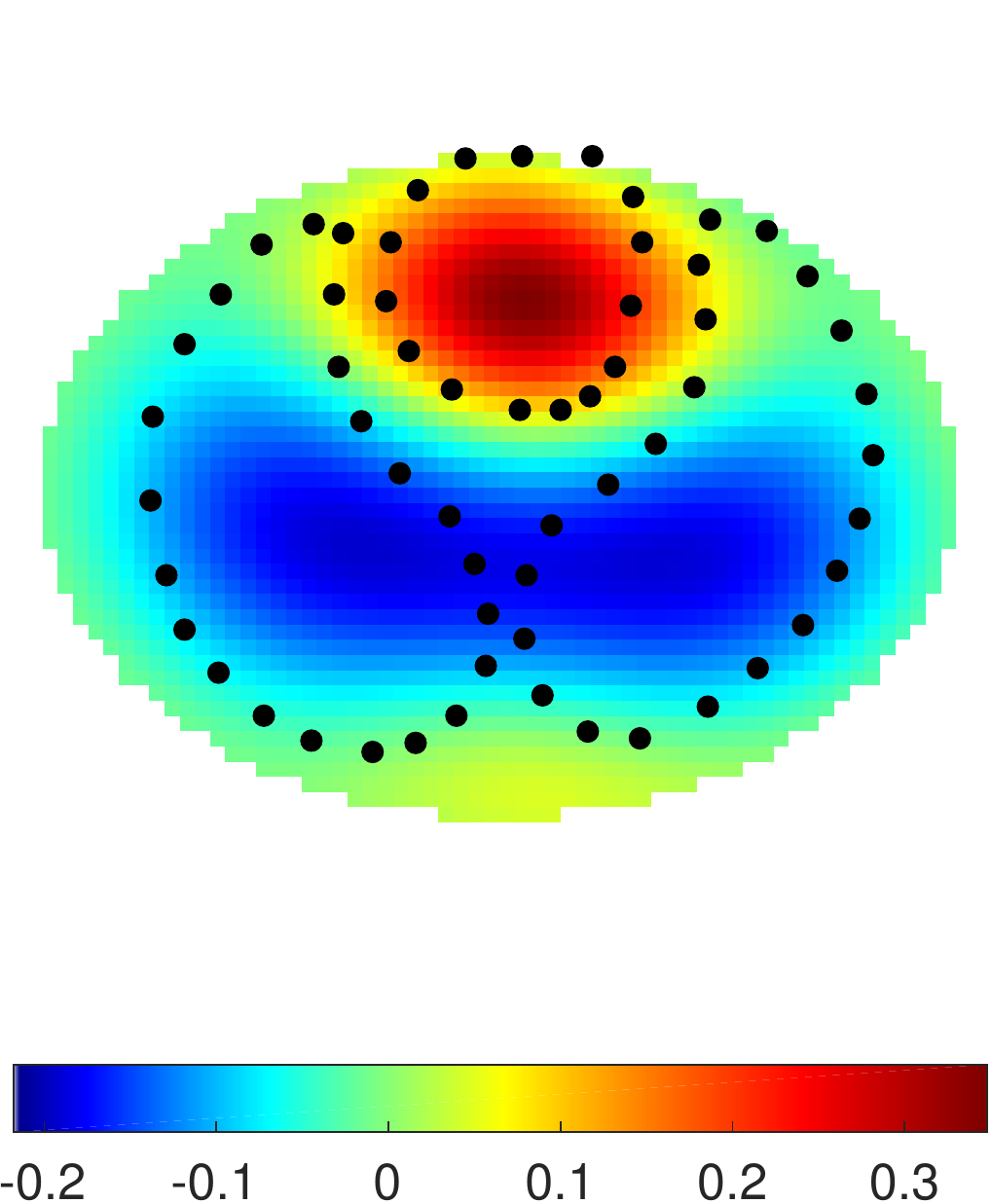}} % Incorrect electrode angles with correct boundary
\put(220,0){\includegraphics[width=80pt]{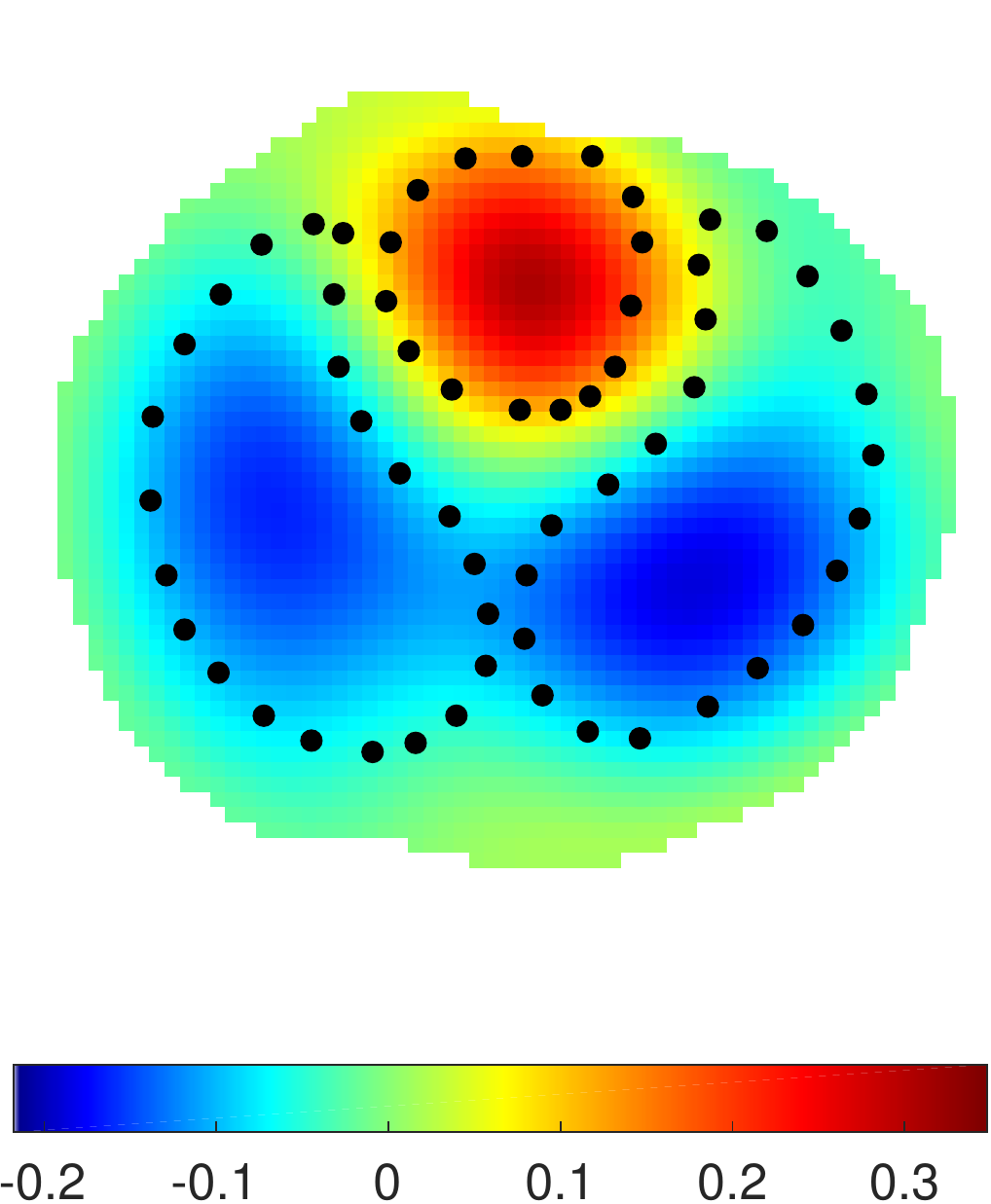}} % True angles and oval1 boundary
%\put(320,0){\includegraphics[width=80pt]{GAM_ACT3_HnL_PsiExpSarah_Chest1_ElecsTrue_DIFF_PriorOverlay_R_jet-eps-converted-to.pdf}} % True angles and oval1 boundary

% Absolute images
\put(20,105){\includegraphics[width=80pt]{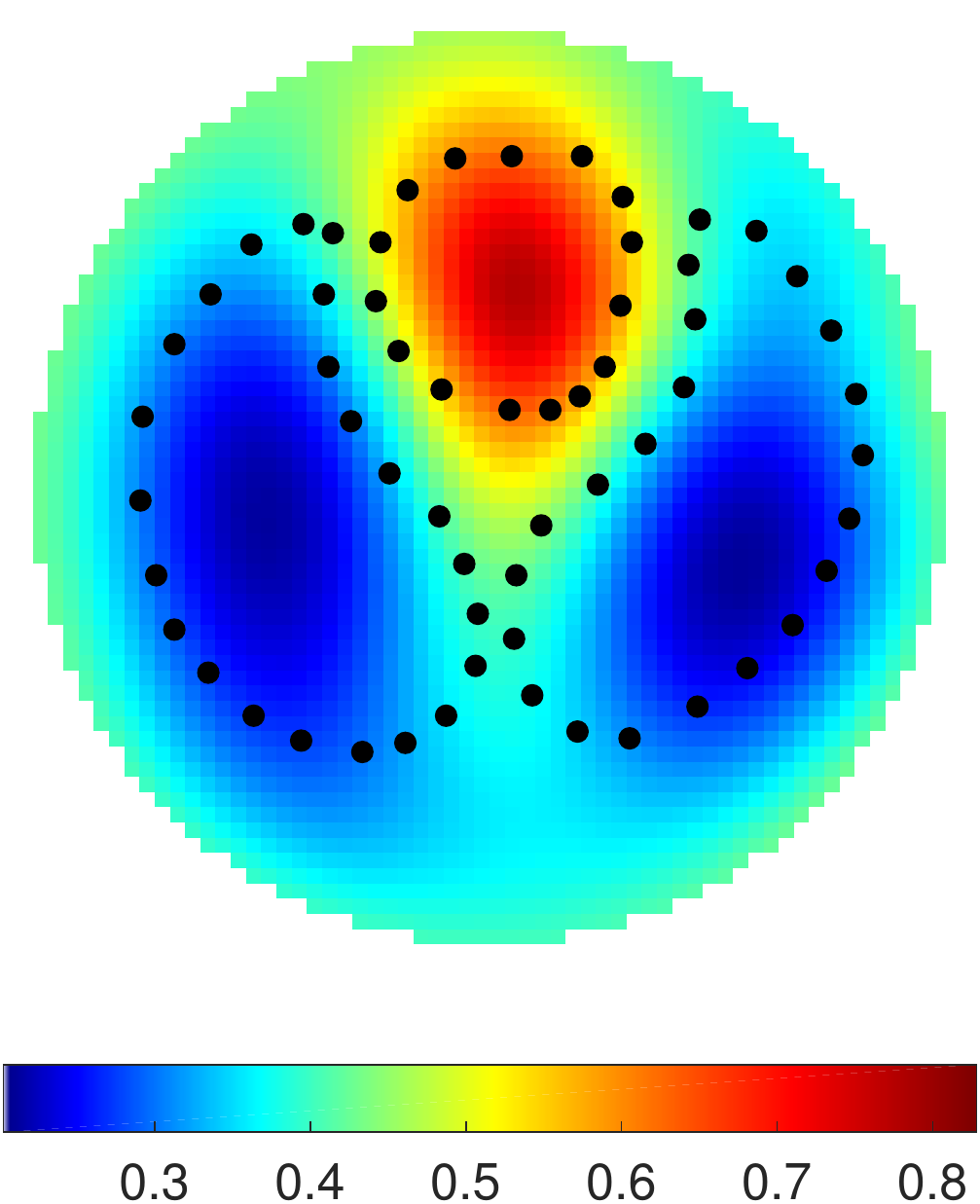}} % True angles and boundary image
\put(120,105){\includegraphics[width=80pt]{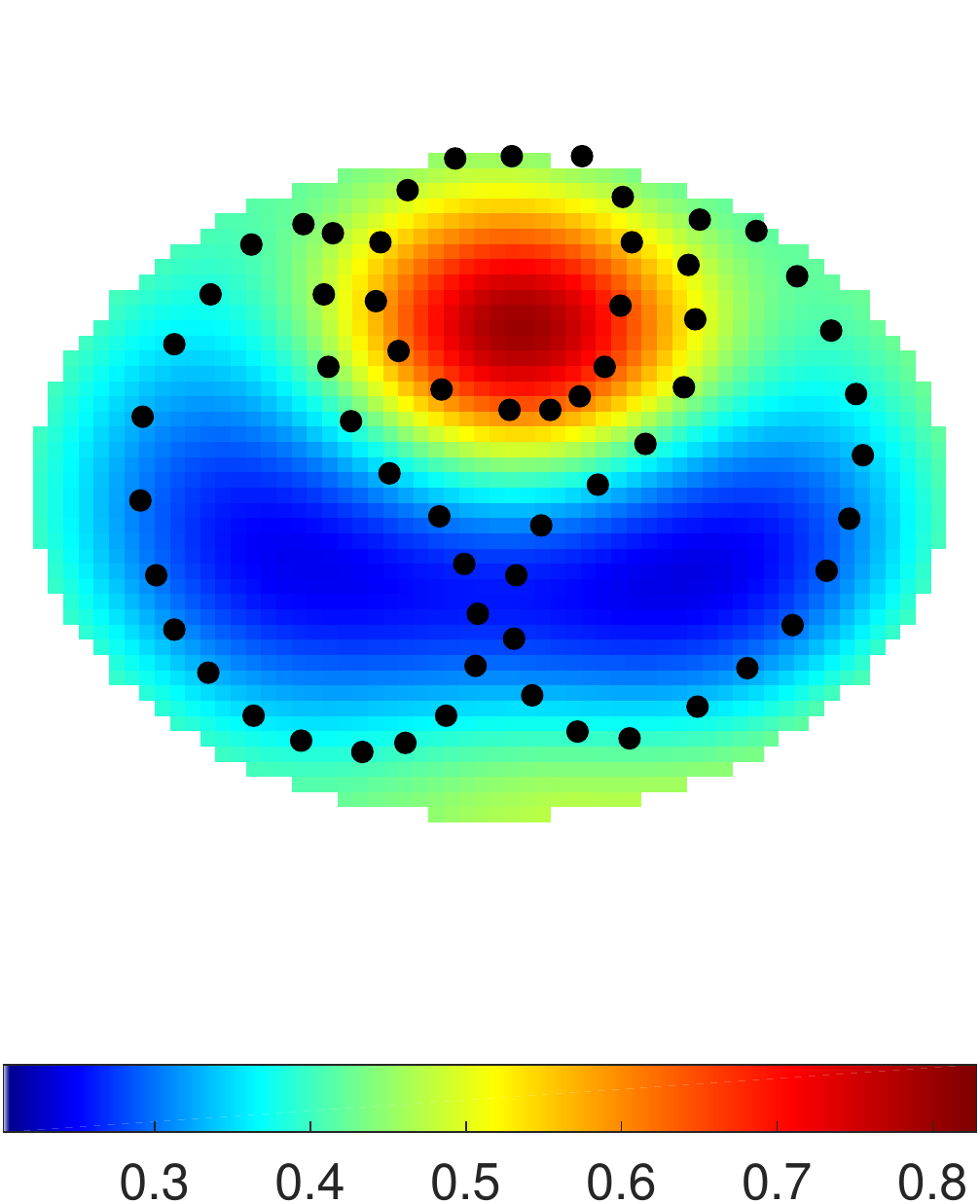}} % Incorrect electrode angles with correct boundary
\put(220,105){\includegraphics[width=80pt]{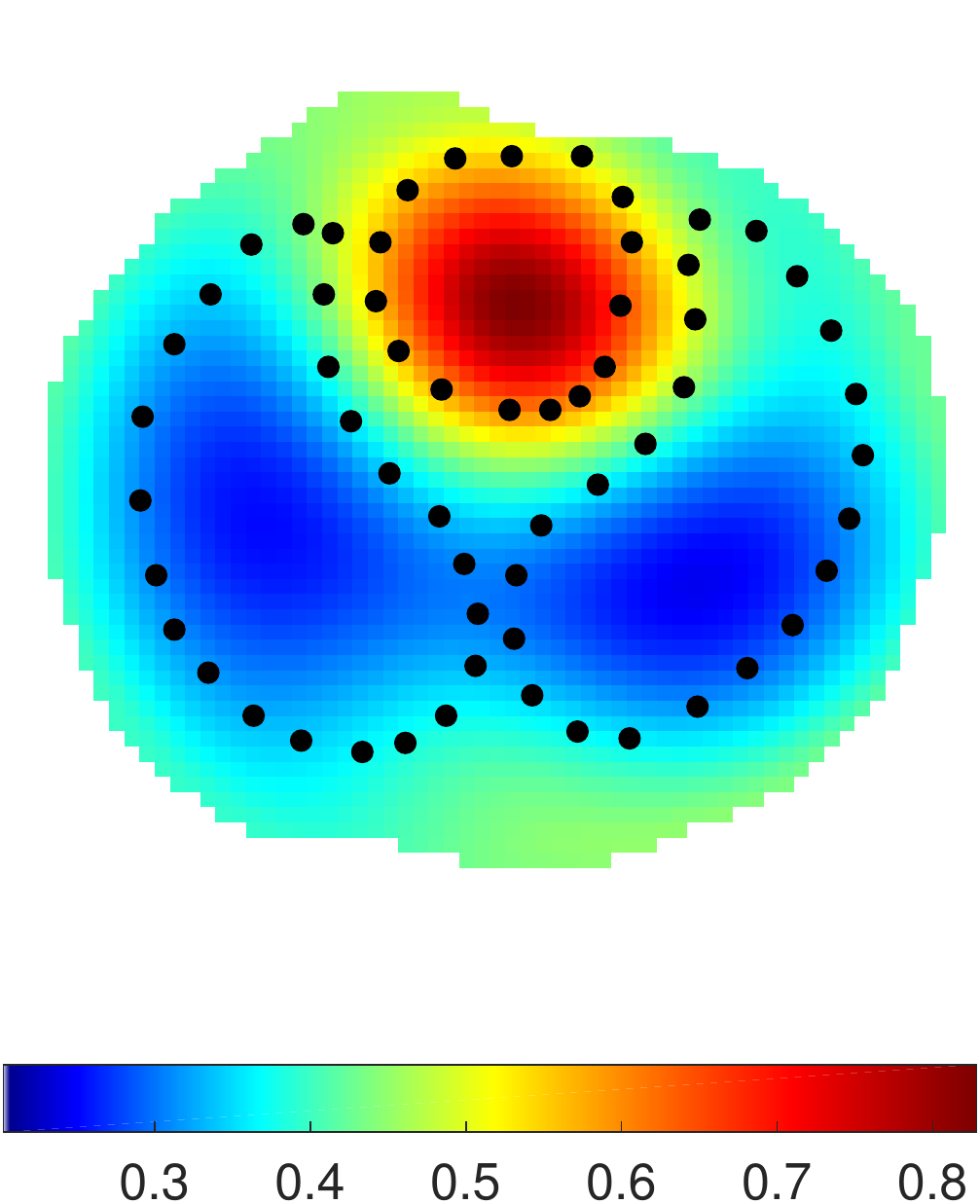}} % True angles and oval1 boundary
%\put(320,105){\includegraphics[width=80pt]{GAM_ACT3_HnL_PsiExpSarah_Chest1_ElecsTrue_ABS_PriorOverlay_R_jet-eps-converted-to.pdf}} % True angles and clicks boundary

% side labels 
\put(0,30){\rotatebox{90}{\scriptsize{\sc Difference}}}
\put(8,42){\rotatebox{90}{\scriptsize{\sc Images}}}

\put(0,135){\rotatebox{90}{\scriptsize{\sc Absolute}}}
\put(8,145){\rotatebox{90}{\scriptsize{\sc Images}}}

% top labels:
\put(35,215){\footnotesize {\sc True Angles}}
\put(30,205){\footnotesize {\sc True Boundary}}

\put(135,215){\footnotesize {\sc True Angles}}
\put(125,205){\footnotesize {\sc Oval Boundary}}

\put(235,215){\footnotesize {\sc True Angles}}
\put(215,205){\footnotesize {\sc Alternative Boundary}}

%\put(335,215){\footnotesize {\sc True Angles}}
%\put(325,205){\footnotesize {\sc Clicks Boundary}}
%\put(238,215){\footnotesize {\sc Med Prior}}
%\put(245,205){\scriptsize $(6,\frac23)$}

%\put(330,215){\footnotesize \textsc{Strong Prior}}
%\put(345,205){\scriptsize $(8,\frac13)$}

\end{picture}
\caption{\label{fig:ACT3_HnL_FranSarah} Comparison of conductivity \trev{(S/m)} reconstructions from the ACT3 data (see Figure~\ref{fig:phantoms}, first) using {\it `Approach~1'} of Section \ref{sec:FranDbar}  for knowledge of true vs. incorrect electrode angles as well as boundary shape.  Absolute images are in row 1\tRev{, plotted on the same color scale}.  Difference images are in row 2\tRev{, plotted on the same color scale.}}
\end{figure}

\begin{figure}[!h]

\begin{picture}(420,230)

% Difference images
\put(20,0){\includegraphics[width=80pt]{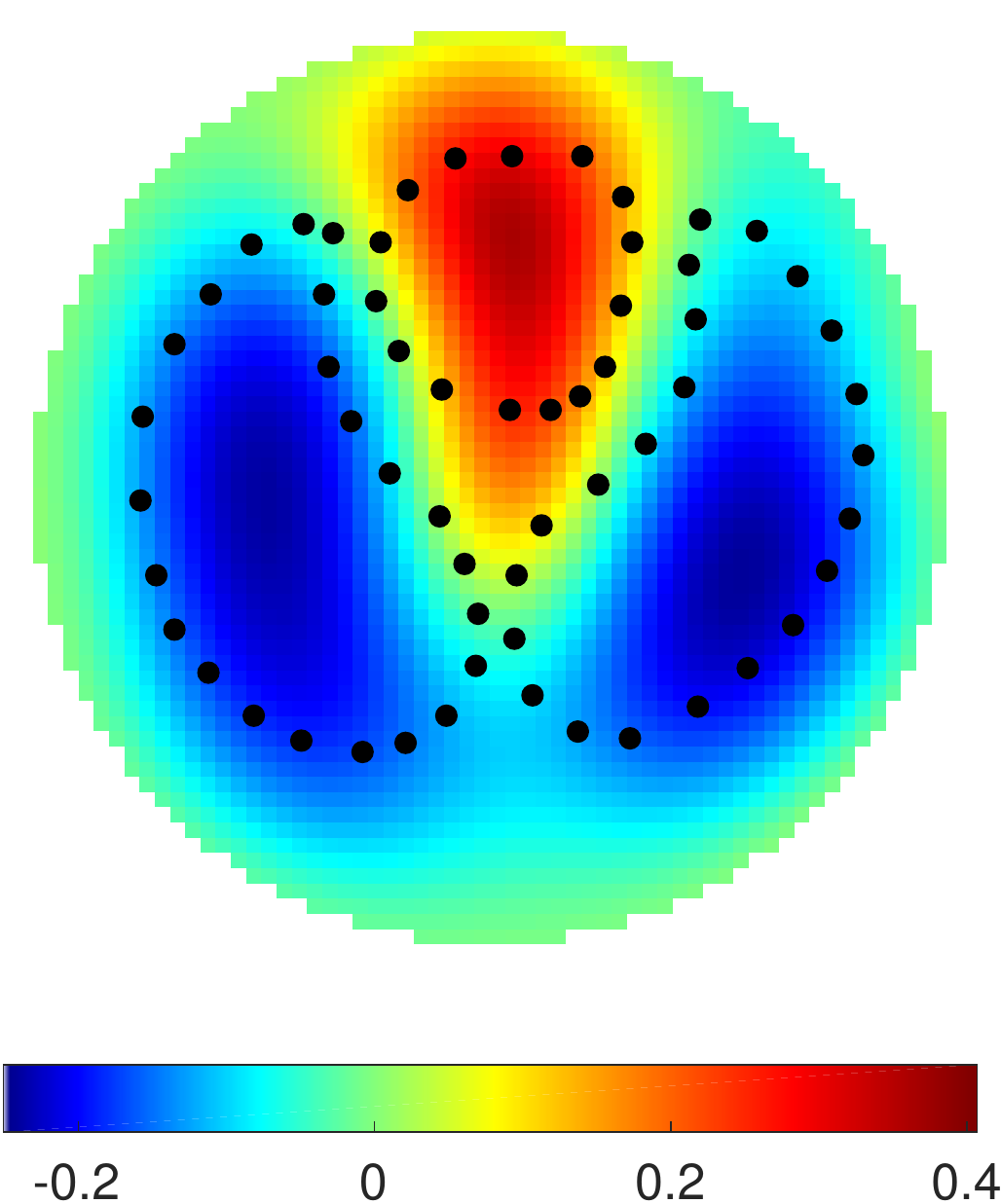}} % True angles and boundary image
\put(120,0){\includegraphics[width=80pt]{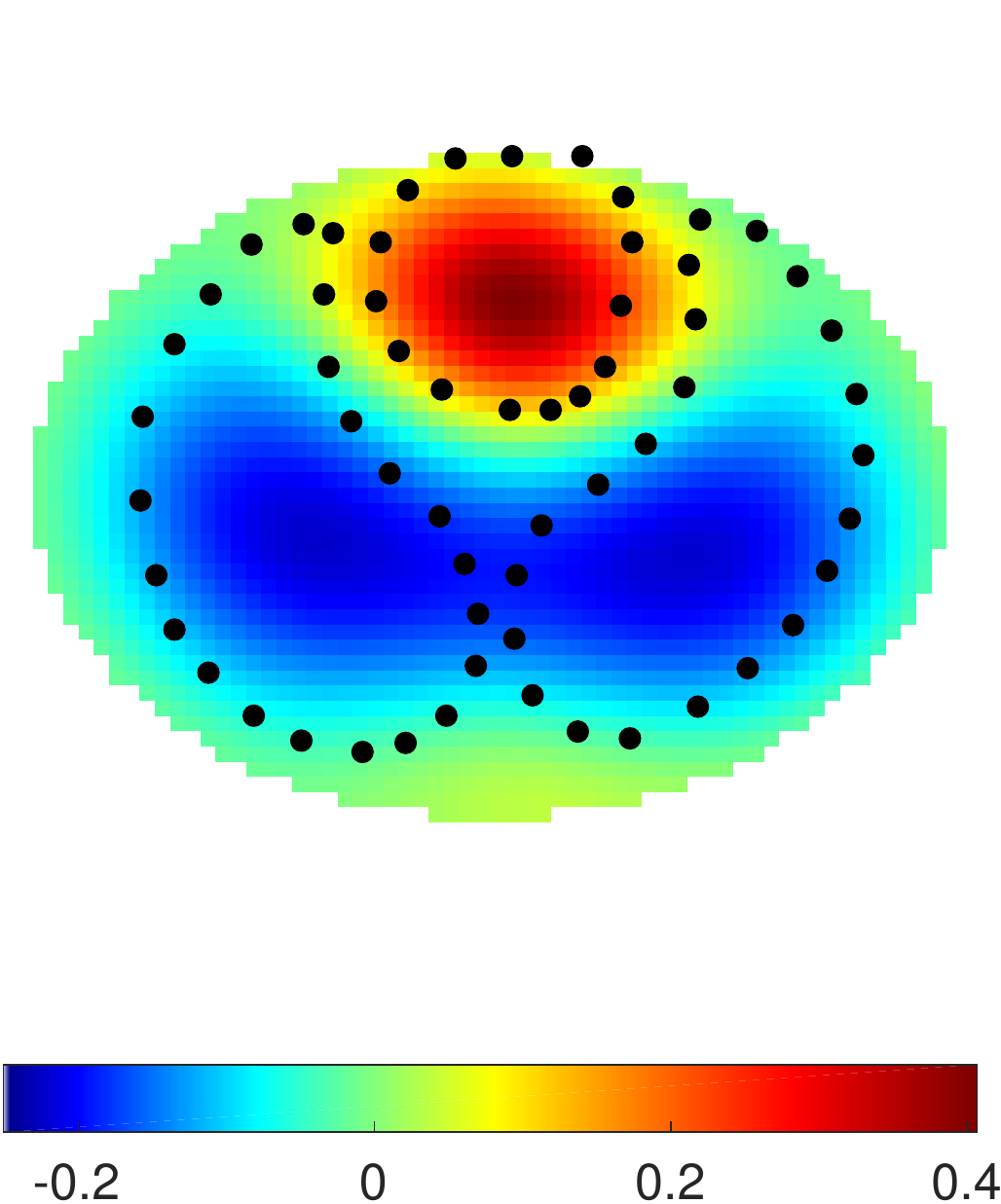}} % Incorrect electrode angles with correct boundary
\put(220,0){\includegraphics[width=80pt]{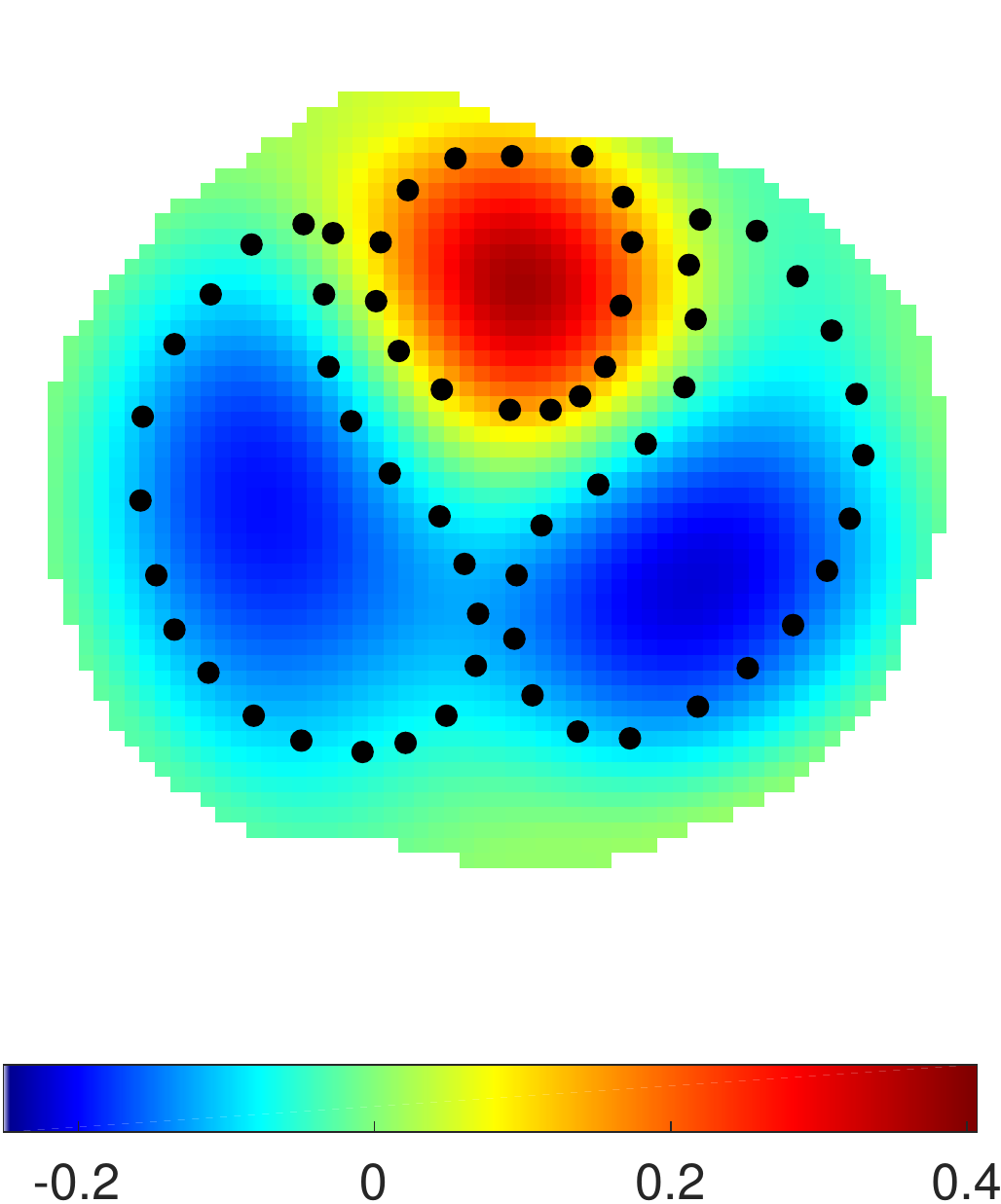}} % True angles and oval1 boundary
\put(320,0){\includegraphics[width=80pt]{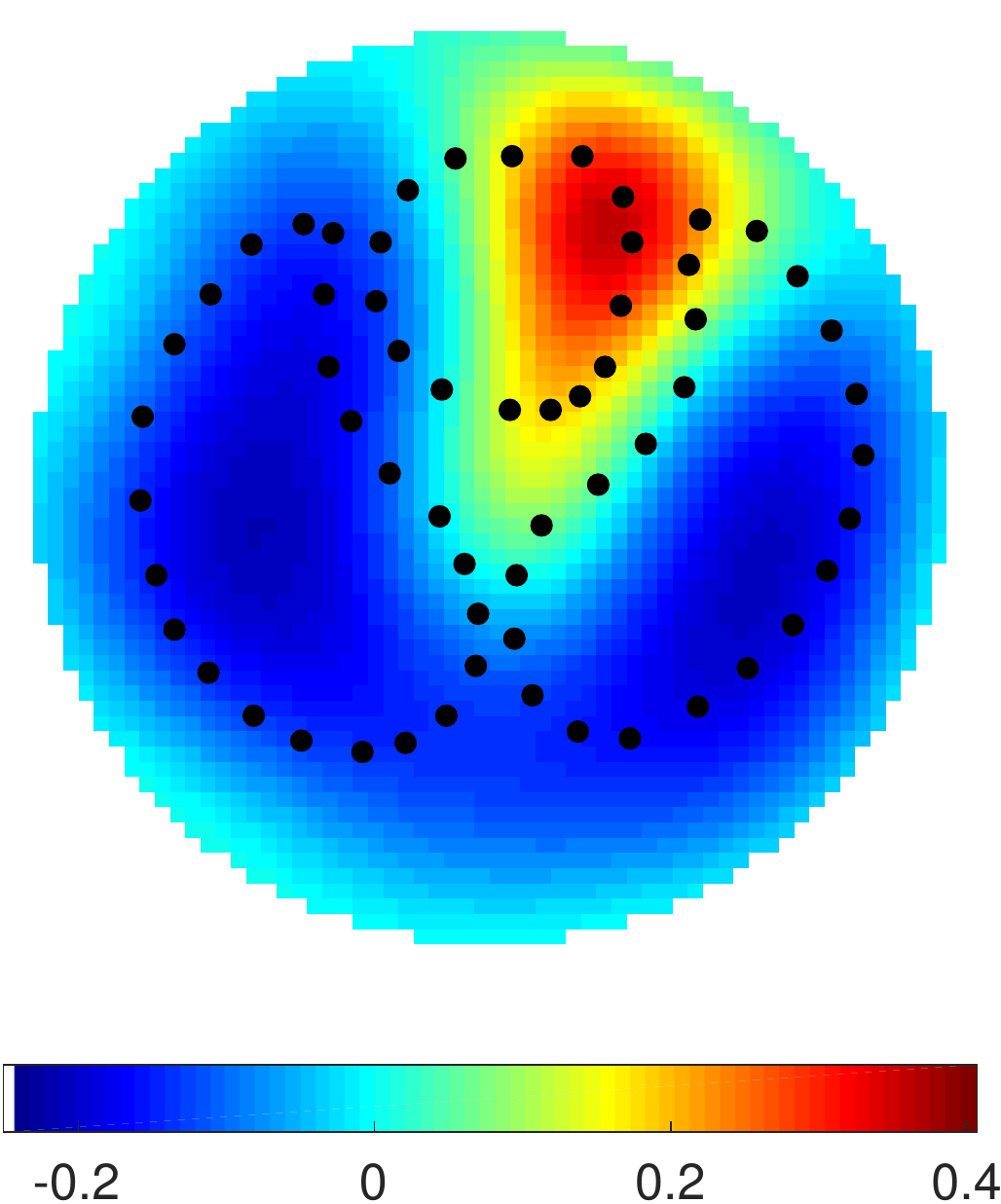}} % True angles and oval1 boundary

% Absolute images
\put(20,105){\includegraphics[width=80pt]{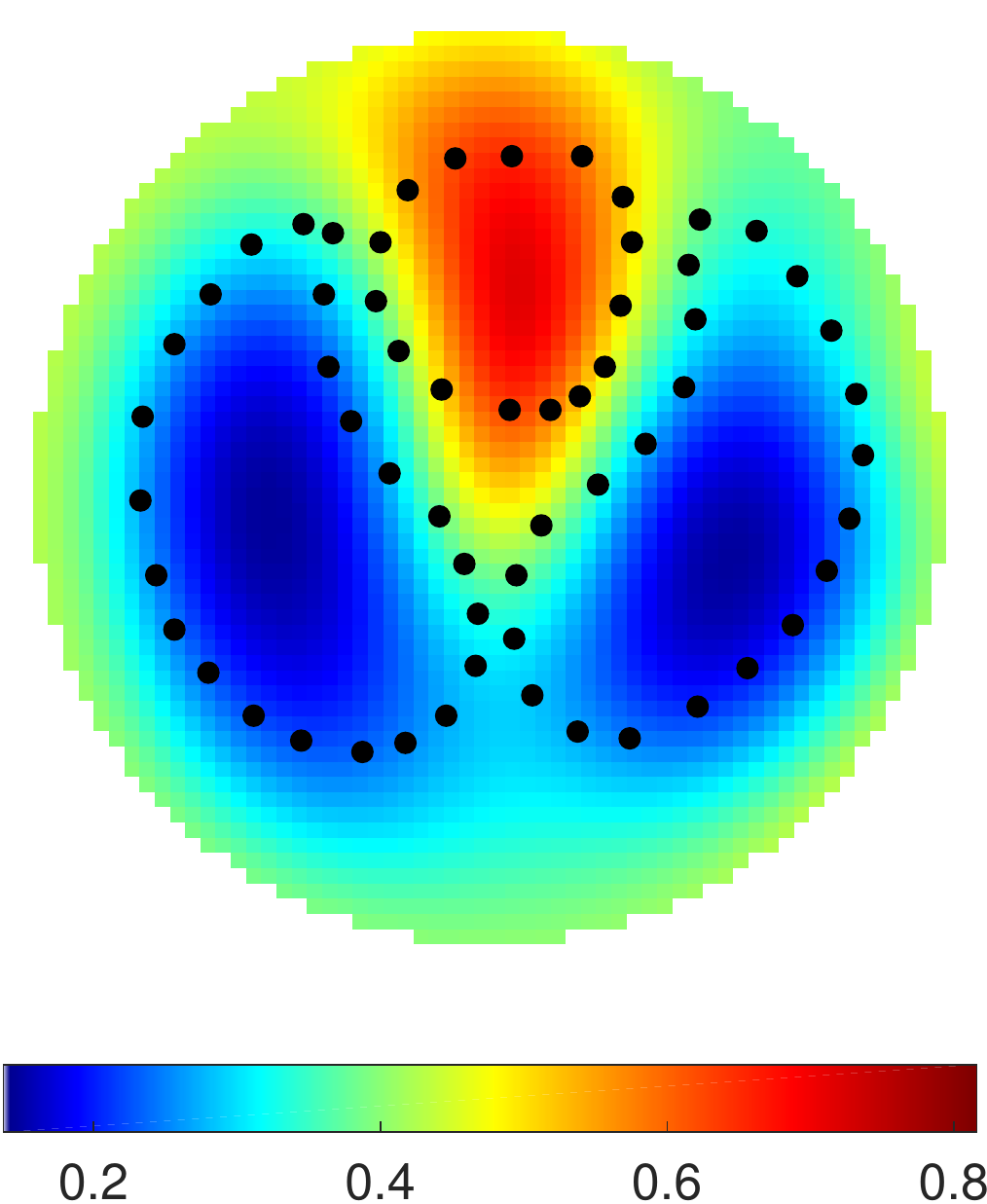}} % True angles and boundary image
\put(120,105){\includegraphics[width=80pt]{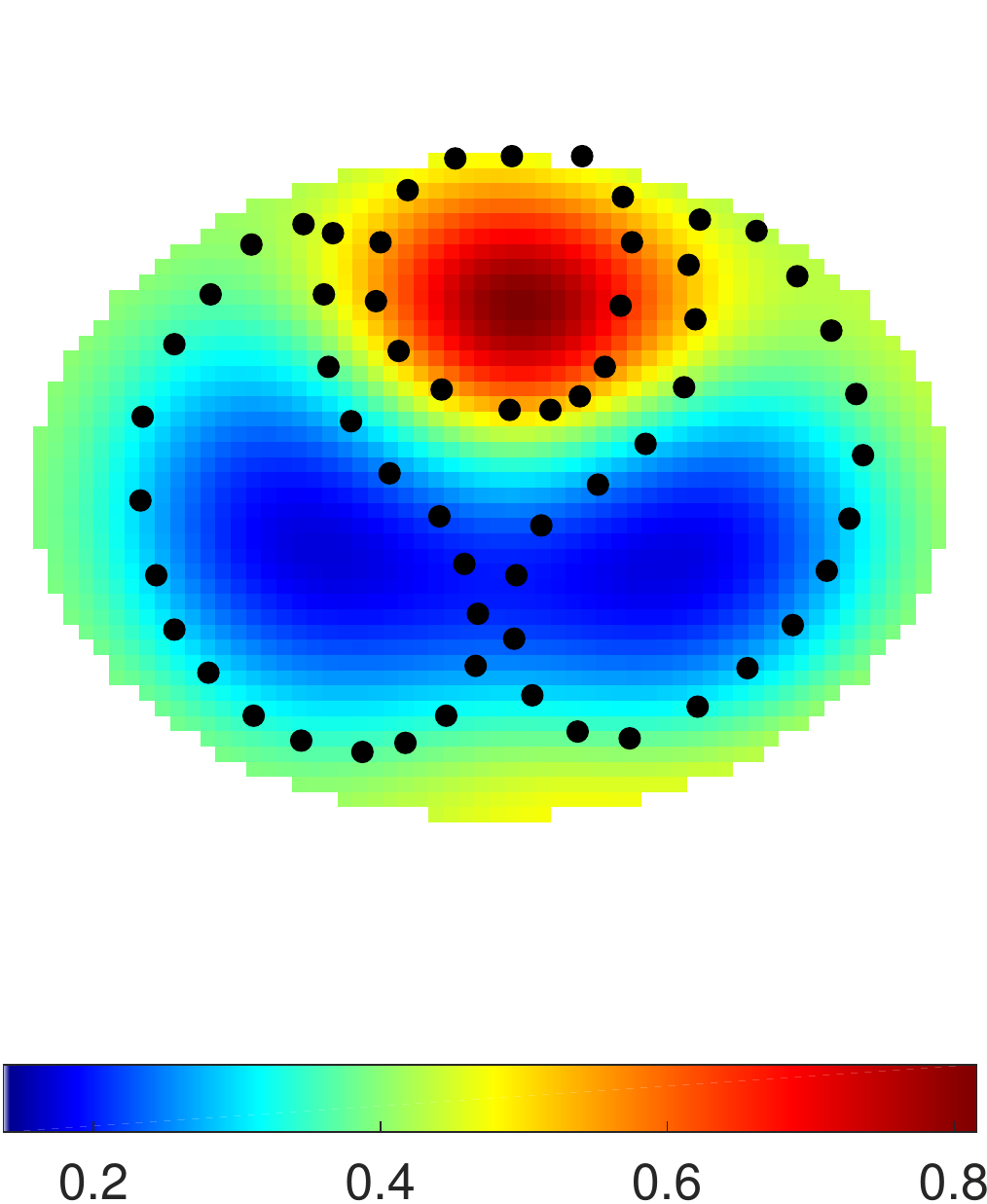}} % Incorrect electrode angles with correct boundary
\put(220,105){\includegraphics[width=80pt]{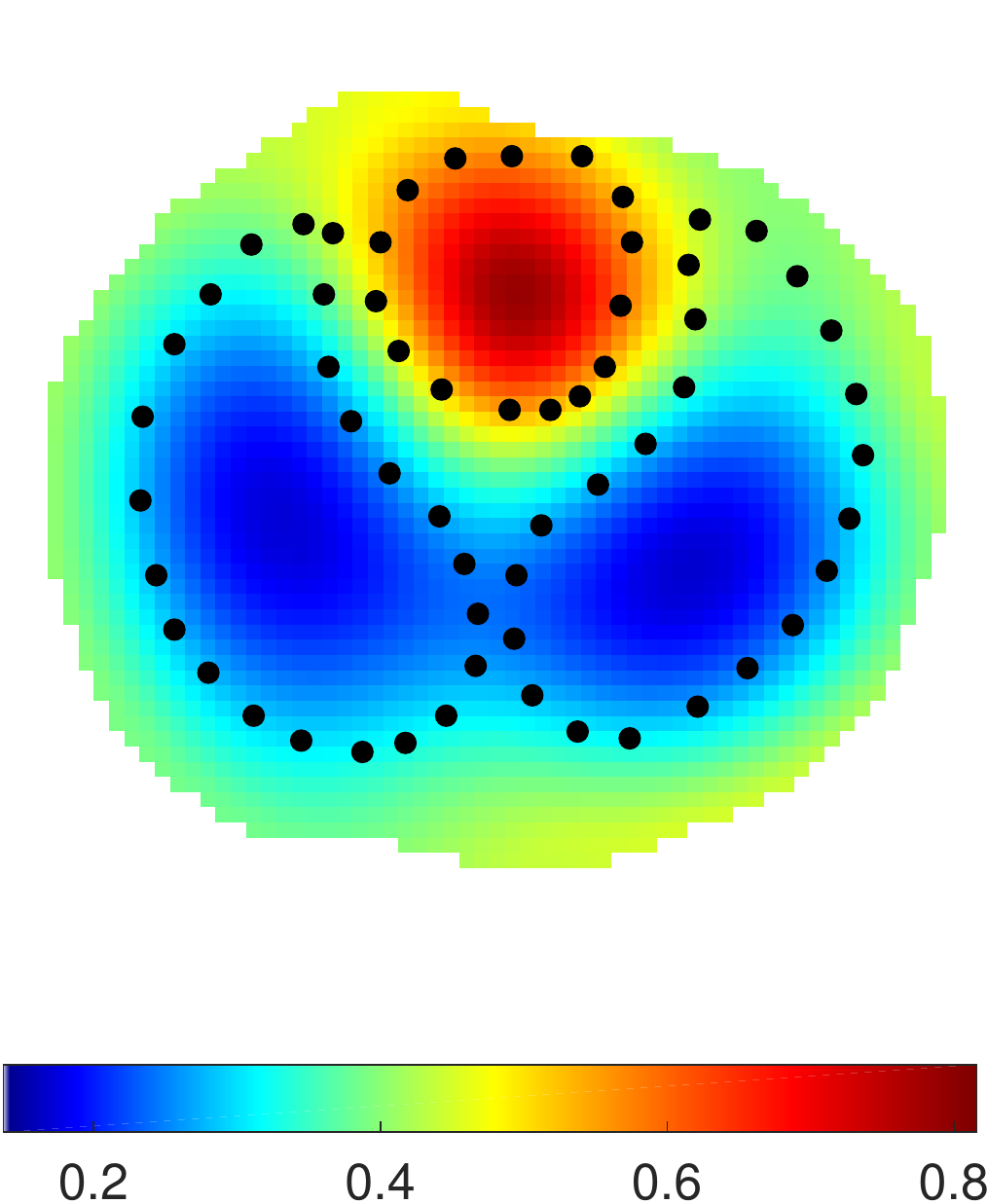}} % True angles and oval1 boundary
\put(320,105){\includegraphics[width=80pt]{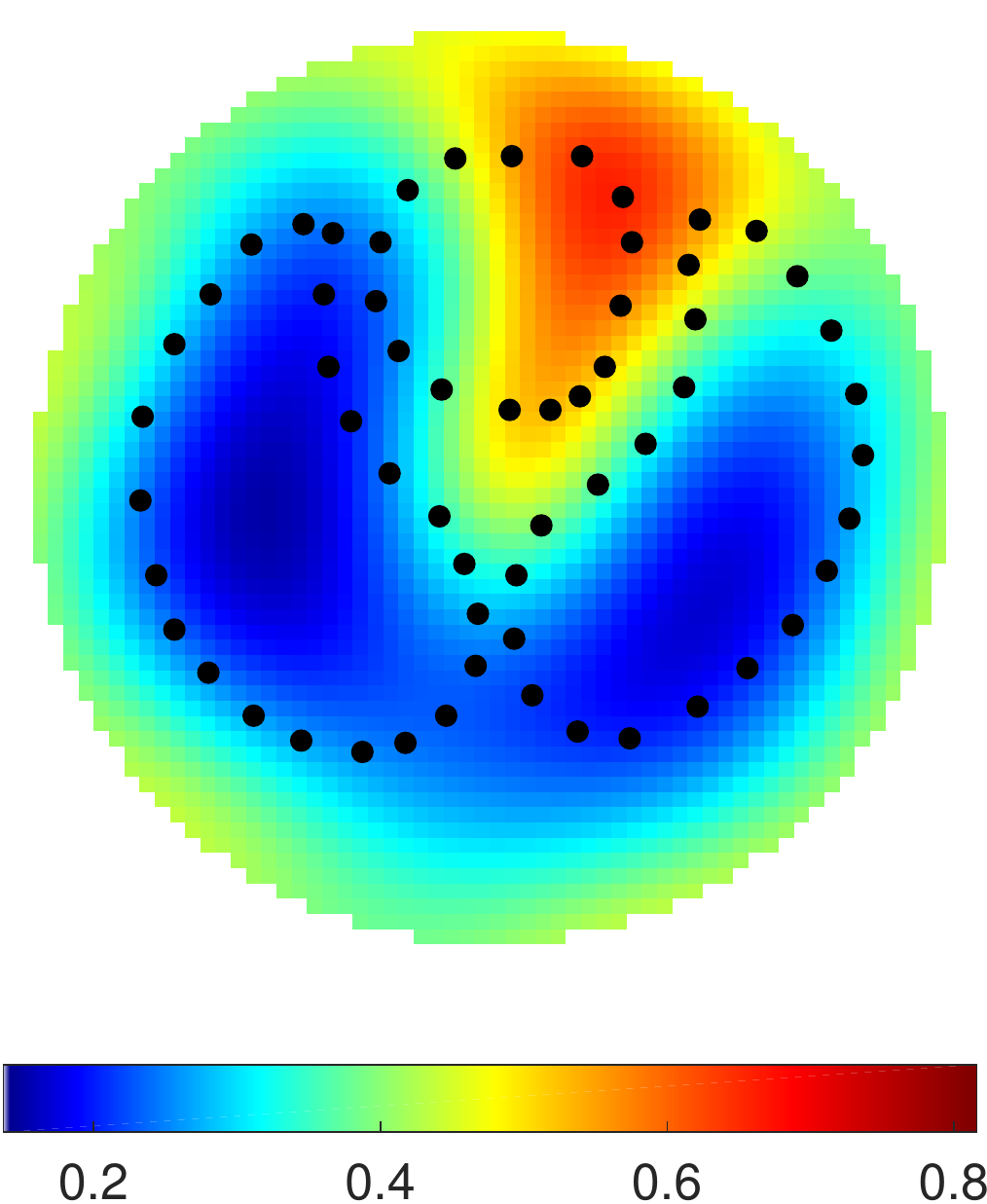}} % True angles and oval1 boundary

% side labels 
\put(0,30){\rotatebox{90}{\scriptsize{\sc Difference}}}
\put(8,42){\rotatebox{90}{\scriptsize{\sc Images}}}

\put(0,135){\rotatebox{90}{\scriptsize{\sc Absolute}}}
\put(8,145){\rotatebox{90}{\scriptsize{\sc Images}}}

%% top labels:
%\put(45,215){\footnotesize {\sc True Angles}}
%\put(35,205){\footnotesize {\sc True Boundary}}
%
%\put(135,215){\footnotesize {\sc True Angles}}
%\put(125,205){\footnotesize {\sc Oval Boundary}}
%
%\put(235,215){\footnotesize {\sc True Angles}}
%\put(225,205){\footnotesize {\sc Clicks Boundary}}
%
%
%%\put(335,215){\footnotesize {\sc True Angles}}
%%\put(325,205){\footnotesize {\sc Clicks Boundary}}
%
%\put(335,215){\footnotesize {\sc Noisy Angles}}
%\put(325,205){\footnotesize {\sc True Boundary}}
%%\put(238,215){\footnotesize {\sc Med Prior}}
%%\put(245,205){\scriptsize $(6,\frac23)$}
%
%%\put(330,215){\footnotesize \textsc{Strong Prior}}
%%\put(345,205){\scriptsize $(8,\frac13)$}

% top labels:
\put(35,215){\footnotesize {\sc True Angles}}
\put(30,205){\footnotesize {\sc True Boundary}}

\put(135,215){\footnotesize {\sc True Angles}}
\put(125,205){\footnotesize {\sc Oval Boundary}}

\put(235,215){\footnotesize {\sc True Angles}}
\put(215,205){\footnotesize {\sc Alternative Boundary}}

\put(330,215){\footnotesize {\sc Noisy Angles}}
\put(327,205){\footnotesize {\sc True Boundary}}

\end{picture}
\caption{\label{fig:ACT3_HnL_FranNata} Comparison of conductivity \trev{(S/m)} reconstructions from the ACT3 data (see Figure~\ref{fig:phantoms}, first) using {\it `Approach~2'} of Section \ref{sec:FranDbar}  for knowledge of true vs. incorrect electrode angles as well as boundary shape.  Absolute images requiring no $\Lambda_1$ are in row 1\tRev{, plotted on the same color scale}.  Difference images are in row 2\tRev{, plotted on the same color scale.}}

\end{figure}
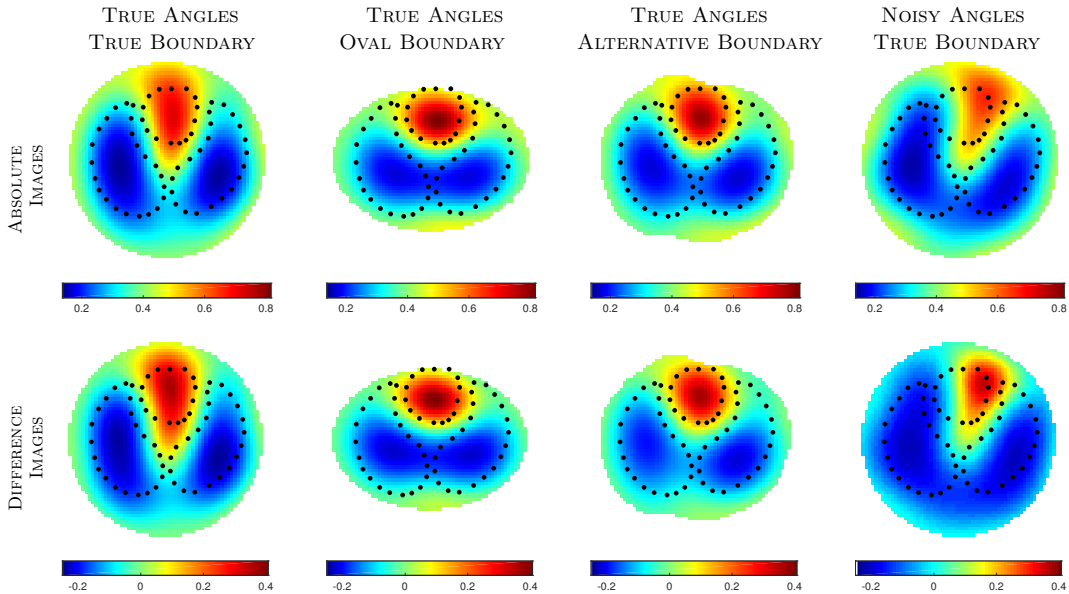

%--------------------------------------------------------
Absolute and time-difference conductivity reconstructions using the three methods are presented in Figures~\ref{fig:ACT3_HnL_Nach}, \ref{fig:ACT3_HnL_FranSarah}, and \ref{fig:ACT3_HnL_FranNata}.  \trev{Note that only Figure~\ref{fig:ACT3_HnL_FranNata} tackles the noisy electrodes with correct boundary shape case.}  This is due to the fact the electrode center locations are so poorly guessed that the resulting electrode edges overlap \trev{making simulation of $\mathbf{L}_1$, required for the $\texp$ and {\it `Approach~1'} methods, impractical.  For the remaining cases, the DN data $\mathbf{L}_1$ was formed using the simulated current and voltage data from solving \eqref{eq:admitt} with $\sigma=1$, subject to boundary conditions defined by the Complete Electrode Model [\cite{Somersalo1992}], with constant non-optimized effective contact impedances of \tRev{0.00000024~$\Omega$m}, using the Finite Element Method for 1) the true boundary and electrode angles (4,227 elements), 2) ovular boundary and correct electrode angles (4,339 elements), and 3) alternative boundary shape and true electrode angles (4,122 elements), with with $A=1$ milliamp amplitude trigonometric current patterns \eqref{eq:trigCPs}.  While difference images with correct boundary shape and incorrect electrode locations were computed for $\tdiff$ and $\Spsidiff$, as no $\mathbf{L}_1$ is needed, the images were very similar to those of {\it `Approach~2'}  and therefore omitted for brevity.  The} highly imprecise electrode location case was included to demonstrate the robustness of {\it `Approach~2'}.  \trev{The reference data for the difference images contained only saline with conductivity $0.424$ S/m.}  Average, max, and min regional pixel values for the heart, left lung, and right lung inclusions for the three D-bar methods are reported in Table~\ref{table:ACT3_Nach}.

%%%%%%%%%%%%%%%%%%%%%%%%%%%%%%%%%%%%%%%%%%%%%%%%%
% ACT3 HnL - Healthy  - DN6o0_S1o0 - Conductivity only - tEXP
%%%%%%%%%%%%%%%%%%%%%%%%%%%%%%%%%%%%%%%%%%%%%%%%%
%%%%%%%%%%%%%%%%%
%%%%%%%%%%%%%%%%%%%%%%%%%%%%%%%%%%%%%%%%%%%%%%%%%
% ACT3-Nach metrics:
%%%%%%%%%%%%%%%%%%%%%%%%%%%%%%%%%%%%%%%%%%%%%%%%%
\begin{table}[h!]

\scriptsize\centering
\caption{Max, min, and average conductivity values \trev{(S/m)} in each of the organ regions in the absolute reconstructions from the ACT3 data (see Figure~\ref{fig:phantoms}, first) using the $\texp$ approximation in Section \ref{sec:NachDbar}, as well as {\it `Approach 1'} \trev{and {\it `Approach 2'} of Section \ref{sec:FranDbar}.}}
\label{table:ACT3_Nach}
\begin{tabular}{|c|ll|ccc|ccc|ccc|ccc|}
\hline
\multicolumn{3}{|c}{\multirow{2}{*}{ }}  &  \multicolumn{3}{|c|}{Correct Boundary} &  \multicolumn{3}{c|}{Oval Boundary} &  \multicolumn{3}{c|}{Alternative Boundary} &  \multicolumn{3}{c|}{Noisy Angles}\\
%&  &  \multicolumn{3}{|c|}{and Elec Angles} &  \multicolumn{3}{c|}{with True  Elec Angles} &  \multicolumn{3}{c|}{and True Elec Angles} \\
\hline
%\cline{2}
Method&Organ & True & AVG &	MAX	& MIN	&	AVG &	MAX	& MIN	&	AVG &	MAX	& MIN&	AVG &	MAX	& MIN	\\
\hline
\multirow{4}{*}{$\texp$} &
Heart  & 0.75 &	0.61 &	0.74 &	0.43 & 0.58  & 0.72 & 	0.41 &		0.60 &	0.73 &	0.45&&&\\
&Left lung & 0.24&0.23 &	0.38 &	0.15 &	0.23 &	0.40 &	0.16 &		0.22 &	0.36 &	0.16 &&&\\
&Right Lung & 0.24 &	0.22 &	0.36 &	0.16 &	0.25 &	0.42 &	0.18 &	0.23 &	0.41 &	0.18 &&&\\
 \cline{2-15}
& \multicolumn{2}{|c|}{\trev{Dynamic Range}} & \multicolumn{3}{|c|}{117\%}& \multicolumn{3}{|c|}{\trev{117\%}}& \multicolumn{3}{|c|}{\trev{117\%}} &&&\\
\hline
\multirow{3}{*}{\it Approach 1} &Heart  & 0.75 &	0.66 &	0.78 &	0.48 &	0.64 &	0.80 &	0.45 &	0.68 &	0.83 &	0.52  &&&\\
&Left lung & 0.24 & 	0.30 &	0.45 &	0.21 &	0.3 & 0.52 &	0.24 &	0.32 &	0.49 &	0.25  &&&\\
&Right Lung &0.24 &	0.28 &	0.42 &	0.21 &	0.31 &	0.46 &	0.25 &	0.31 &	0.47 &	0.25 &&&\\
\cline{2-15}
& \multicolumn{2}{|c|}{\trev{Dynamic Range}} & \multicolumn{3}{|c|}{\trev{111\%}}& \multicolumn{3}{|c|}{\trev{110\%}}& \multicolumn{3}{|c|}{\trev{114\%}} &&&\\
\hline
\multirow{3}{*}{\it Approach 2} &Heart  & 0.75 &	0.61 &	0.72 &	0.39 &		0.65 &	0.82 &	0.47 &		0.65 & 	0.79 & 	0.49 & 0.49 & 0.66 & 0.25  \\
&Left lung & 0.24 & 	0.24 &	0.42 &	0.15 &		0.28 &	0.52 &	0.18 &		0.26 &	0.45 &	0.17 & 0.26 &	0.54 &	0.17 \\
&Right Lung  &0.24 &	0.22 &	0.38 &	0.15 &		0.26 &	0.47 &	0.17 &		0.24 &	0.44 &	0.17 &	0.21 &	0.34 &	0.15 \\
\cline{2-15}
& \multicolumn{2}{|c|}{\trev{Dynamic Range}} & \multicolumn{3}{|c|}{\trev{112\%}}& \multicolumn{3}{|c|}{\trev{126\%}}& \multicolumn{3}{|c|}{\trev{122\%}} &\multicolumn{3}{|c|}{\trev{100\%}} \\
\hline
\hline
\end{tabular}
\end{table}

%%%%%%%%%%%%%%%%%%%%%%%%%%%%%%%%%%%%%%%%%%%%%%%%%
% ACT3 HnL - Healthy  - DN6o0_S1o0 - Conductivity only - Sarah Scat
%%%%%%%%%%%%%%%%%%%%%%%%%%%%%%%%%%%%%%%%%%%%%%%%%
%%%%%%%%%%%%%%%%%

%%%%%%%%%%%%%%%%%
%\clearpage

Figures~\ref{fig:ACT3_HnL_Nach}, \ref{fig:ACT3_HnL_FranSarah}, and \ref{fig:ACT3_HnL_FranNata} demonstrate that the three D-bar methods of Section~\ref{sec:methods} produce absolute and difference images quite similar in resolution and reconstructed conductivity values.   In each case, approximating the circular tank by an ovular boundary has the effect of compressing the reconstructed heart and lungs.  Perturbing the electrode positions has a rotating effect on the images, but still yields easily recognizable reconstructions of the targets.   For uniform perturbations of the electrode angles\trev{,} the resulting images appear uniformly rotated.  

The results in Table~\ref{table:ACT3_Nach} from the $\texp$ approximation show that reconstructions using the correct boundary shape have a maximum value in the heart region within 1.5\% of the true value, while the minimum values in the lung region are 37.5\% and 33.3\% lower than the true values.  However, the average value in each lung was within 4.2\% of the true value.  On the oval and alternative boundaries, the conductivity values remained very close to the reconstructed values on the correct boundary shape, with a very slight increase in error.  The dynamic ranges for the reconstructions on the true boundary, oval boundary, and alternative boundary were 117\%, 108\%, and 112\%, respectively.

The results in Table~\ref{table:ACT3_Nach} from {\it `Approach 1'}, Section~\ref{sec:FranDbar}\trev{,} show that reconstructions using the correct boundary shape have a maximum value in the heart region within 4\% of the true value, while the minimum values in the lung region are 12.5\% lower than the true values.  The average value in each lung was 25\% and 16.7\% above the true value.  Errors were slightly larger on the oval and alternative boundaries.  The dynamic ranges for the reconstructions on the true boundary, oval boundary, and alternative boundary were 111\%, 110\%, and 114\%, respectively.

The results in Table~\ref{table:ACT3_Nach} from {\it `Approach 2'}, Section~\ref{sec:FranDbar}\trev{,} resulted in reconstructions with a maximum value in the heart region within 4\% of the true value when using the correct boundary shape, with the minimum values in the lung region 37.5\% lower than the true values.  The average value in each lung was exactly correct in the left lung and within 8.3\% of the true value in the right lung.  Errors were comparable on the oval and alternative boundaries, except slightly larger in the heart region.  The largest errors in the heart region were seen in the case of the incorrect boundary locations (``noisy angles''), which had a 12\% relative error in the reconstructed maximum value from the true value.  The dynamic ranges for the reconstructions on the true boundary, oval boundary, alternative boundary, and correct boundary with incorrect electrode locations were 112\%, 126\%,  122\%, and 100\%, respectively.

%\clearpage
%--------------------------------------------------------
\subsection{Reconstructions from ACT4 Data}\label{sec:results_ACT4}
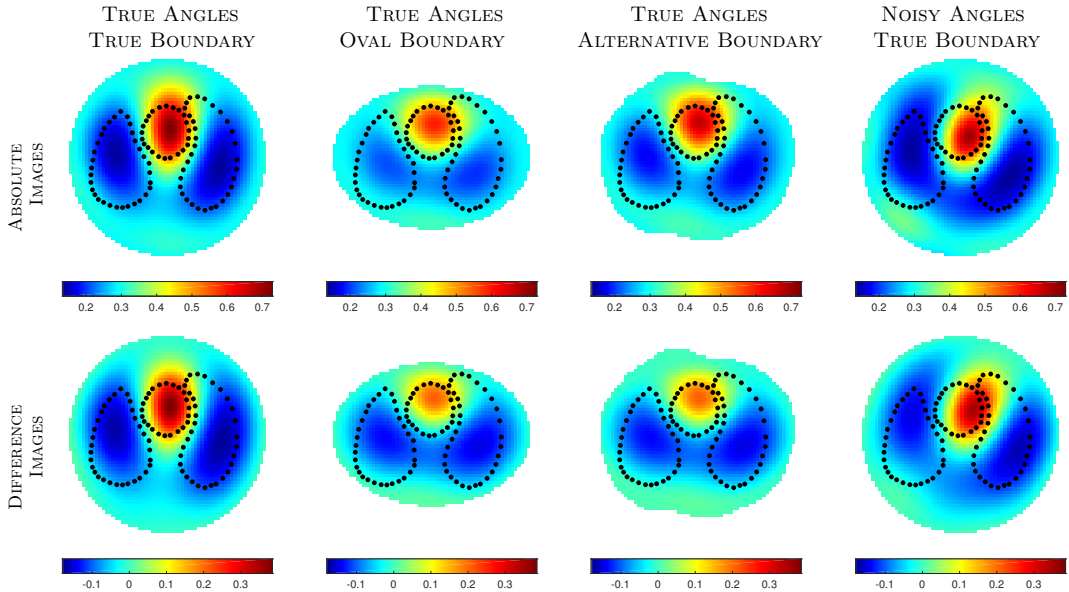
\begin{figure}[!h]

\begin{picture}(420,220)

% Difference images
\put(20,0){\includegraphics[width=80pt]{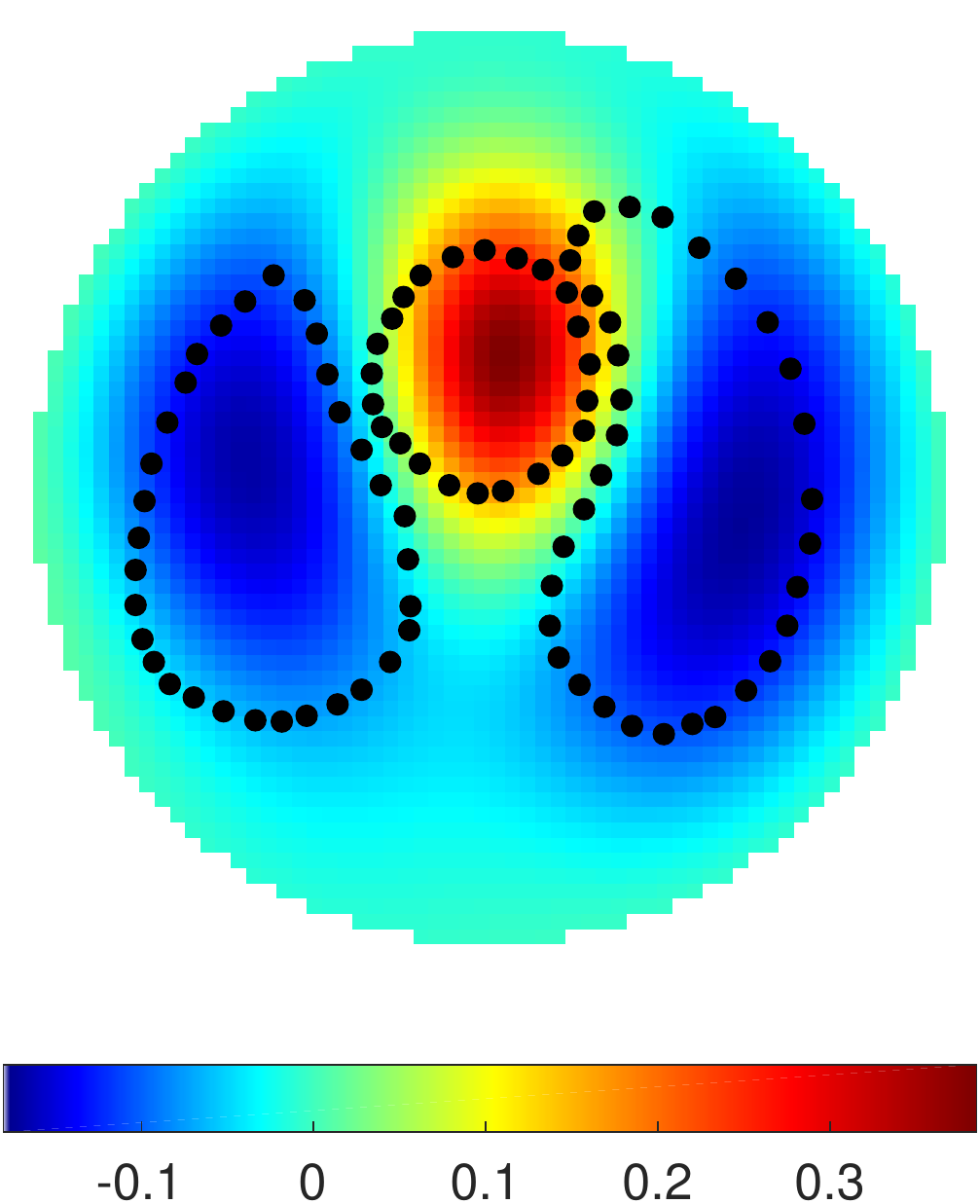}} % True angles and boundary image

\put(120,0){\includegraphics[width=80pt]{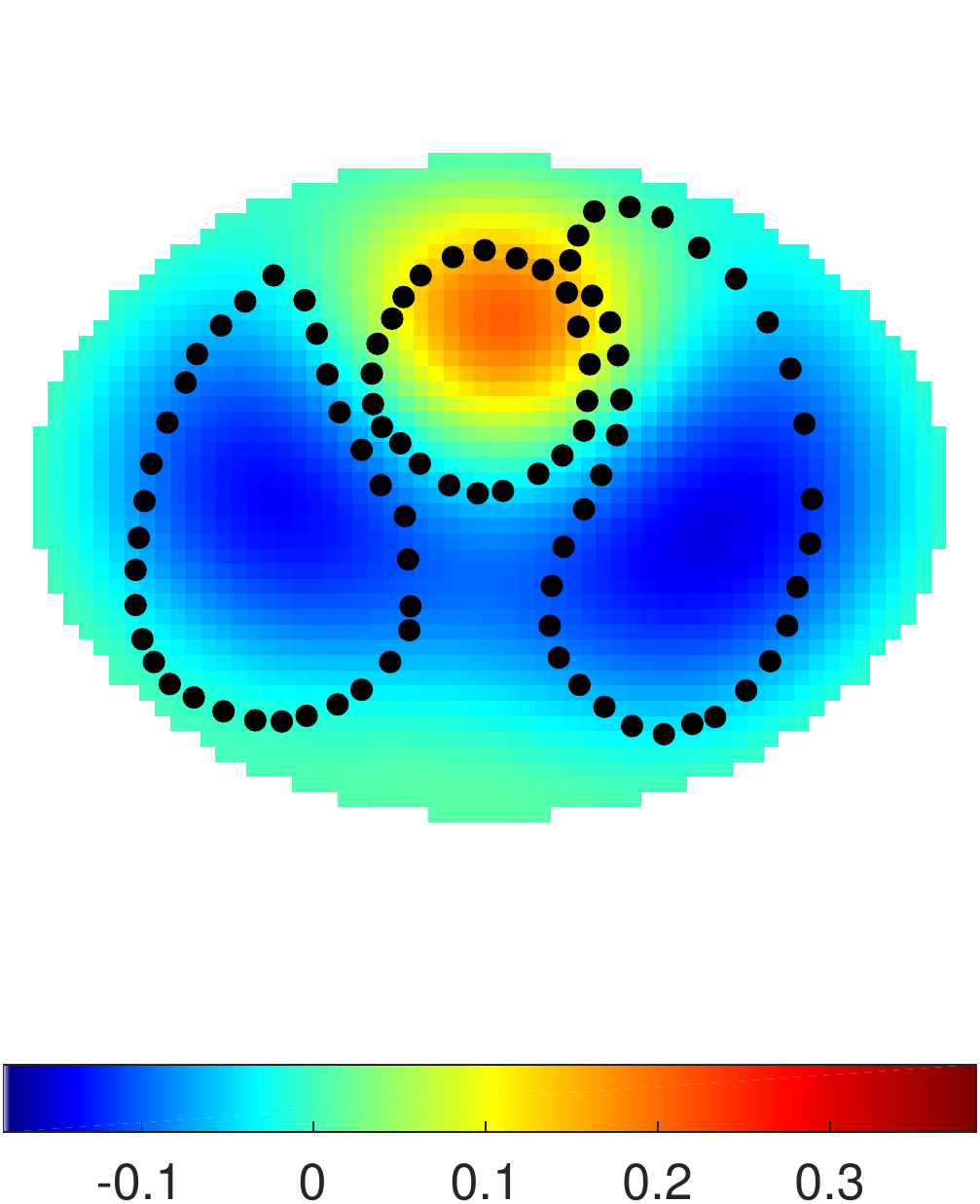}} % True angles and oval1 boundary
\put(220,0){\includegraphics[width=80pt]{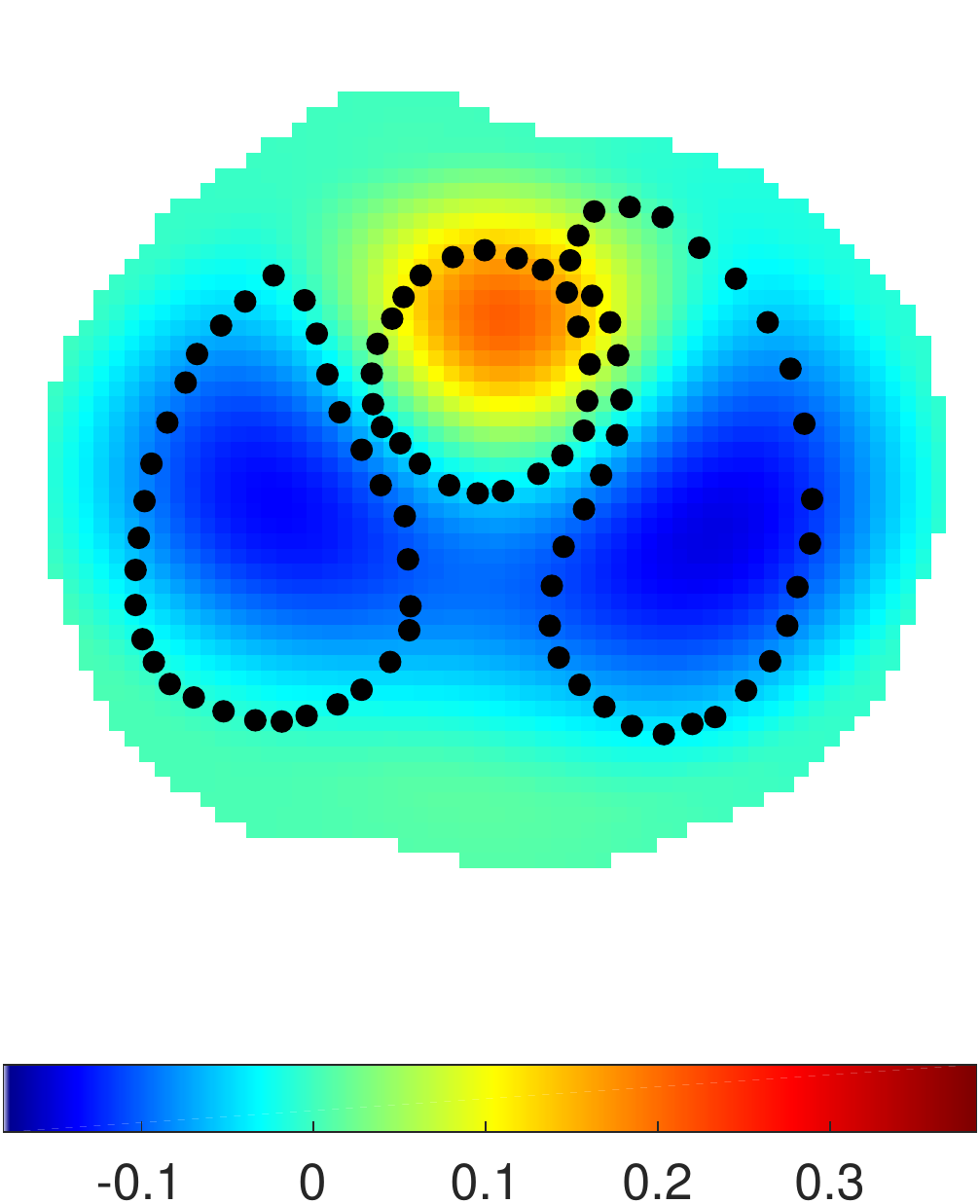}} % True angles and clicks boundary
\put(320,0){\includegraphics[width=80pt]{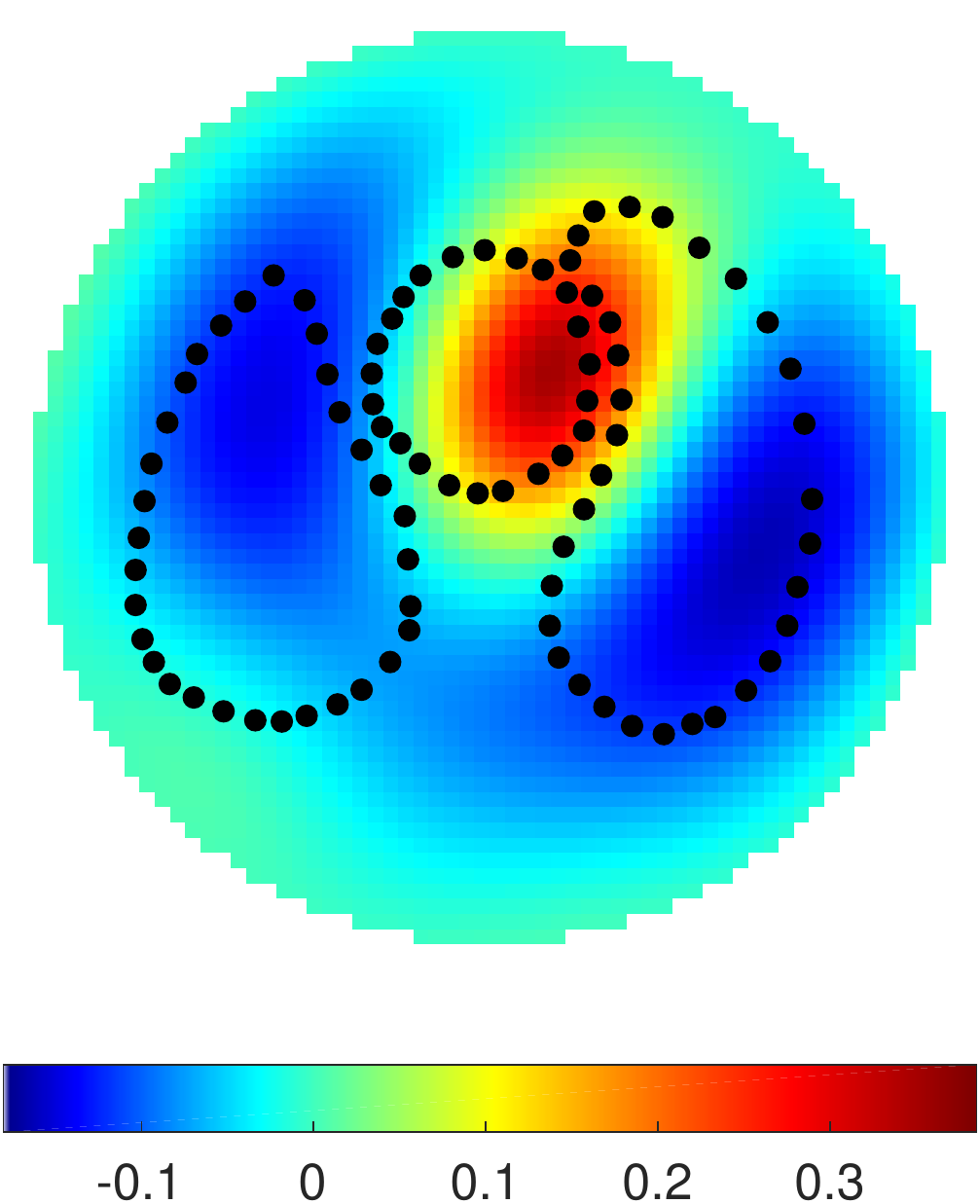}} % Incorrect electrode angles with correct boundary

% Absolute images
\put(20,105){\includegraphics[width=80pt]{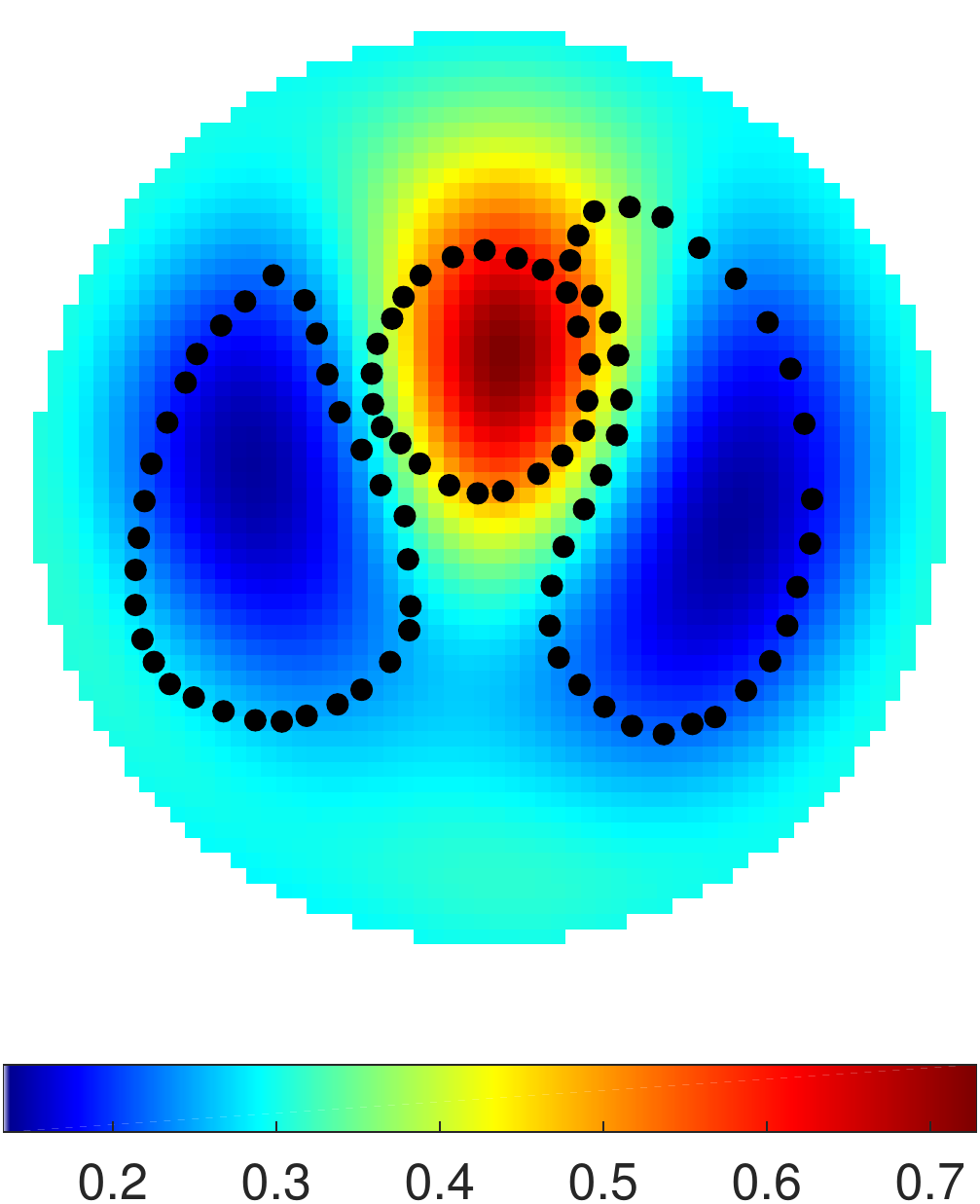}} % True angles and boundary image
\put(120,105){\includegraphics[width=80pt]{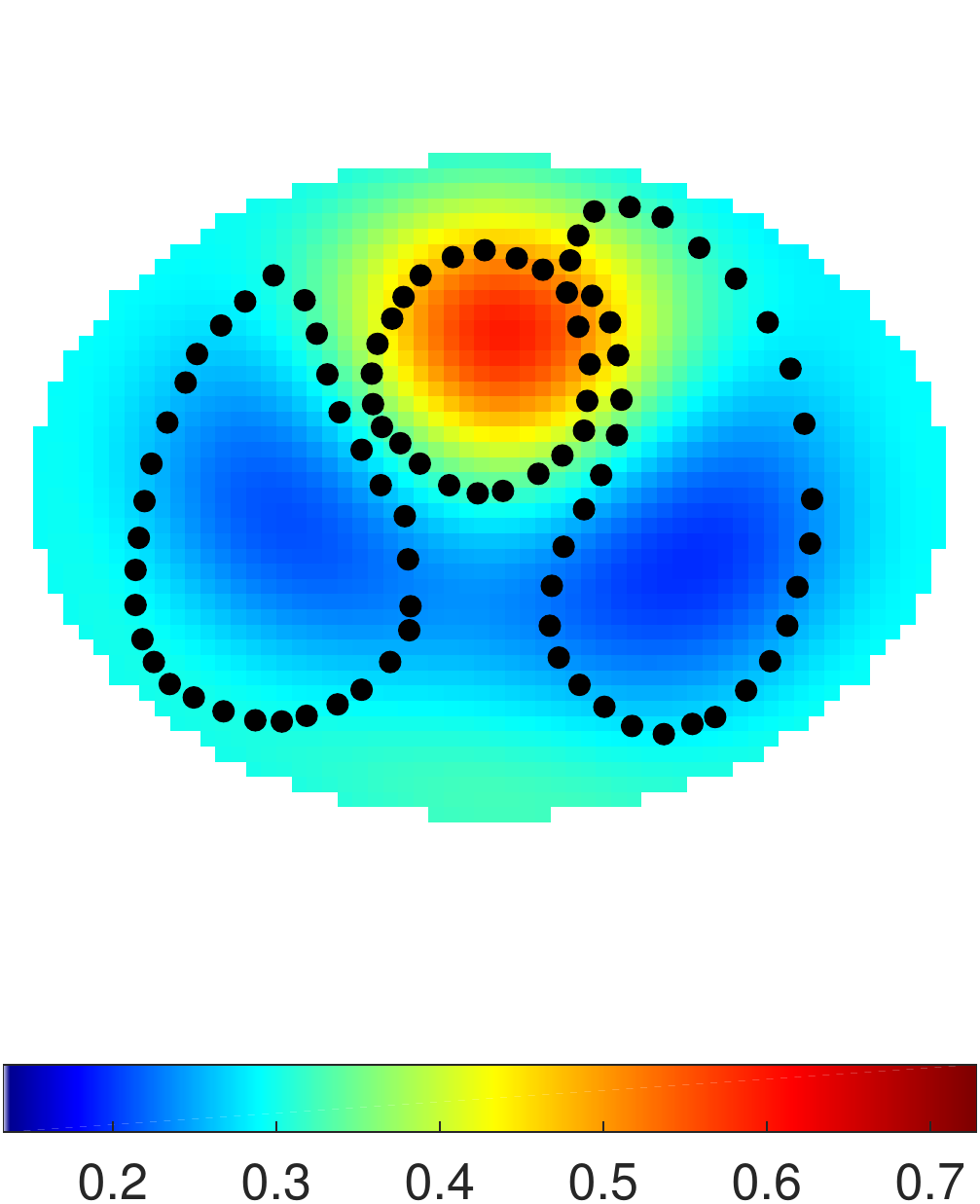}} % True angles and oval1 boundary
\put(220,105){\includegraphics[width=80pt]{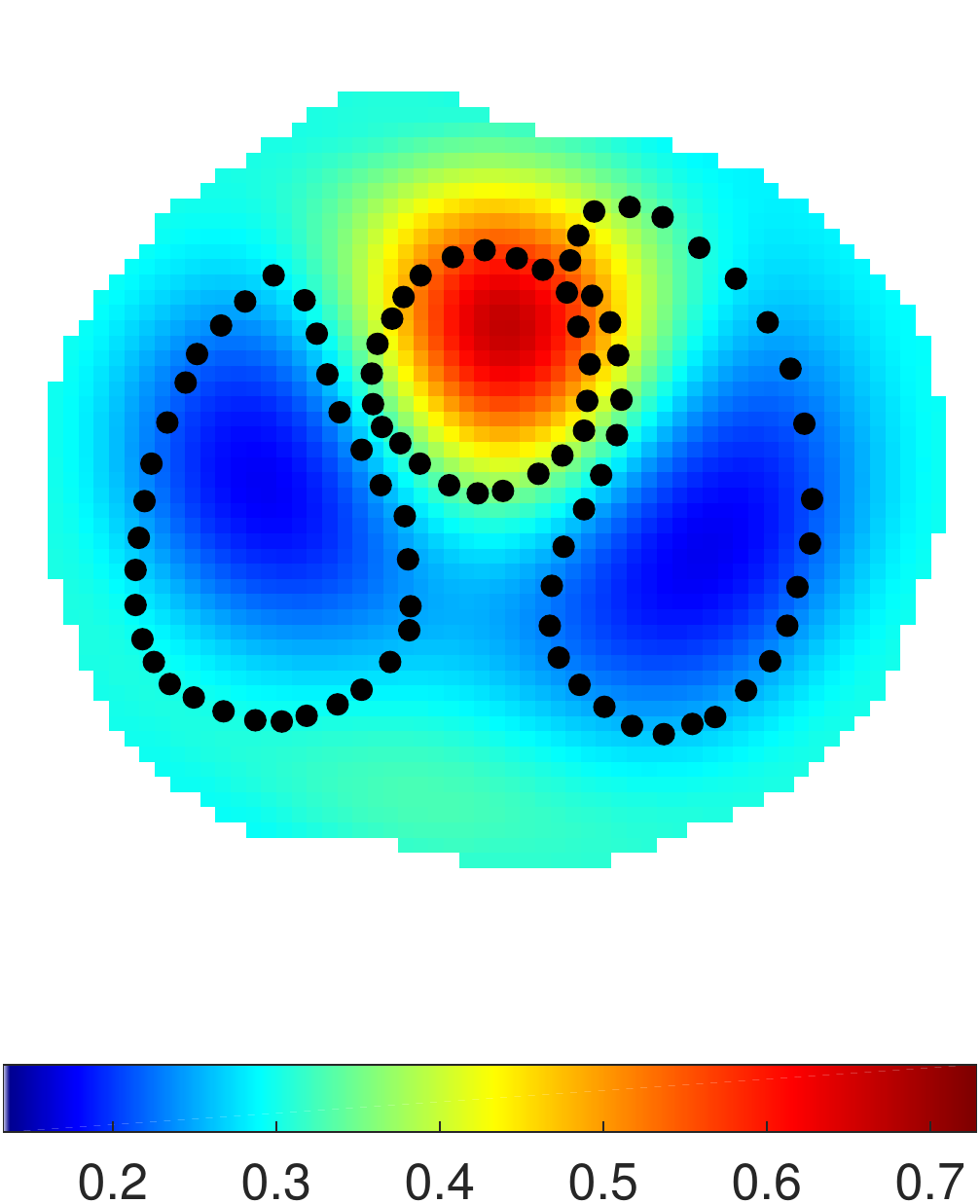}} % True angles and clicks boundary
\put(320,105){\includegraphics[width=80pt]{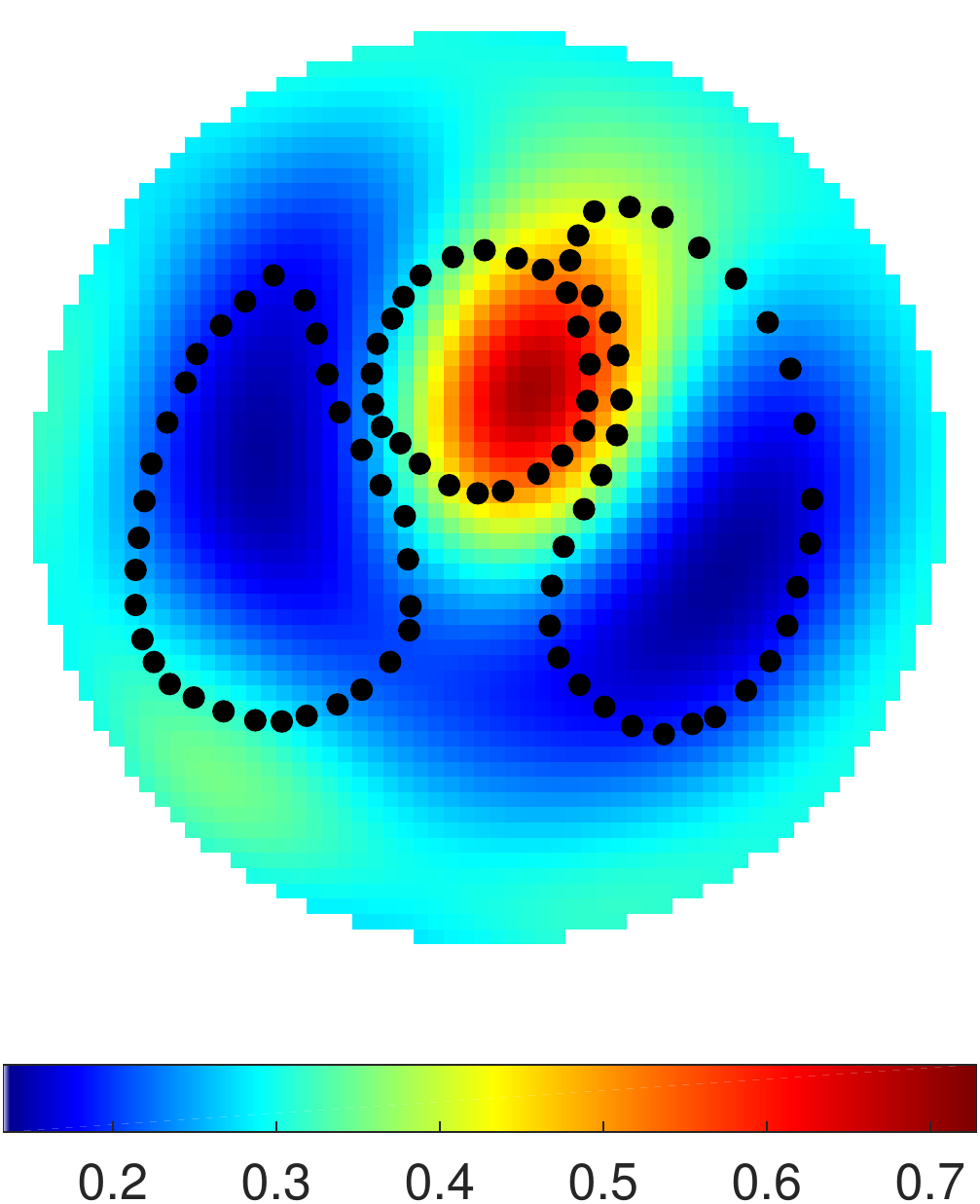}} % Incorrect electrode angles with correct boundary

% side labels 
\put(0,30){\rotatebox{90}{\scriptsize{\sc Difference}}}
\put(8,42){\rotatebox{90}{\scriptsize{\sc Images}}}

\put(0,135){\rotatebox{90}{\scriptsize{\sc Absolute}}}
\put(8,145){\rotatebox{90}{\scriptsize{\sc Images}}}

% top labels:
\put(35,215){\footnotesize {\sc True Angles}}
\put(30,205){\footnotesize {\sc True Boundary}}

\put(135,215){\footnotesize {\sc True Angles}}
\put(125,205){\footnotesize {\sc Oval Boundary}}

\put(235,215){\footnotesize {\sc True Angles}}
\put(215,205){\footnotesize {\sc Alternative Boundary}}

\put(330,215){\footnotesize {\sc Noisy Angles}}
\put(327,205){\footnotesize {\sc True Boundary}}
\end{picture}
\caption{\label{fig:ACT4_HLSA_PE_FranNata} Comparison of conductivity  \trev{(S/m)} reconstructions from ACT4  data with \trev{the healthy heart and lungs} phantom (see Figure~\ref{fig:phantoms}, second).  Results are compared for knowledge of true vs. incorrect electrode angles as well as boundary shape.  Absolute images requiring no $\Lambda_1$ are in row 1\tRev{, plotted on the same color scale}.  Difference images are in row 2\tRev{, plotted on the same color scale.}}
\end{figure}

\begin{figure}[!h]

\begin{picture}(420,220)

% Difference images
\put(20,0){\includegraphics[width=80pt]{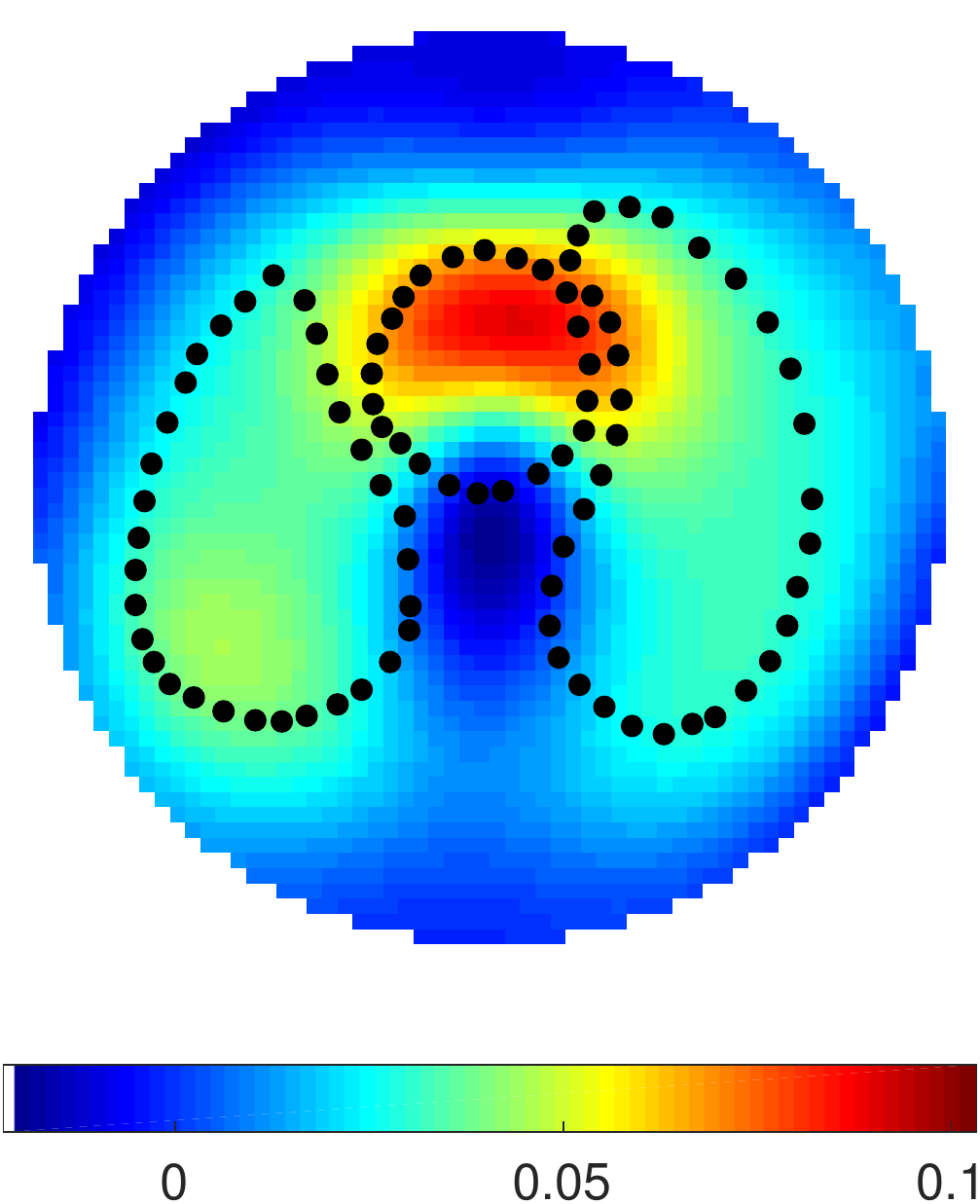}} % True angles and boundary image
\put(120,0){\includegraphics[width=80pt]{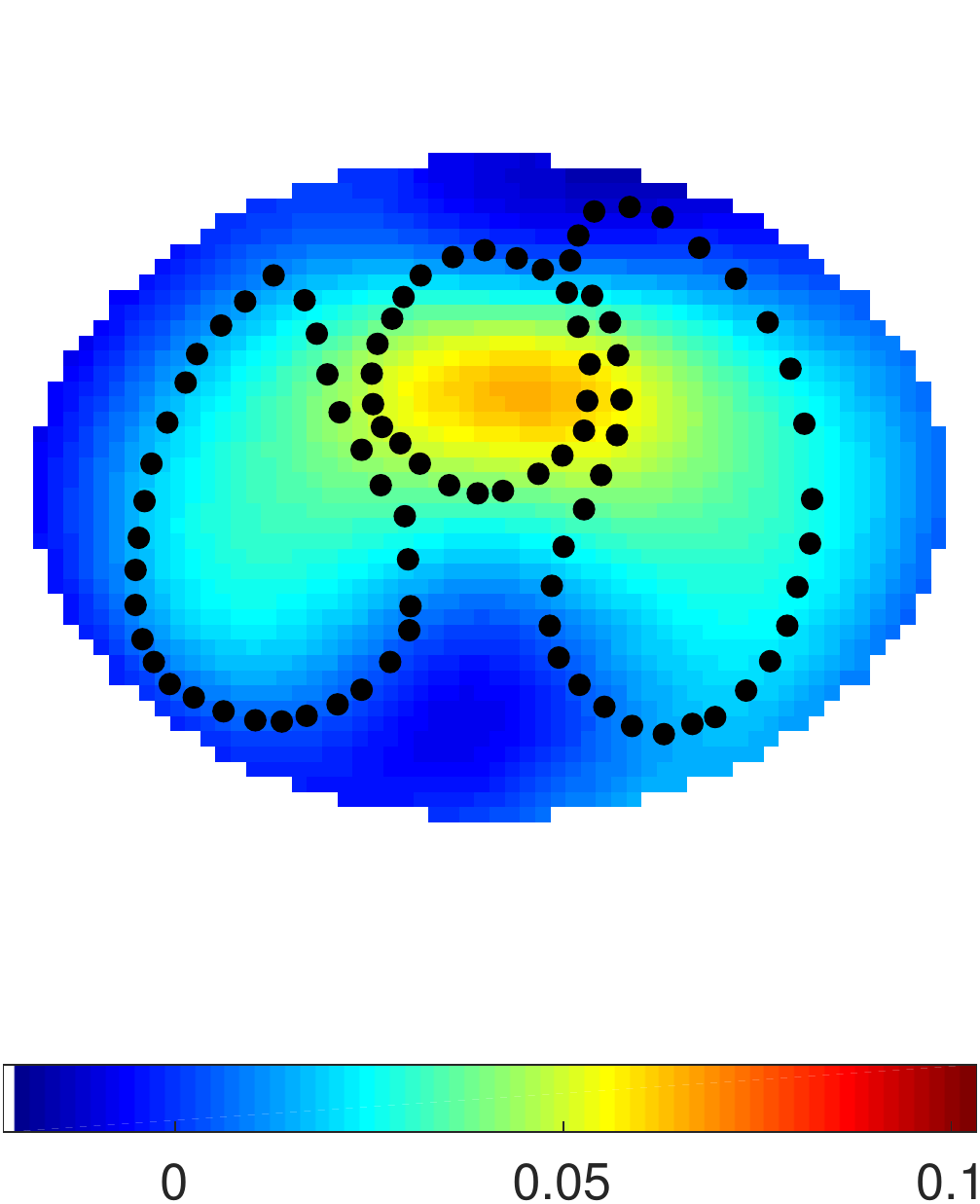}} % True angles and oval1 boundary
\put(220,0){\includegraphics[width=80pt]{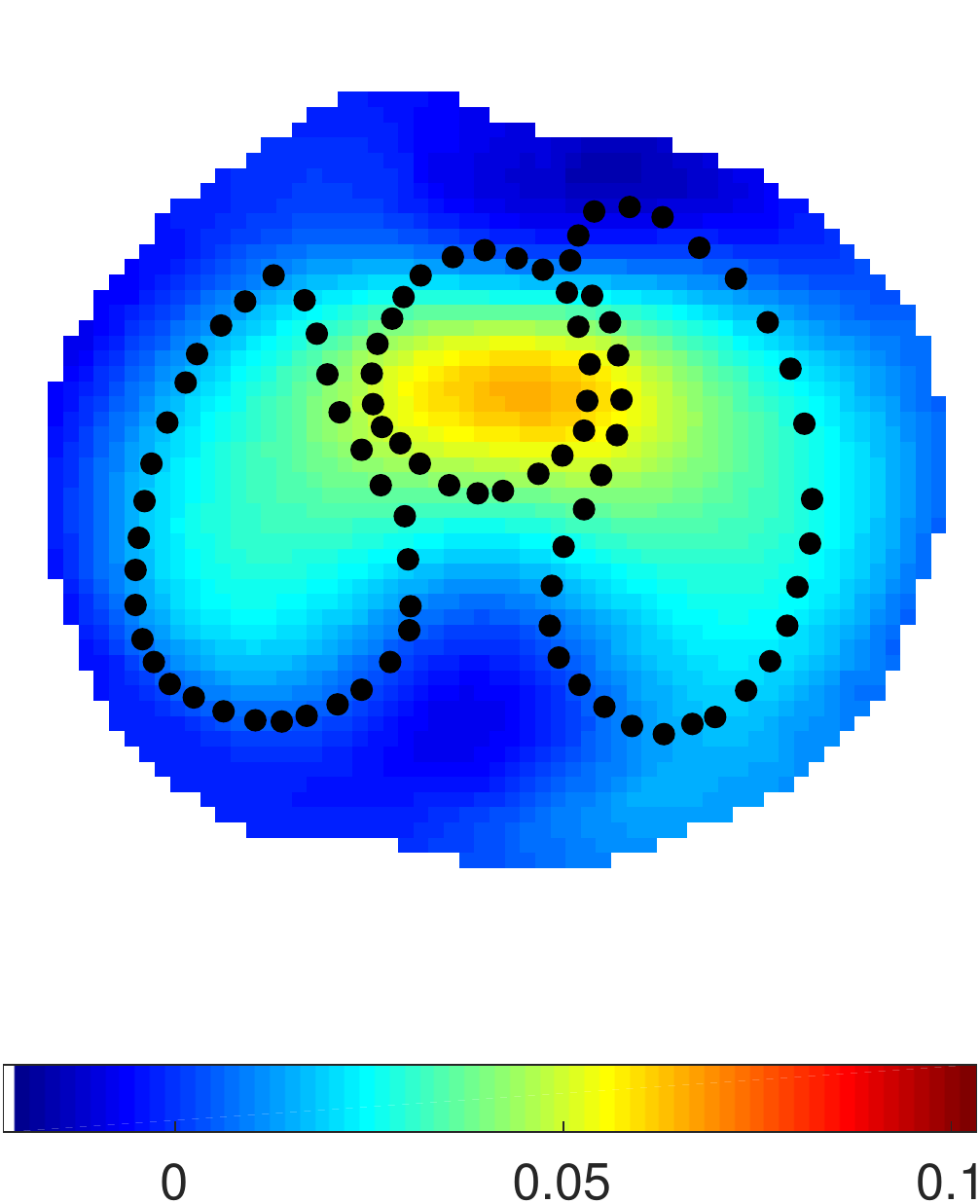}} % True angles and clicks boundary
\put(320,0){\includegraphics[width=80pt]{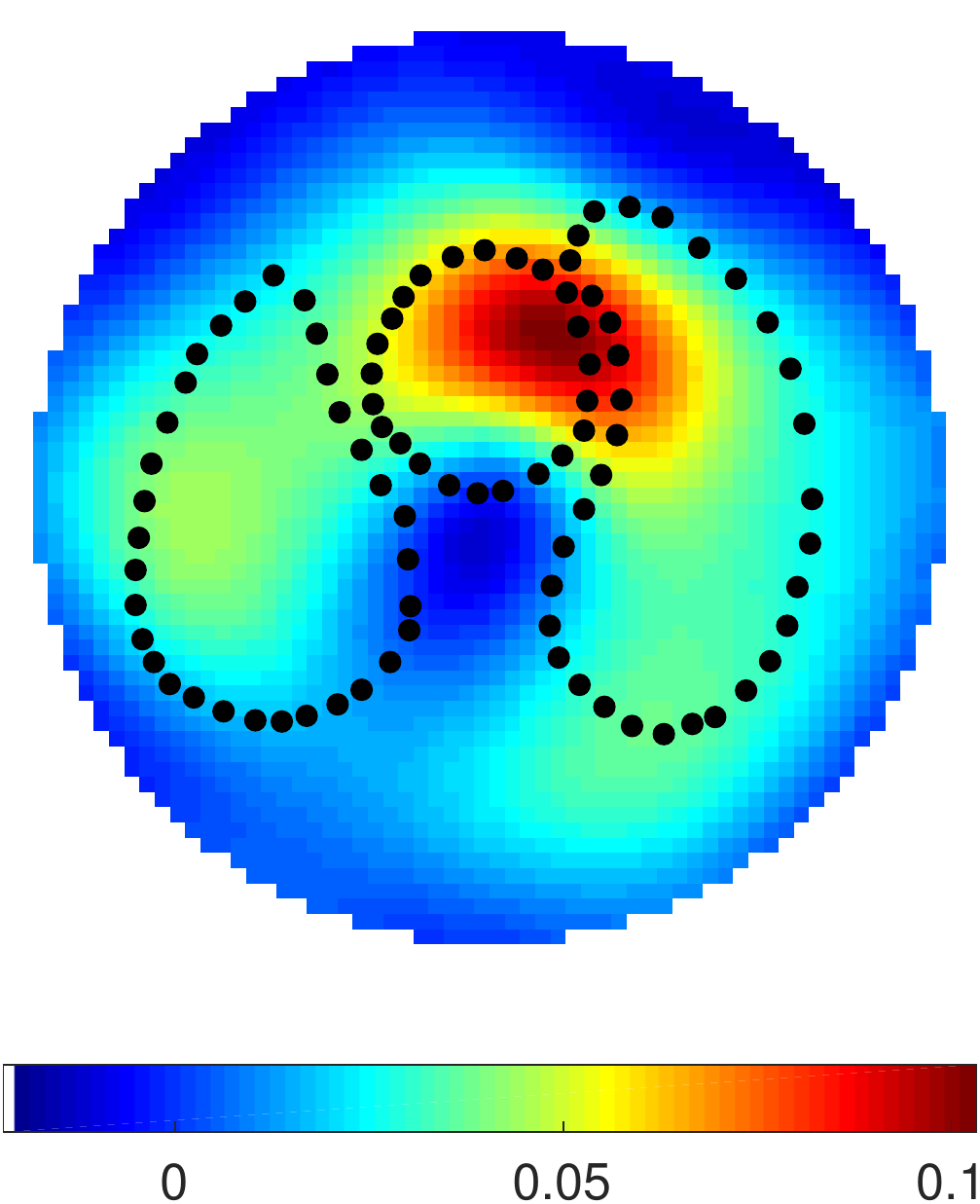}} % Incorrect electrode angles with correct boundary

% Absolute images
\put(20,105){\includegraphics[width=80pt]{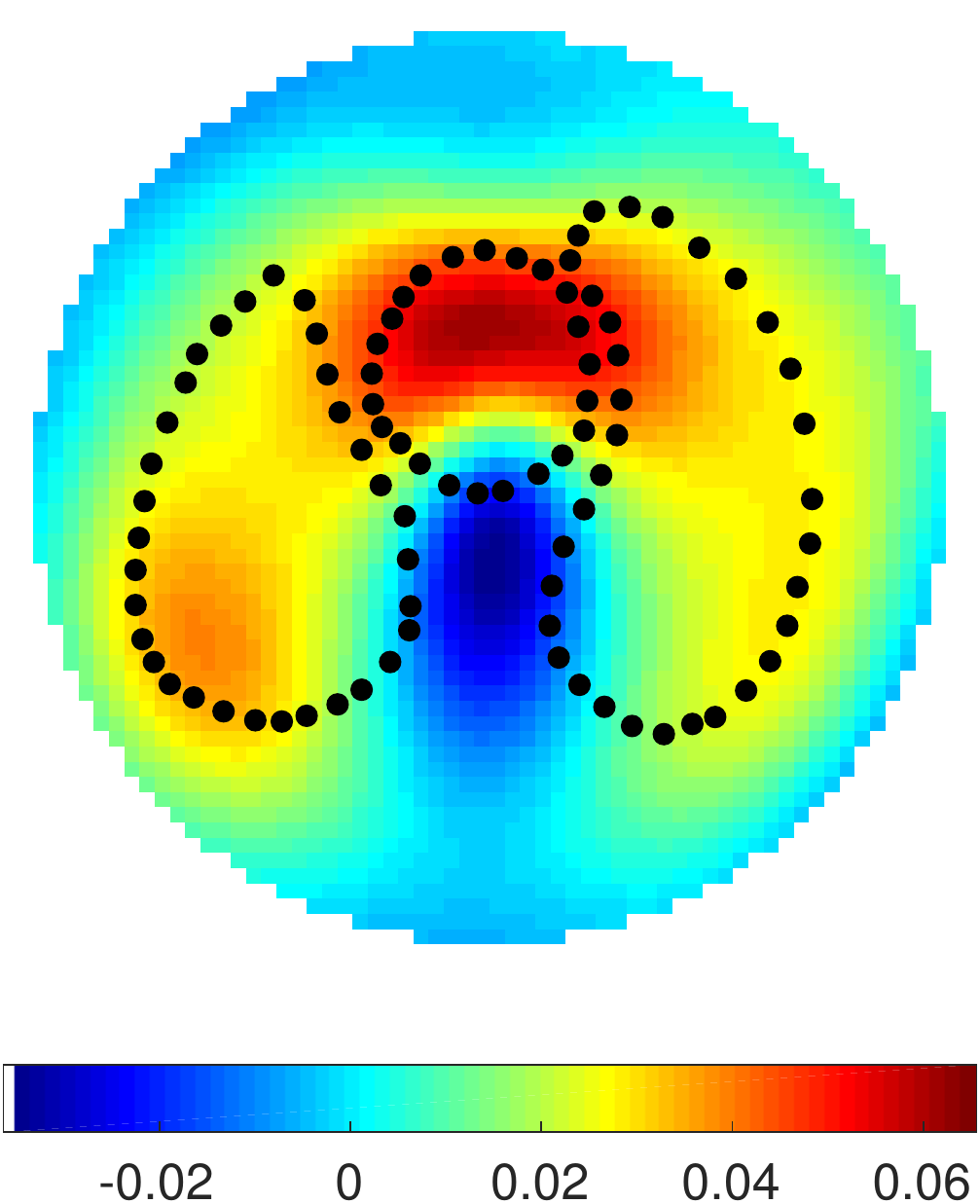}} % True angles and boundary image
\put(120,105){\includegraphics[width=80pt]{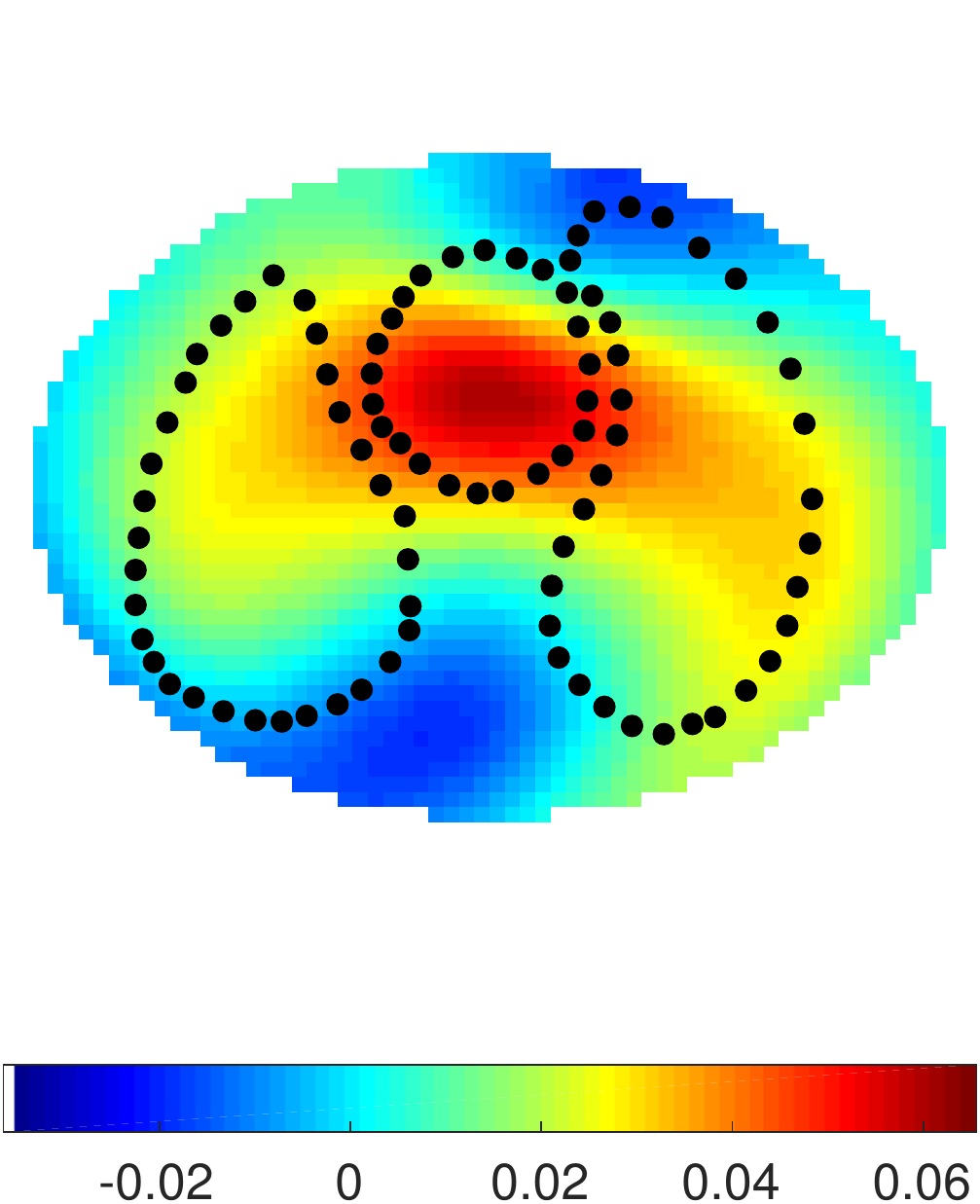}} % True angles and oval1 boundary
\put(220,105){\includegraphics[width=80pt]{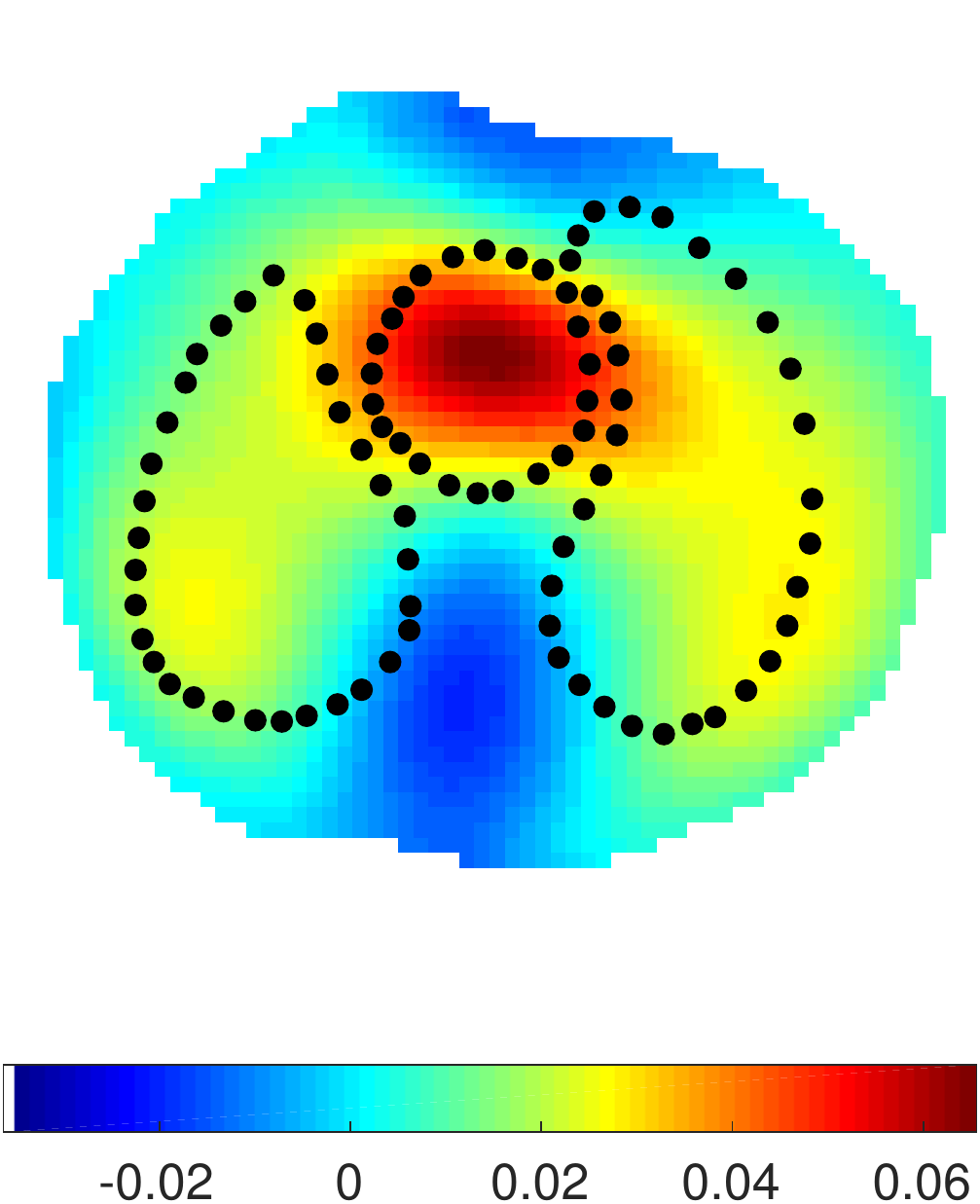}} % True angles and clicks boundary
\put(320,105){\includegraphics[width=80pt]{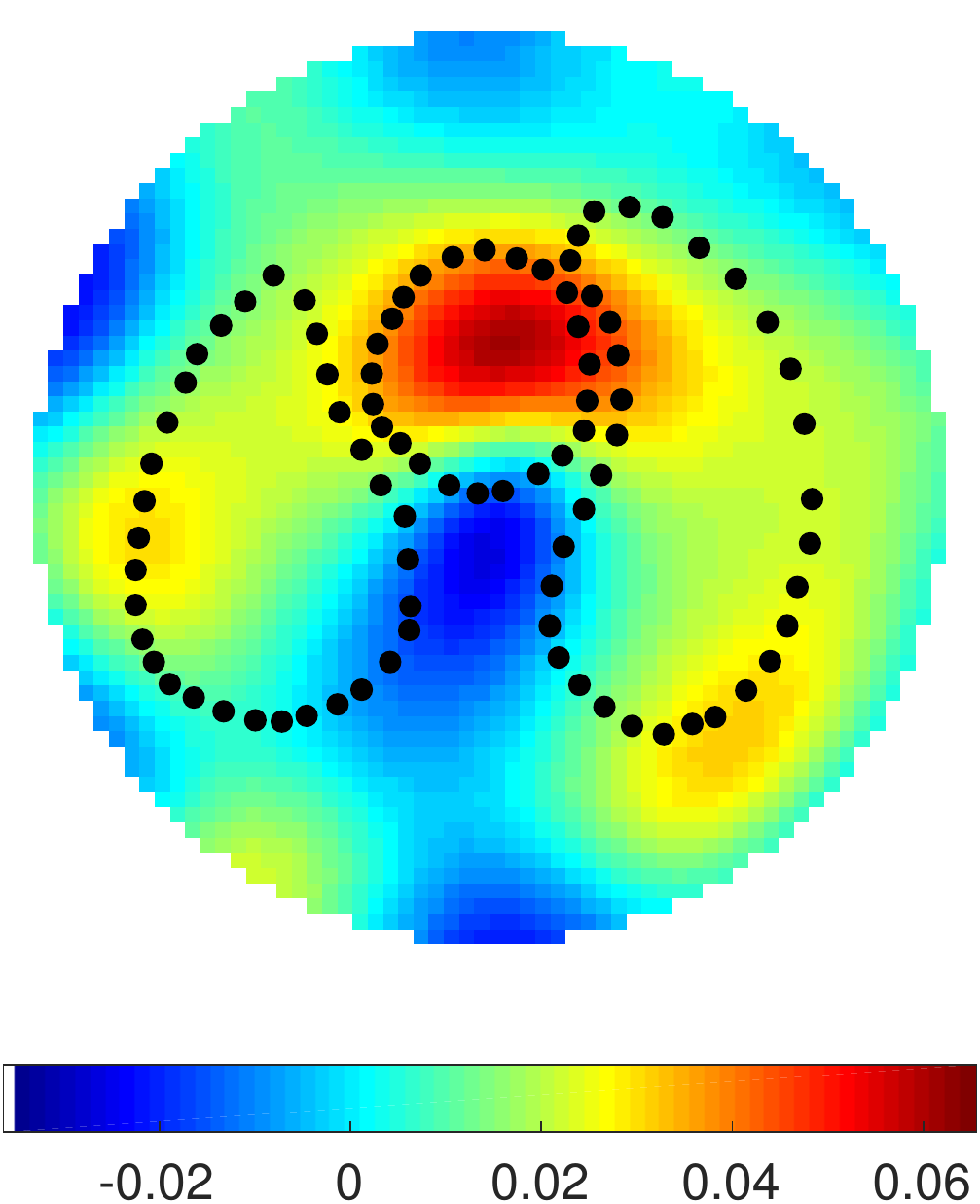}} % Incorrect electrode angles with correct boundary

% side labels 
\put(0,30){\rotatebox{90}{\scriptsize{\sc Difference}}}
\put(8,42){\rotatebox{90}{\scriptsize{\sc Images}}}

\put(0,135){\rotatebox{90}{\scriptsize{\sc Absolute}}}
\put(8,145){\rotatebox{90}{\scriptsize{\sc Images}}}

% top labels:
\put(35,215){\footnotesize {\sc True Angles}}
\put(30,205){\footnotesize {\sc True Boundary}}

\put(135,215){\footnotesize {\sc True Angles}}
\put(125,205){\footnotesize {\sc Oval Boundary}}

\put(235,215){\footnotesize {\sc True Angles}}
\put(215,205){\footnotesize {\sc Alternative Boundary}}

\put(330,215){\footnotesize {\sc Noisy Angles}}
\put(327,205){\footnotesize {\sc True Boundary}}
%\put(238,215){\footnotesize {\sc Med Prior}}
%\put(245,205){\scriptsize $(6,\frac23)$}

%\put(330,215){\footnotesize \textsc{Strong Prior}}
%\put(345,205){\scriptsize $(8,\frac13)$}

\end{picture}
\caption{\label{fig:ACT4_HLSA_PE_FranNata_perm}  Comparison of susceptivity $\omega\epsilon$  \trev{(S/m)} reconstructions from ACT4  data with \trev{the healthy heart and lungs} phantom (see Figure~\ref{fig:phantoms}, second).  Results are compared for knowledge of true vs. incorrect electrode angles as well as boundary shape.  Absolute images requiring no $\Lambda_1$ are in row 1\tRev{, plotted on the same color scale}.  Difference images are in row 2\tRev{, plotted on the same color scale.}}

\end{figure}

%--------------------------------------------------------
Reconstructions for the remainder of the manuscript focus solely on the {\it 'Approach~2'} D-bar method due to it's ability to form absolute images without a need to simulate $\Lambda_1$ measurements, as well as flexibility for conductivity or admittivity imaging.   \trev{The reference data for the difference images contained only saline with conductivity $0.3$ S/m.}  
%The data collected using  ACT4 system contains targets of agar and graphite simulating a phantom chest with heart, lungs, spine, and aorta, where the right (DICOM) lung has a portion (away from the heart) replaced with `heart' agar/graphite phantom to simulate a pleural effusion.  
Conductivity $(\sigma)$ images for the ACT4 data are displayed in Figure~\ref{fig:ACT4_HLSA_PE_FranNata}, while susceptivity ($\omega\epsilon$) images are shown in Figure~\ref{fig:ACT4_HLSA_PE_FranNata_perm}\trev{,} for the true and incorrect boundaries and electrode locations. Regional averages, maxes, and mins are reported in Table~\ref{table:ACT4_FranNata_COND}.

%%%%%%%%%%%%%%%%%
%

%%%%%%%%%%%%%%%%%%%%%%%%%%%%%%%%%%%%%%%%%%%%%%%%%
% ACT4 HL Perm only
%%%%%%%%%%%%%%%%%%%%%%%%%%%%%%%%%%%%%%%%%%%%%%%%%
\begin{table}[h!]

\scriptsize\centering
\caption{Max, min, and average conductivity and susceptivity values  \trev{(S/m)} in each of the organ regions in the absolute reconstructions from the ACT4 data (see Figure~\ref{fig:phantoms}, second) using {\it `Approach~2'} of Section~\ref{sec:FranDbar}.}
\label{table:ACT4_FranNata_COND}
\begin{tabular}{|c|ll|ccc|ccc|ccc|ccc|}
\hline
\multicolumn{3}{|c}{\multirow{2}{*}{ }}  &  \multicolumn{3}{|c|}{Correct Boundary} &  \multicolumn{3}{c|}{Oval Boundary} &  \multicolumn{3}{c|}{Alternative Boundary} &  \multicolumn{3}{c|}{Noisy Angles}\\
%&  &  \multicolumn{3}{|c|}{and Elec Angles} &  \multicolumn{3}{c|}{with True  Elec Angles} &  \multicolumn{3}{c|}{and True Elec Angles} &  \multicolumn{3}{c|}{Noisy Angles}\\
\hline
%\cline{2}
&Organ & True & AVG &	MAX	& MIN	&	AVG &	MAX	& MIN	&	AVG &	MAX	& MIN	&	AVG &	MAX	& MIN	\\
\hline
\multirow{4}{*}{\rotatebox{90}{Cond.} }&Heart  &\trev{0.68} &	\trev{0.58}	&\trev{0.73}	&\trev{0.32}	&\trev{0.48}	&\trev{0.60}	&\trev{0.33}	&\trev{0.52}	&\trev{0.67}	&\trev{0.32}	&\trev{0.51}	&\trev{0.70}	&\trev{0.25}\\
&Left lung & \trev{0.057} & \trev{0.20	} &	\trev{0.30}	&\trev{0.14}  	&\trev{0.25}	&\trev{0.32}	&\trev{0.21}	&\trev{0.23}	&\trev{0.30}	&\trev{0.18}	&\trev{0.20}	&\trev{0.32}	&\trev{0.14}\\
&\trev{Right Lung} &\trev{0.057} &\trev{0.22}	&\trev{0.51}&	\trev{0.14}	&	\trev{0.27}&	\trev{0.48}	&\trev{0.20}	&	\trev{0.25}	&\trev{0.52}	&\trev{0.17}	&	\trev{0.24}&	\trev{0.52}	&\trev{0.14}\\
%&Right Lung PE  &0.78 &	0.44 &	0.55 &	0.24 &	0.41 & 0.48 &	0.28 &	0.42 & 	0.50 &	0.26 &	0.31 &	0.46 &	0.16 \\
\cline{2-15}
& \multicolumn{2}{|c|}{\trev{Dynamic Range}} & \multicolumn{3}{|c|}{\tRev{94\%}}& \multicolumn{3}{|c|}{\tRev{64\%}}& \multicolumn{3}{|c|}{\tRev{79\%}} &\multicolumn{3}{|c|}{\tRev{89\%}} \\
\hline
\multirow{4}{*}{\rotatebox{90}{Susc.} }&Heart  & \trev{0.05} &	 \trev{0.041}	& \trev{0.062}	& \trev{-0.016}	& \trev{0.043}	& \trev{0.061}	& \trev{0.004}	 &\trev{0.048}	& \trev{0.066}	& \trev{0.020}		& \trev{0.040}	& \trev{0.062}	 &\trev{-0.009}\\
&Left lung &  \trev{0.011} & 	\trev{0.027}	&\trev{0.039}	&\trev{-0.007}	&\trev{0.019}	&\trev{0.041}	&\trev{-0.009}	&\trev{0.019}	&\trev{0.032}	&\trev{-0.009}	&\trev{0.015}	&\trev{0.030}	&\trev{-0.015}\\
&\trev{Right Lung} &  \trev{0.011}  	&\trev{0.025}	&\trev{0.050} 	&\trev{-0.017}	&\trev{0.022}	&\trev{0.047}	&\trev{-0.015}	&\trev{0.021}	&\trev{0.043}	&\trev{-0.005}	&\trev{0.021}	&\trev{0.044}	&\trev{-0.009}\\
%&Right Lung PE  & 0.09 &	0.081 &	0.126 &	0.025 &	0.053 &	0.112 &	0.017 &	0.057 &	0.096 &	0.018 &	0.042 &	0.084 &	0.013 \\
\cline{2-15}
& \multicolumn{2}{|c|}{\trev{Dynamic Range}} & \multicolumn{3}{|c|}{\tRev{194\%}}& \multicolumn{3}{|c|}{\tRev{160\%}}& \multicolumn{3}{|c|}{\tRev{171\%}} &\multicolumn{3}{|c|}{\tRev{175\%}} \\
\hline
\hline
\end{tabular}
\end{table}

\trev{A slight rotational effect is observed in the reconstruction with the noisy angles.  The reconstructions of heart and lung regions are blurred together in the susceptivity reconstructions.  This is in part due to the fact that all targets have susceptivity above the background (0~S/m).  Nevertheless, the heart is reconstructed as the most susceptive target in all images as seen in the reconstructed values in Table~\ref{table:ACT4_FranNata_COND}.  Compounding the visual challenges, the susceptivity of the \trev{heart} is about \trev{five} times larger than that of the lungs.  This is a feature of the experimental values rather than a flaw in the algorithm, e.g., see [\cite{Hamilton2013}] for an example of a case with background susceptivity values between those of the targets.  Nevertheless, the main features in the ACT4 images are distinguishable in both the absolute and difference images.}%.  In the reconstructions on the incorrect boundary shape, the heart and effusion blur together slightly but behave similarly in both absolute and difference images.

%\tRev{!!!!***UPDATE the dynamic ranges in the text****!!!} 
\trev{The \trev{conductivity} results in the tables compare most favorably to the true values by considering the maximum values in the \trev{heart region} to the true values and the minimum values in the \trev{lung regions} to the true values.  Using these entries from Table~\ref{table:ACT4_FranNata_COND}, the relative errors for the reconstructions of conductivity on the correct boundary shape are \trev{7.4\%} for the heart \trev{and 153.3\% for the lungs}.  The conductivity of the lungs was consistently overestimated across boundary shapes and electrode positions, whereas the conductivity of the heart region was relative stable.  The relative errors (heart, lungs) for the reconstructions are worst on the oval boundary (12.0\%, 246.0\%) with relative errors of (1.5\%, 206.0\%) for the alternative boundary and (2.7\%, 149.6\%) for the correct boundary with incorrect electrode locations.  The relative errors in susceptivity using the \trev{average regional values} from Table~\ref{table:ACT4_FranNata_COND} were 17.4\% for the heart and 136\% for the lung for the correct boundary shape and electrode angles and (13.0\%, 88.3\%) for the oval boundary, (4.8\%, 82.3\%) for the alternative boundary, and (20.2\%, 64.3\%) for the correct boundary with incorrect electrode locations.  Although the susceptivity cannot physically be negative, the algorithm returned negative values in \trev{some regions}.  The dynamic ranges for the conductivity values on the true boundary, oval boundary, alternative boundary, and correct boundary with incorrect electrode locations were \tRev{94\%, 64\%,  79\%, and 89\%}, respectively.  The dynamic ranges for the susceptivity values on the true boundary, oval boundary, alternative boundary, and correct boundary with incorrect electrode locations were \tRev{194\%, 160\%,  171\%, and 175\%}, respectively\trev{, due to the reconstruction of negative susceptivity values.}}

\trev{We remark that the ACT4 dataset contains high contrast targets of much greater contrast than those of the ACT3 dataset discussed above.  For the ACT4 example, the most conductive object (heart) was nearly 12 times as conductive as the least conductive object (lungs).  By comparison, the heart was approximately 3 times as conductive as the lungs for the ACT3 case.  Recovering such high contrast targets is inherently challenging and will require further examination.}

%\clearpage
%%%%%%%%%%%%%%%%%
\subsection{Reconstructions from ACE1 data}\label{sec:results_ACE1}
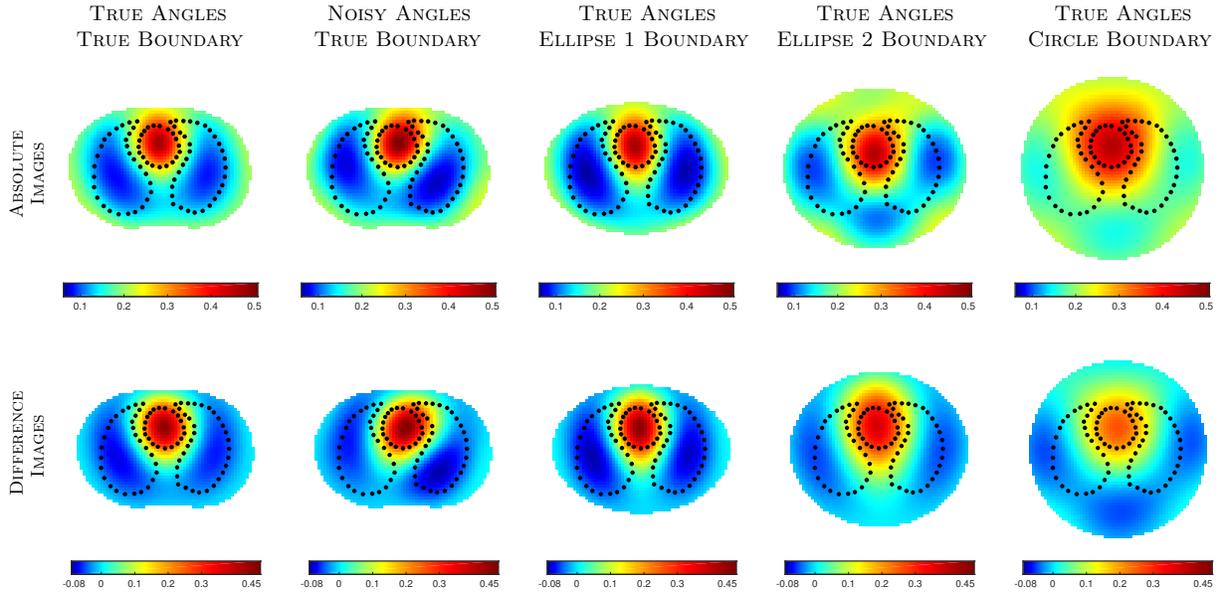
\begin{figure}[!h]

\begin{picture}(420,230)

% Difference images
\put(0,0){\includegraphics[width=75pt]{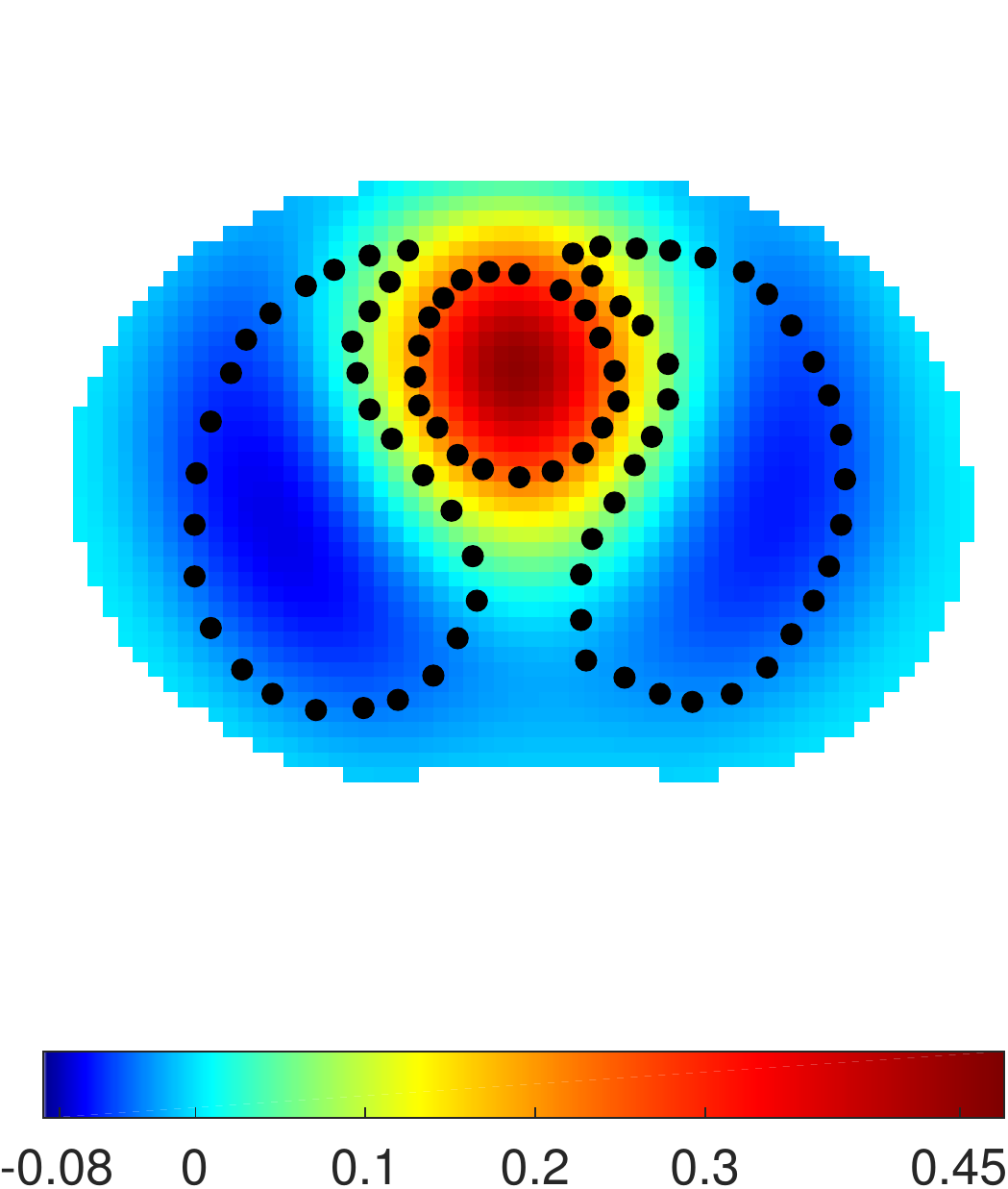}}
\put(90,0){\includegraphics[width=75pt]{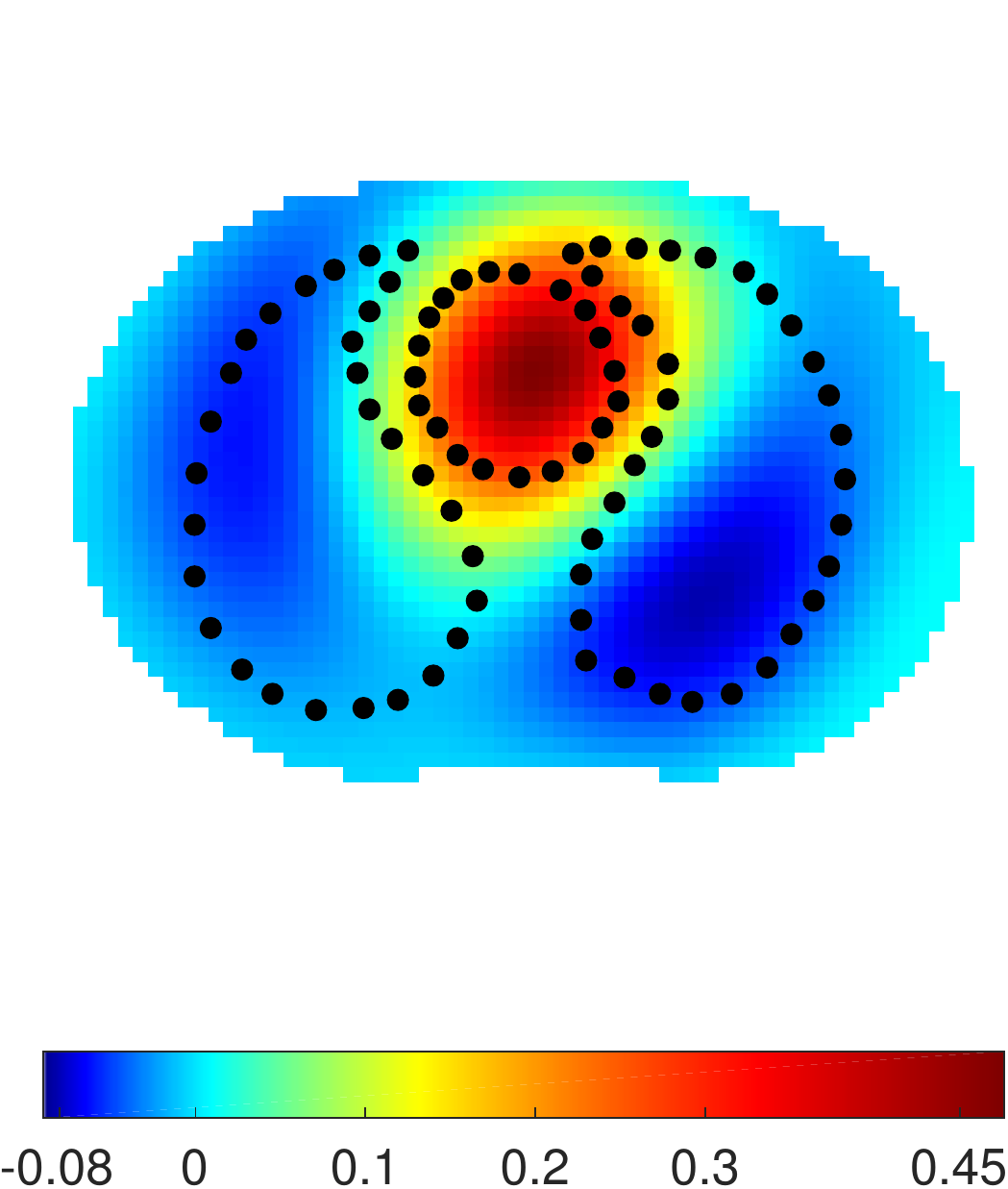}}
\put(180,0){\includegraphics[width=75pt]{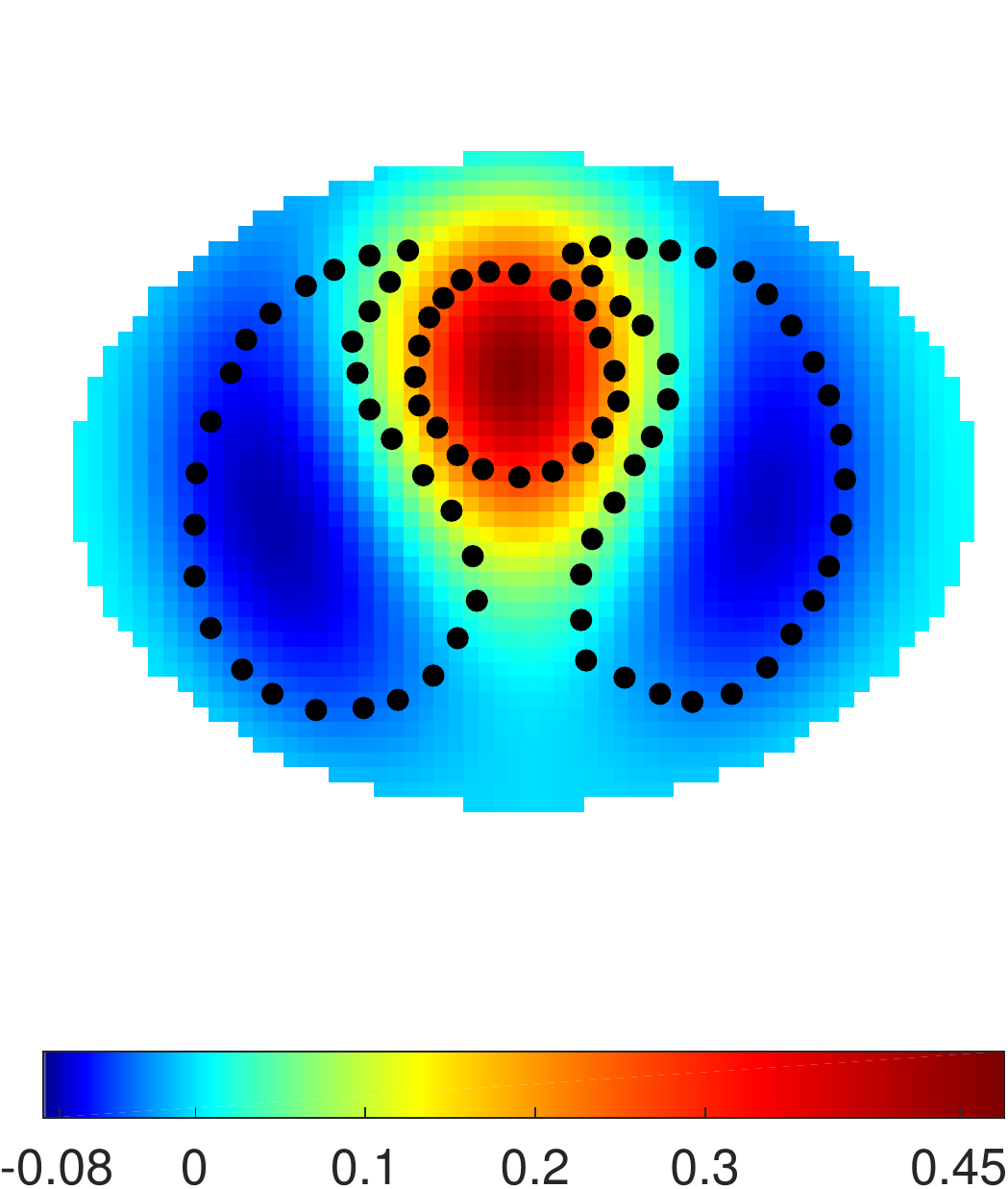}}
\put(270,0){\includegraphics[width=75pt]{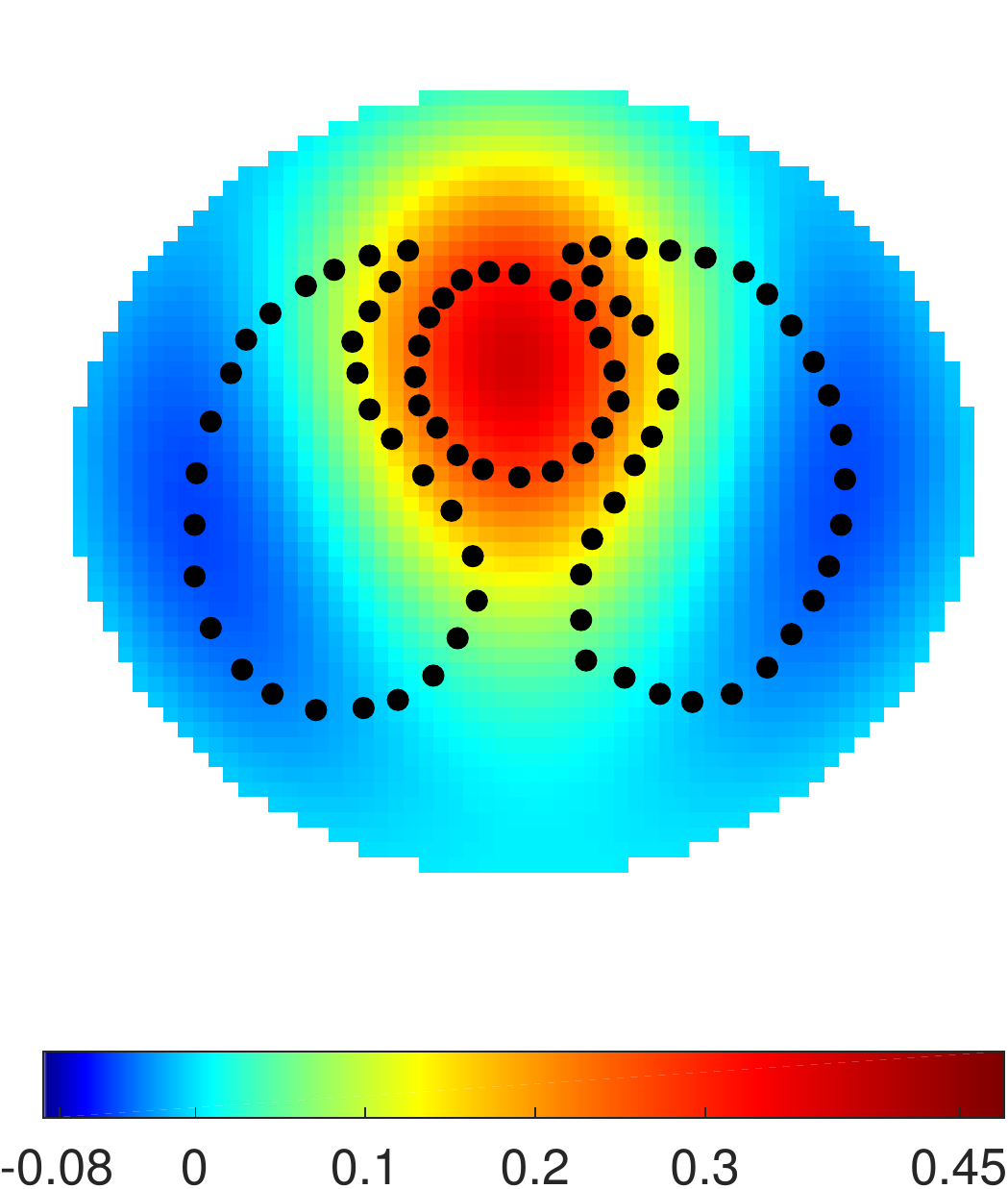}}
\put(360,0){\includegraphics[width=75pt]{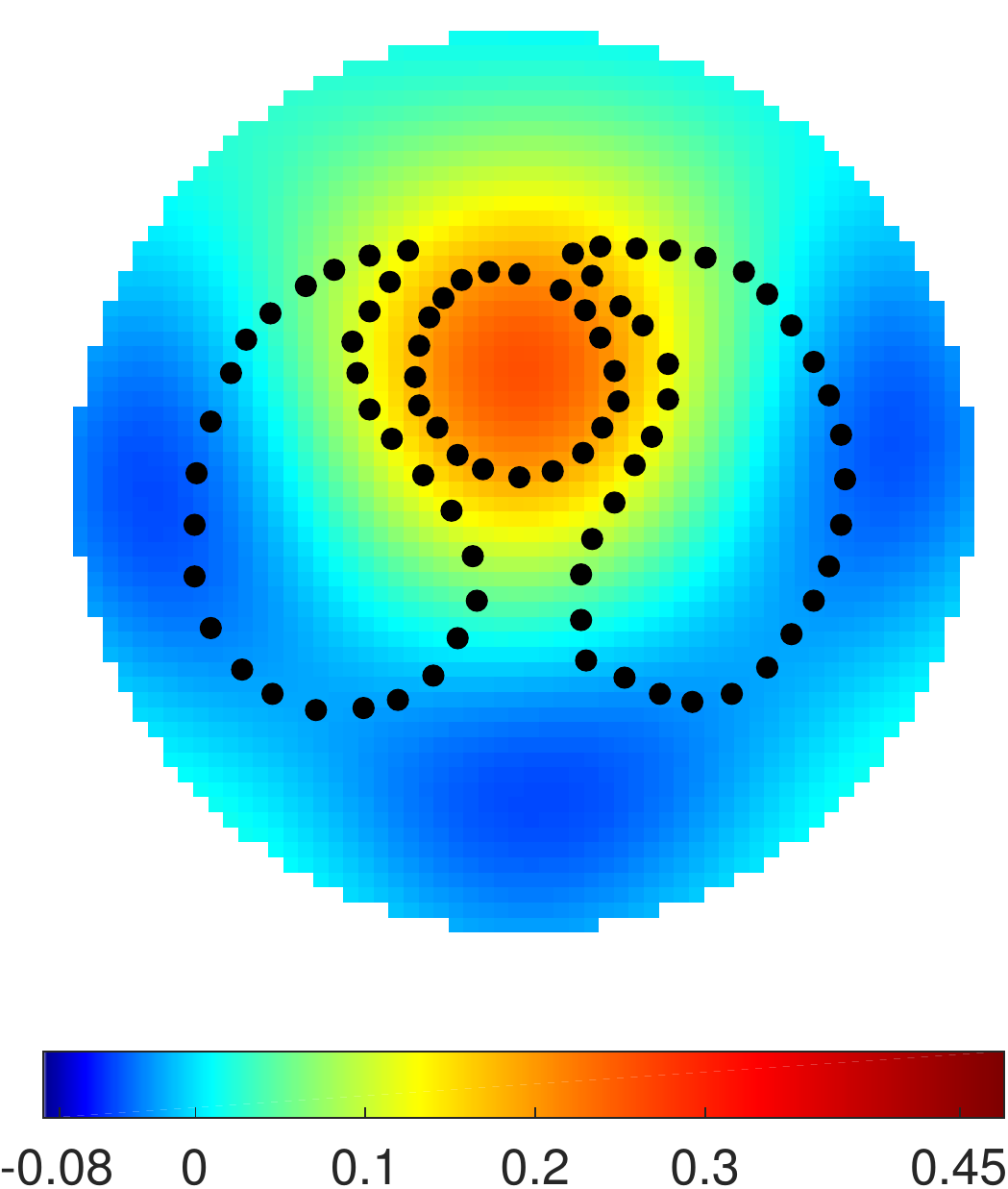}}

% Absolute images
\put(0,105){\includegraphics[width=75pt]{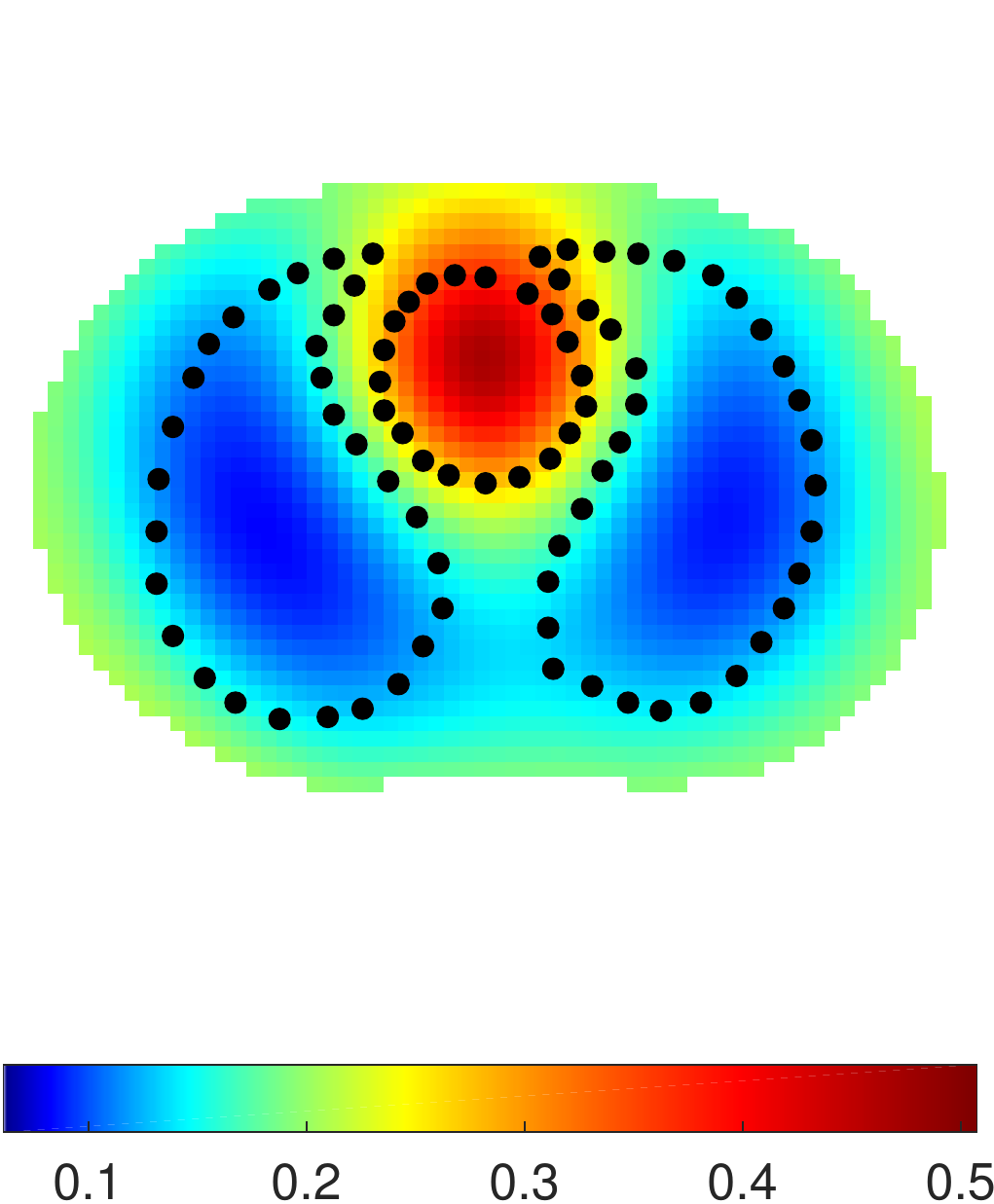}}
\put(90,105){\includegraphics[width=75pt]{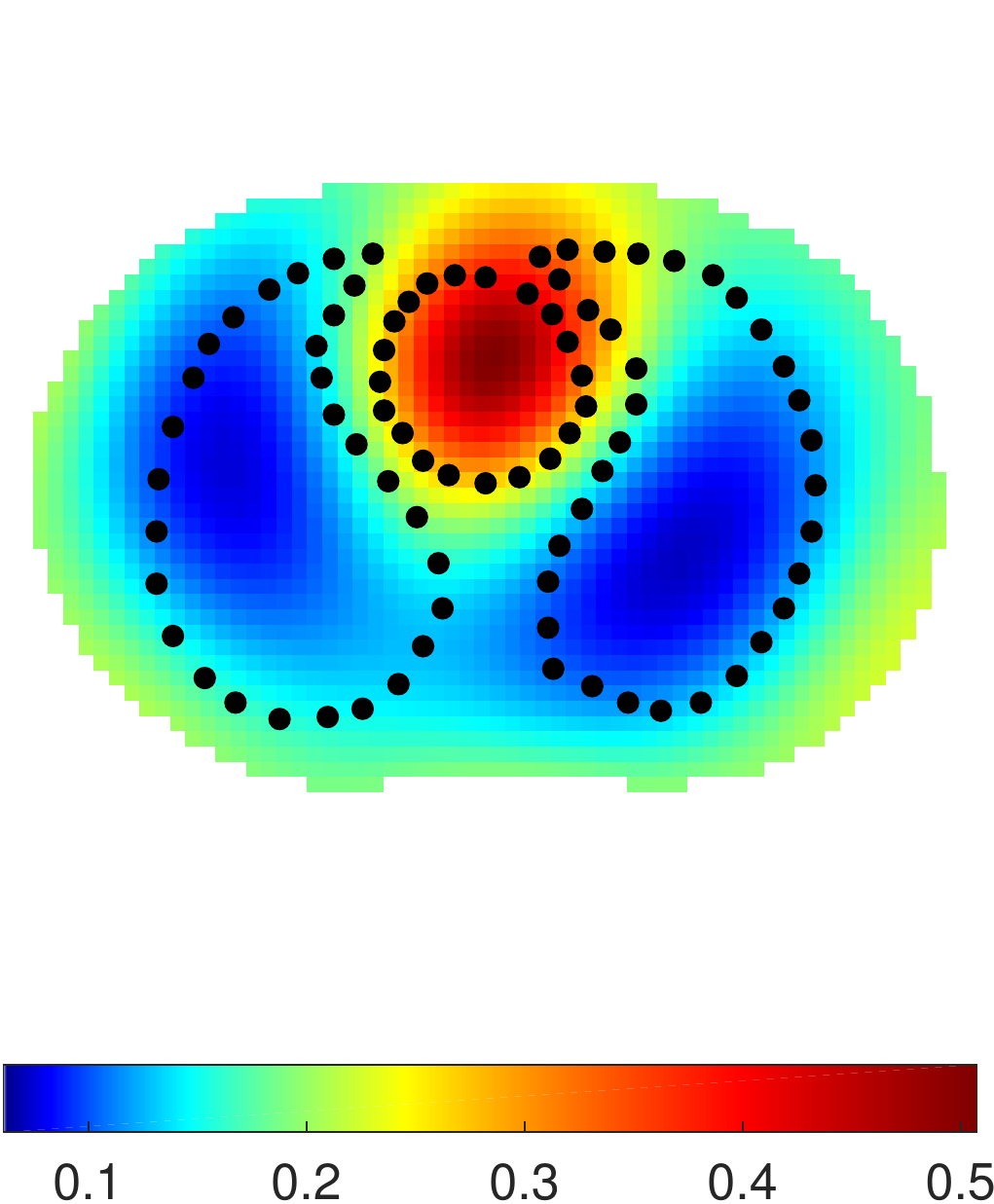}}
\put(180,105){\includegraphics[width=75pt]{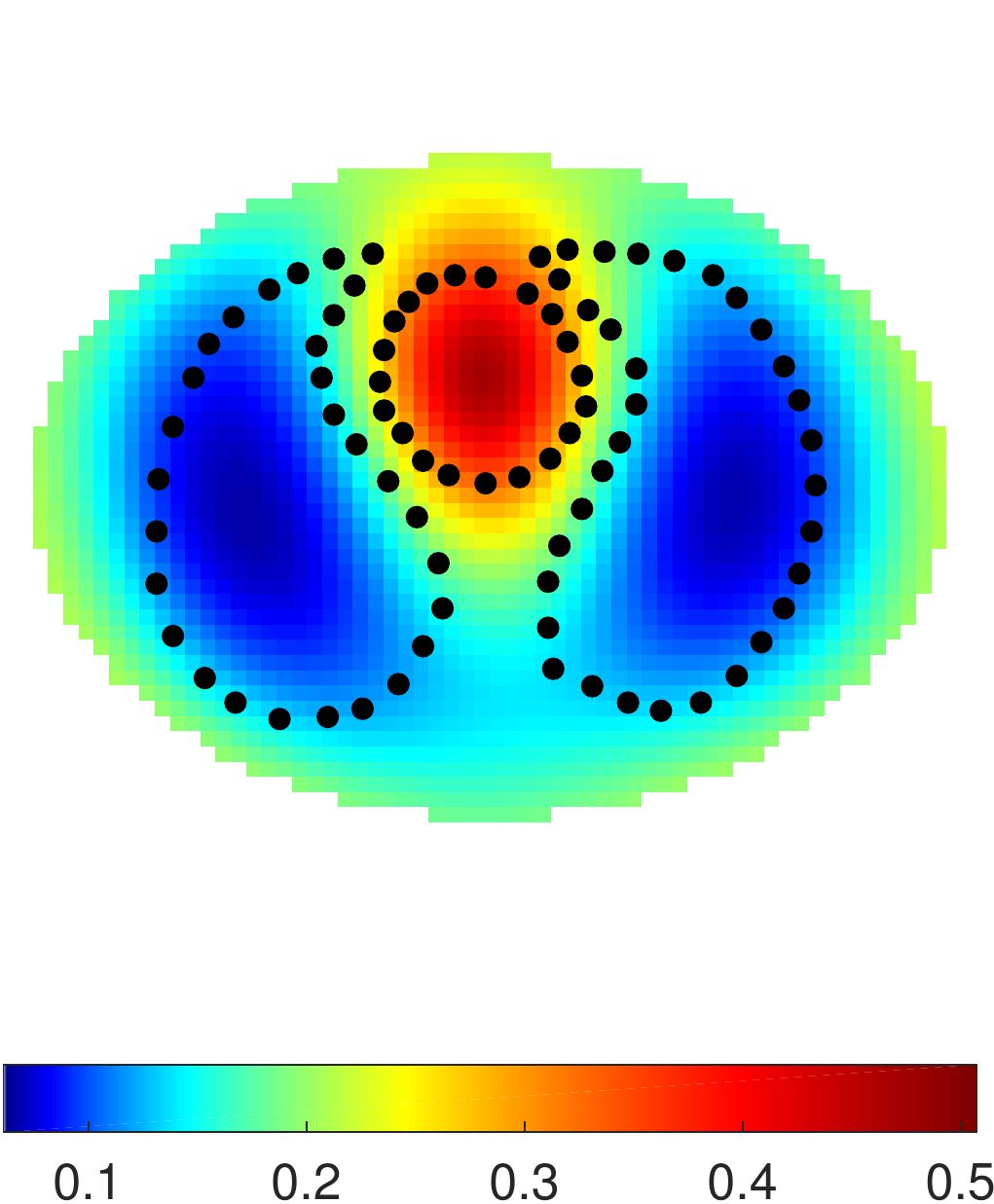}}
\put(270,105){\includegraphics[width=75pt]{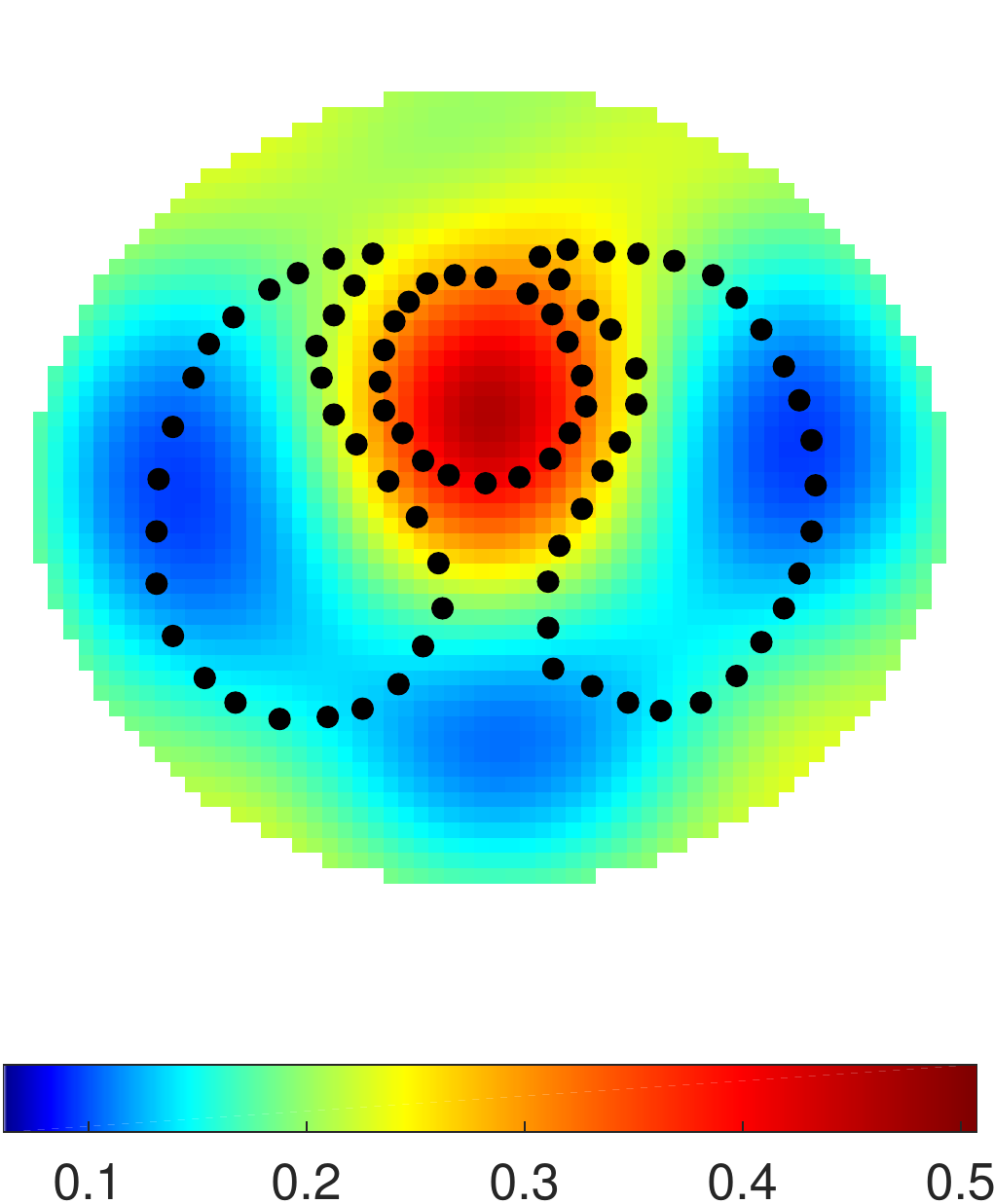}}
\put(360,105){\includegraphics[width=75pt]{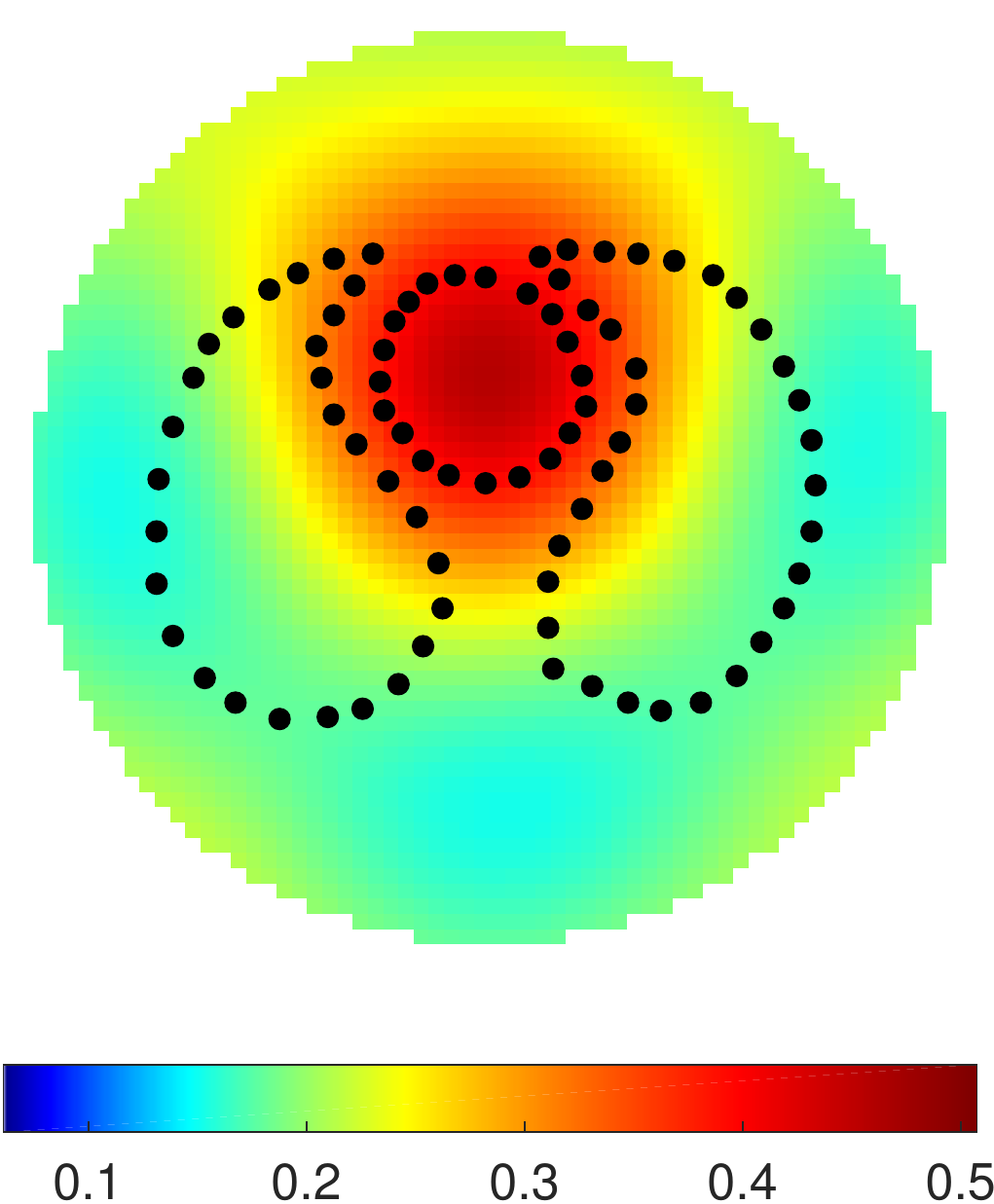}}

% side labels 
\put(-20,35){\rotatebox{90}{\scriptsize{\sc Difference}}}
\put(-12,42){\rotatebox{90}{\scriptsize{\sc Images}}}

\put(-20,140){\rotatebox{90}{\scriptsize{\sc Absolute}}}
\put(-12,145){\rotatebox{90}{\scriptsize{\sc Images}}}

% top labels:
\put(10,215){\footnotesize {\sc True Angles}}
\put(5,205){\footnotesize {\sc True Boundary}}

\put(100,215){\footnotesize {\sc Noisy Angles}}
\put(95,205){\footnotesize {\sc True Boundary}}

\put(195,215){\footnotesize {\sc True Angles}}
\put(180,205){\footnotesize {\sc Ellipse 1 Boundary}}

\put(285,215){\footnotesize {\sc True Angles}}
\put(270,205){\footnotesize {\sc Ellipse 2 Boundary}}

\put(375,215){\footnotesize {\sc True Angles}}
\put(365,205){\footnotesize {\sc Circle Boundary}}

\end{picture}
\caption{\label{fig:ACE1_HL_FranNata_pula0} Comparison of conductivity \trev{(S/m)} reconstructions for the agar heart and lungs ACE1 phantom (see, Figure~\ref{fig:phantoms}, third).  Results are compared for knowledge of true vs. incorrect electrode angles as well as boundary shape.  Absolute images requiring no $\Lambda_1$ are in row 1\tRev{, plotted on the same color scale}.  Difference images are in row 2\tRev{, plotted on the same color scale.}
 The reconstructions for the `Circular' boundary were performed on a $k$-disc of radius $3$ rather than 4 for stability.}
\end{figure}

%%%%%%%%%%%%%%%%%
\trev{Figure}~\ref{fig:ACE1_HL_FranNata_pula0} displays the conductivity reconstructions from the \trev{ACE1 data. Reconstructions} were computed using  {\it `Approach~2'} with the true and incorrect electrode locations and boundary shapes shown in Figure~\ref{fig:CSU_elecSetups}.   \trev{The reference data for the difference images contained only saline with conductivity $0.2$ S/m.}  Regional max, min, and average \trev{values are reported} in Table~\ref{table:ACE1_FranNata_Agar_skip0}.% and in Table~\ref{table:ACE1_FranNata_MelonCuc} for the melon and cucumber targets.

%\tRev{Reconstructions for the remainder of the manuscript focus solely on the {\it 'Approach~2'} D-bar method due to it's ability to form absolute images without a need to simulate $\Lambda_1$ measurements, as well as flexibility for conductivity or admittivity imaging.} \trev{Figure}~\ref{fig:ACE1_HL_FranNata_pula0} displays the conductivity reconstructions from the \trev{ACE1 data. Reconstructions} were computed using  {\it `Approach~2'} with the true and incorrect electrode locations and boundary shapes shown in Figure~\ref{fig:CSU_elecSetups}.   \trev{The reference data for the difference images contained only saline with conductivity $0.2$ S/m.}  Regional max, min, and average \trev{values are reported} in Table~\ref{table:ACE1_FranNata_Agar_skip0}.% and in Table~\ref{table:ACE1_FranNata_MelonCuc} for the melon and cucumber targets.

%%%%%%%%%%%%%%%%%%%%%%%%%%%%%%%%%%%%%%%%%%%%%%%%%
% ACE1 HL -  Conductivity only - NATA style - pula 0- NONLINEAR SCALING OF COLORBARS
%%%%%%%%%%%%%%%%%%%%%%%%%%%%%%%%%%%%%%%%%%%%%%%%%
%%%%%%%%%%%%%%%%%

%%%%%%%%%%%%%%%%%%%%%%%%%%%%%%%%%%%%%%%%%%%%%%%%%
\begin{table}[h!]
\scriptsize\centering
\caption{Max, min, and average conductivity values \trev{(S/m)} in each of the organ regions in the reconstructions from the ACE1 data (see Figure~\ref{fig:phantoms}, third) using {\it `Approach~2'} of Section \ref{sec:FranDbar}.}
\label{table:ACE1_FranNata_Agar_skip0}
\begin{tabular}{|ll|ccc|ccc|ccc|ccc|ccc|}
\hline
\multicolumn{2}{|c}{\multirow{2}{*}{ }}  &  \multicolumn{3}{|c|}{Correct Boundary} &  \multicolumn{3}{c|}{Noisy Angles} &  \multicolumn{3}{c|}{\trev{Ellipse 1}} &  \multicolumn{3}{c|}{\trev{Ellipse 2}} &   \multicolumn{3}{c|}{Circle} \\
%&  &  \multicolumn{3}{|c|}{and Elec Angles} &  \multicolumn{3}{c|}{Noisy  Angles} &  \multicolumn{3}{c|}{True Elec Angles} &  \multicolumn{3}{c|}{True Elec Angles} &  \multicolumn{3}{c|}{True Elec Angles} \\
\hline
%\cline{2}
Organ & True & AVG &	MAX	& MIN	&	AVG &	MAX	& MIN	&	AVG &	MAX	& MIN	&	AVG &	MAX	& MIN	&	AVG &	MAX	& MIN \\
\hline
Heart  & 0.45 &	0.37&	0.47&	0.27	&	0.39	&0.51&	0.27		&0.38	&0.47	&0.27	&	0.39&	0.46	&0.29	&	0.43&	0.46&	0.39\\
Left lung & 0.09 & 	0.13	&0.29	&0.09		&0.12	&0.34	&0.07		&0.11	&0.27	&0.07		&0.16	&0.29	&0.10	&	0.23	&0.38	&0.17\\
Right Lung   & 0.09 &0.12&	0.23	&0.08		&0.12	&0.19	&0.08		&0.11	&0.21&	0.07	&	0.15&	0.29	&0.10		&0.22&	0.33	&0.16\\
\hline
\multicolumn{2}{|c|}{\trev{Dynamic Range}} & \multicolumn{3}{|c|}{\trev{107.5\%}}& \multicolumn{3}{|c|}{\trev{120.7\%}}& \multicolumn{3}{|c|}{\trev{112.5\%}} &\multicolumn{3}{|c|}{\trev{102\%}} &\multicolumn{3}{|c|}{\trev{84\%} }\\
\hline
\hline
\end{tabular}
\end{table}

Reconstructions from the ACE1 \trev{data further} demonstrate the method {\it `Approach 2'} of Section~\ref{sec:FranDbar} is very robust to errors in boundary shape and electrode position \trev{(Figure~\ref{fig:ACE1_HL_FranNata_pula0}, Table~\ref{table:ACE1_FranNata_Agar_skip0})}.  The noisy angles again have the effect of rotating the reconstructed image, and the elliptical and circular boundaries cause the reconstructed lungs to be pushed toward the boundary, and the heart toward the center, with the effect becoming more pronounced as the boundary becomes more circular.  As the domain becomes more circular, a low-conductivity artefact appears opposite the heart and between the lungs, which can most likely be attributed to the distortion of the lungs toward the boundary.  

As was the case for the ACT4 data, the conductivity values in the tables compare most favorably to the true values by considering the maximum values in the heart to the true values and the minimum values in the lung to the true values.  Using these entries from Table~\ref{table:ACE1_FranNata_Agar_skip0}, the relative errors for the reconstructions of conductivity on the correct boundary shape are 4.4\% for the heart, 0\% for the left lung, and 11\% for the right lung.   The errors for the incorrect boundary shapes for the heart were 13.3\% for the noisy angles, 4.4\% on \trev{Ellipse 1} boundary, and 2.2\% on \trev{Ellipse 2} boundary and the circle.  The errors in the lungs regions were 22\% and 11\% for the noisy angles, 22\% for \trev{Ellipse 1} boundary, 11\% for \trev{Ellipse 2} boundary, and 89\% and 78\% for the circle, due to the distortion toward the boundary.  
The dynamic ranges for the conductivity values of the agar targets using the true boundary, true boundary with incorrect electrode locations, \trev{Ellipse 1} boundary, \trev{Ellipse 2} boundary, and circular boundary were  107.5\%, 120.7\%, 112.5\%, 102\%, and 84\%, respectively. % These are for skip 0
% skip 2 results:  90\%, 142\%,  90\%, 101\%, and 86\%, respectively. 

\begin{figure}[!h]
\begin{picture}(420,230)

% Difference images
\put(0,0){\includegraphics[width=75pt]{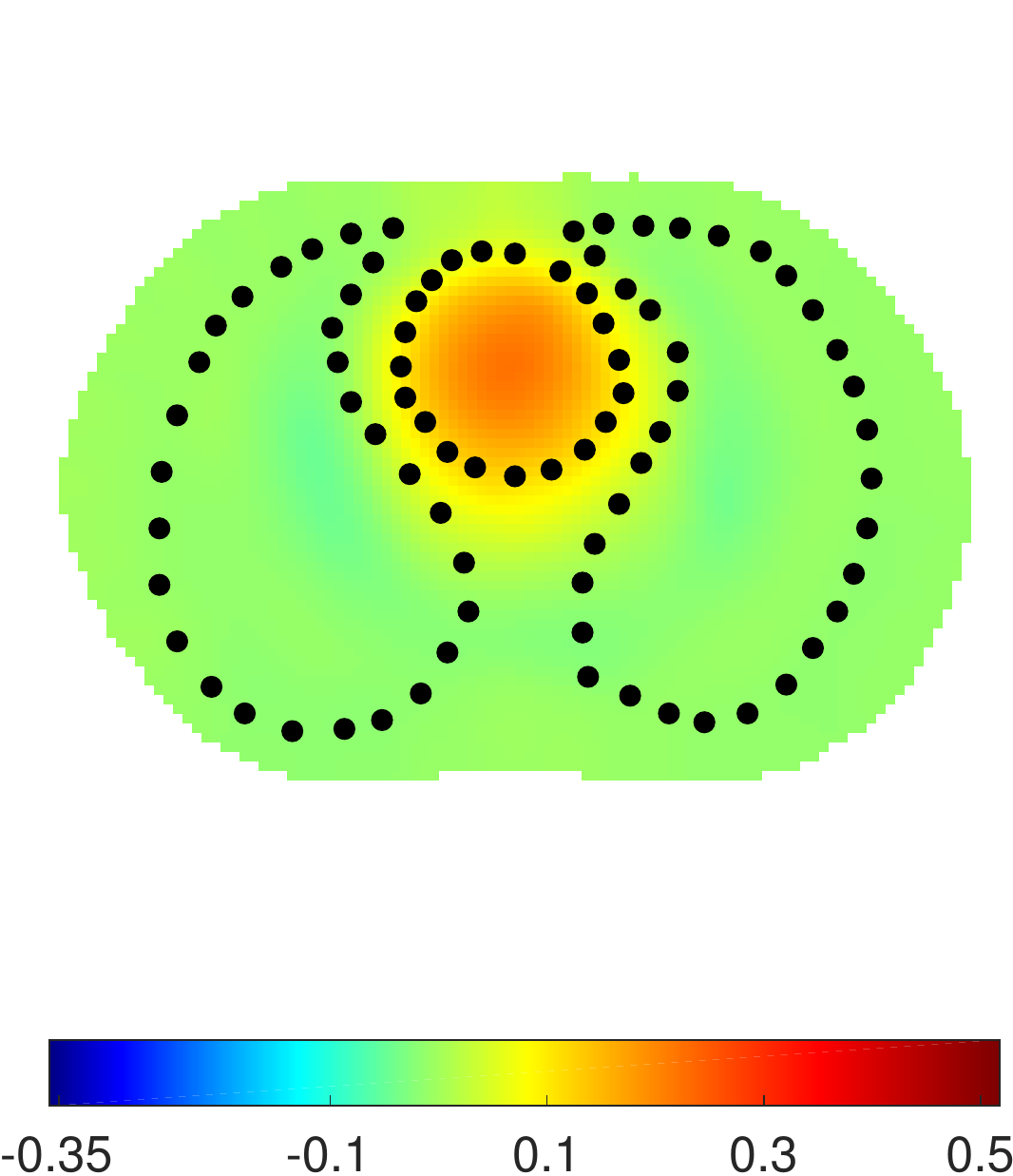}}
\put(90,0){\includegraphics[width=75pt]{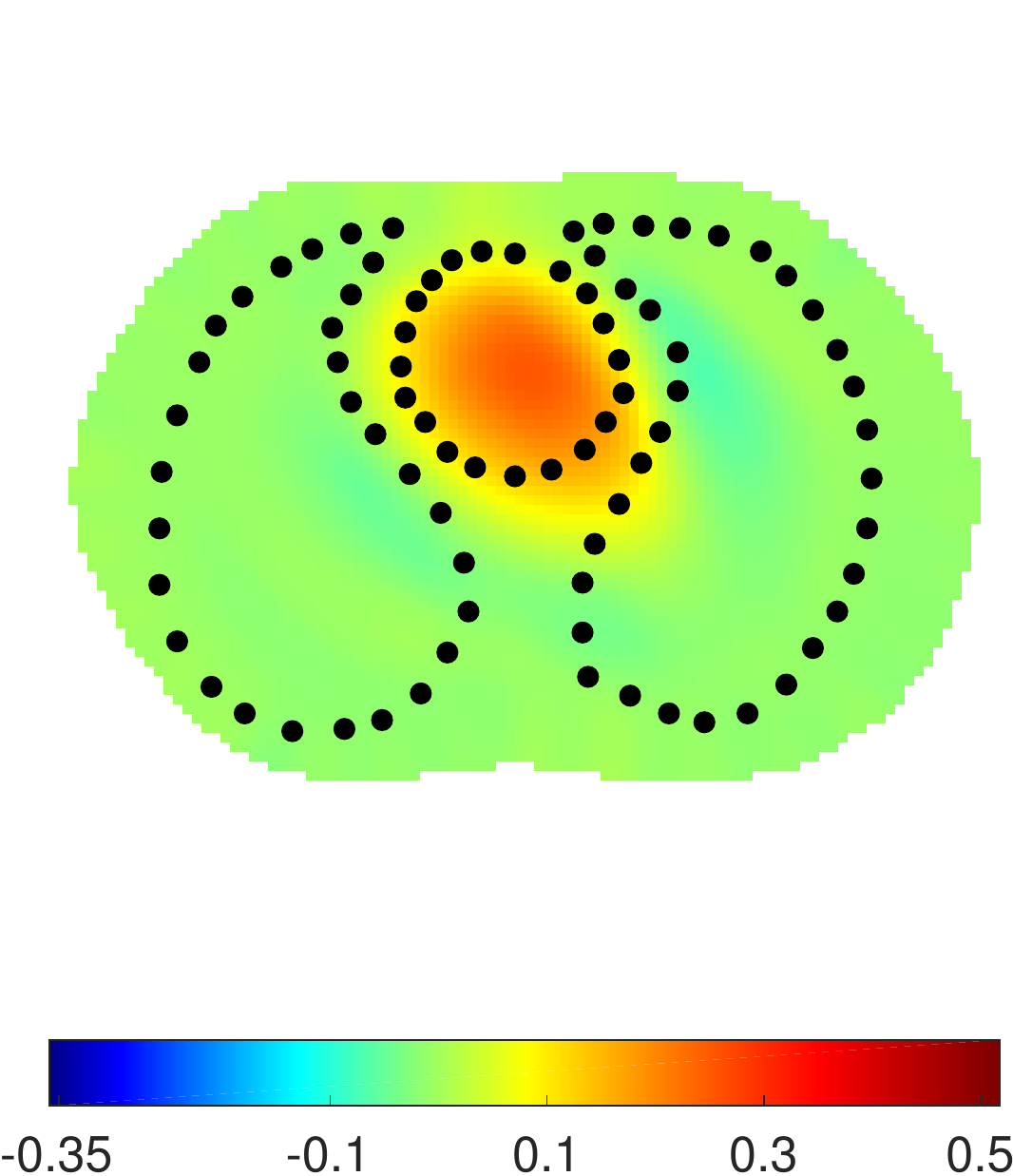}}
\put(180,0){\includegraphics[width=75pt]{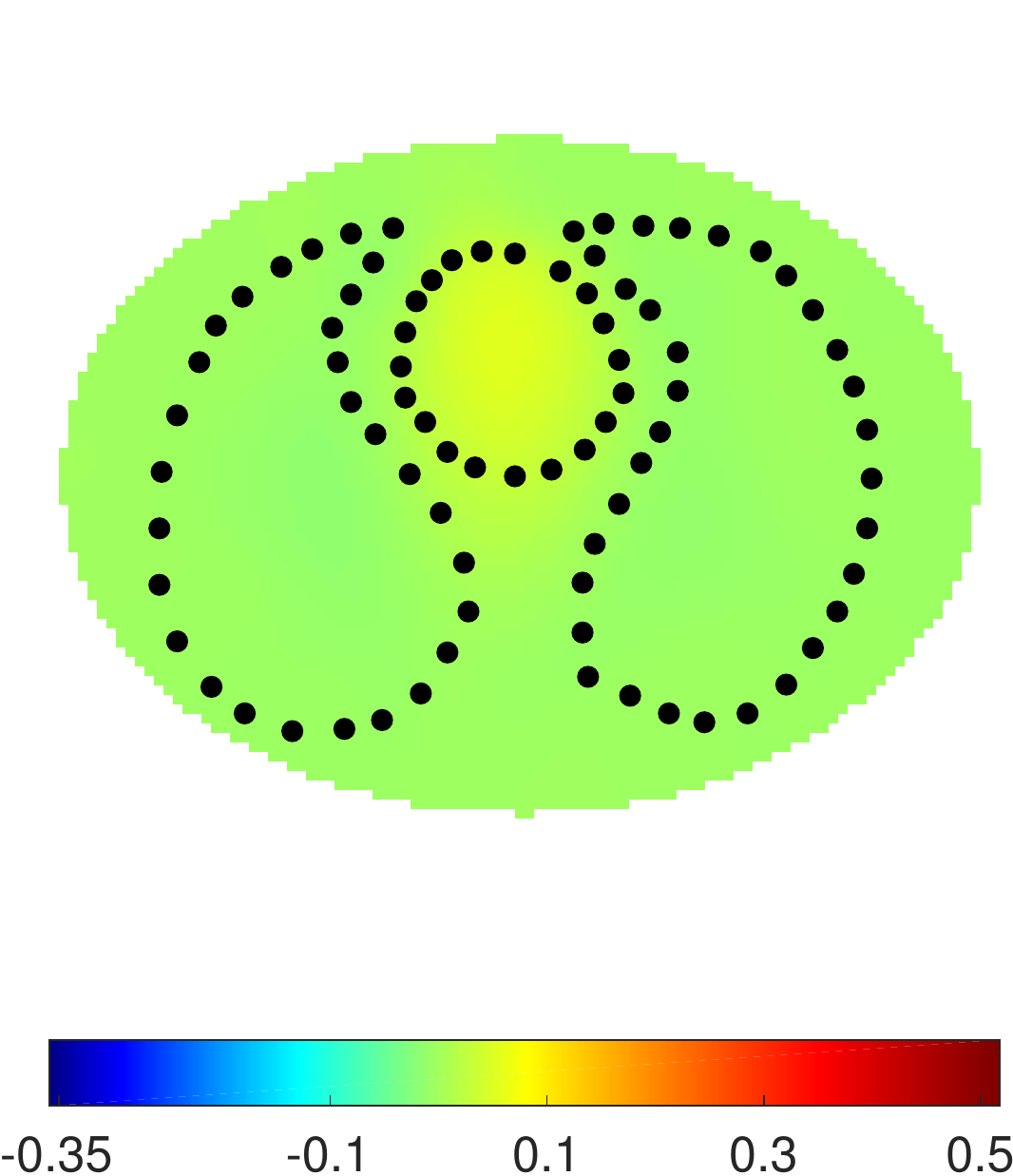}}
\put(270,0){\includegraphics[width=75pt]{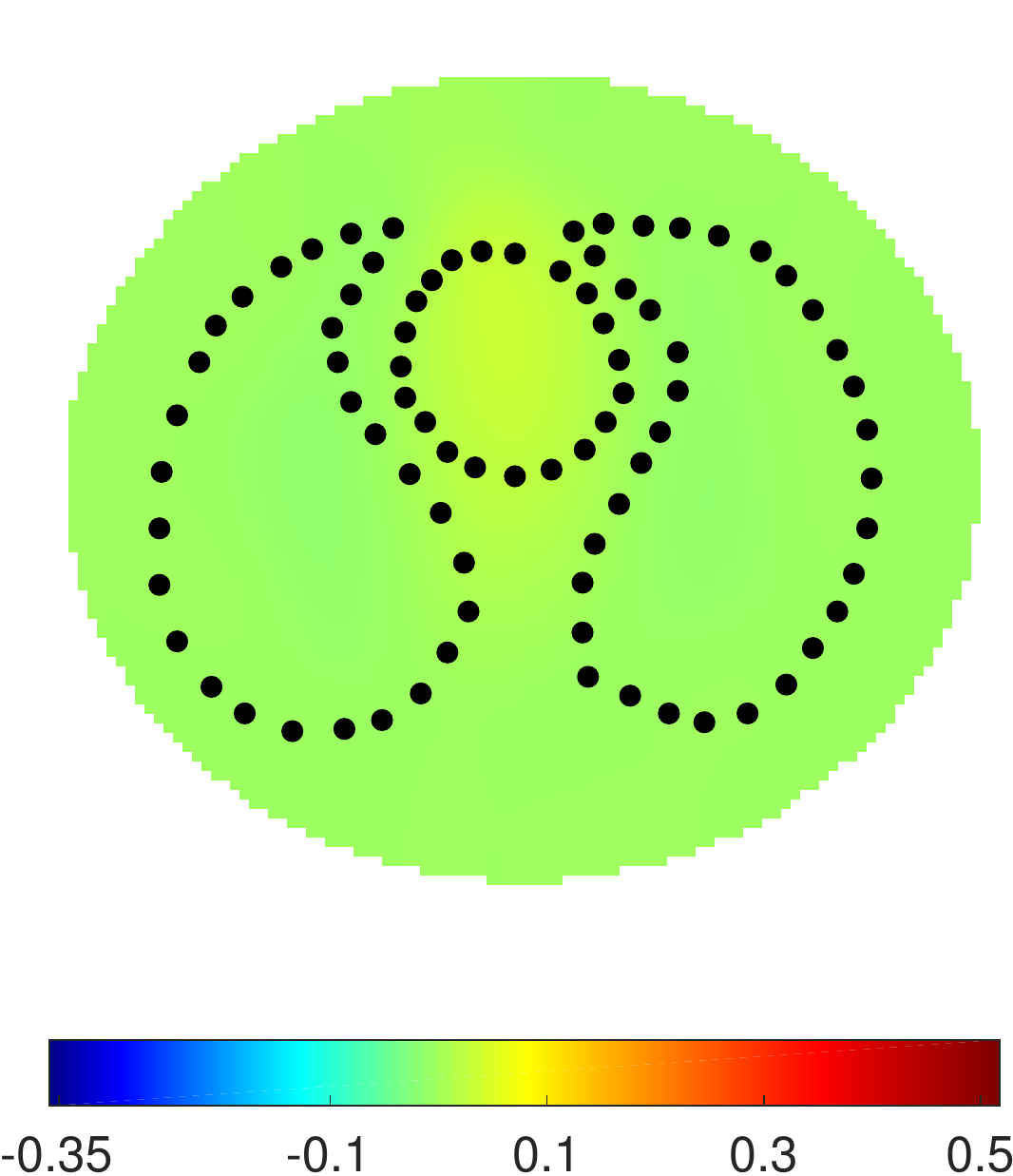}}
\put(360,0){\includegraphics[width=75pt]{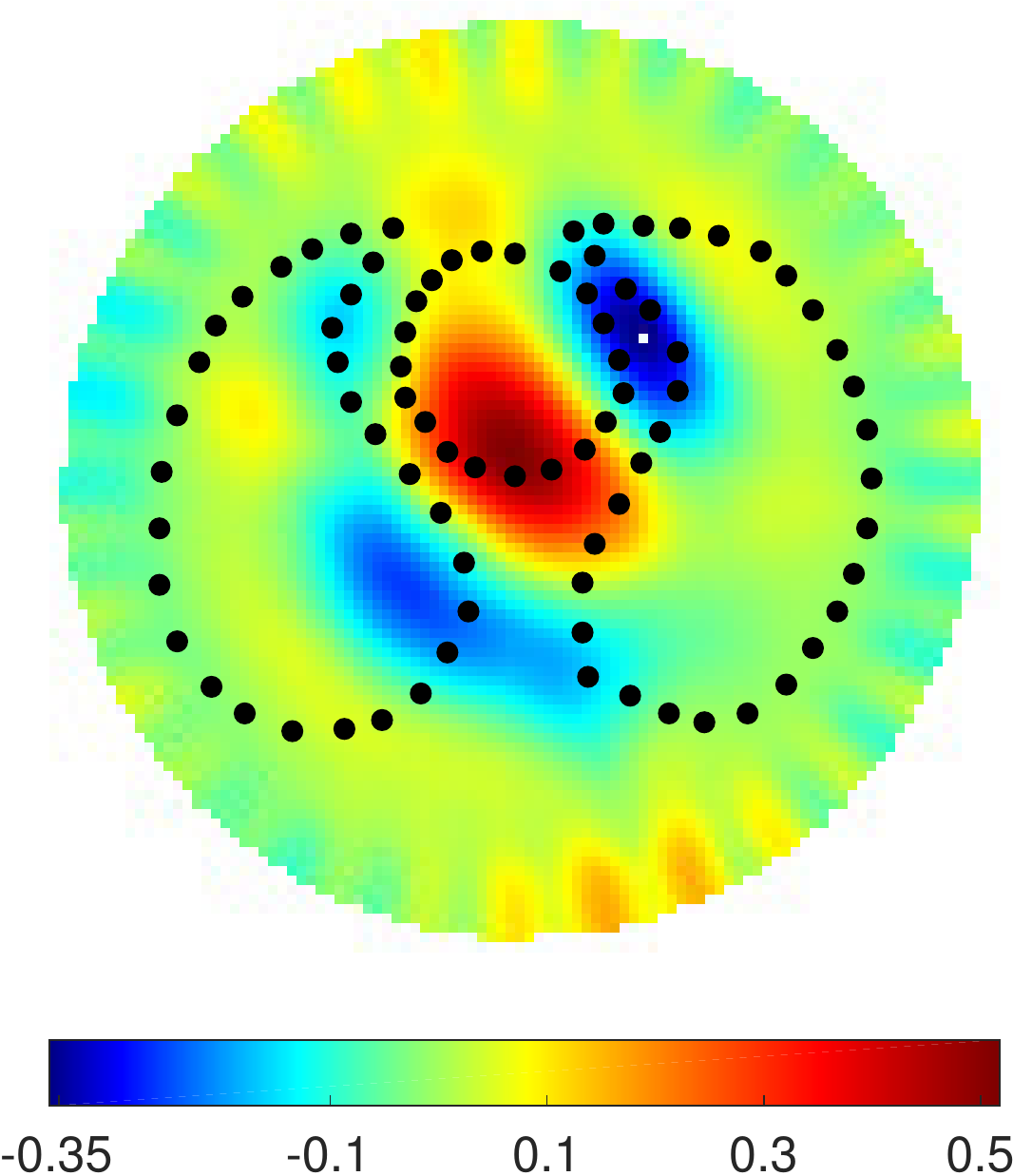}}

% Absolute images
\put(0,105){\includegraphics[width=75pt]{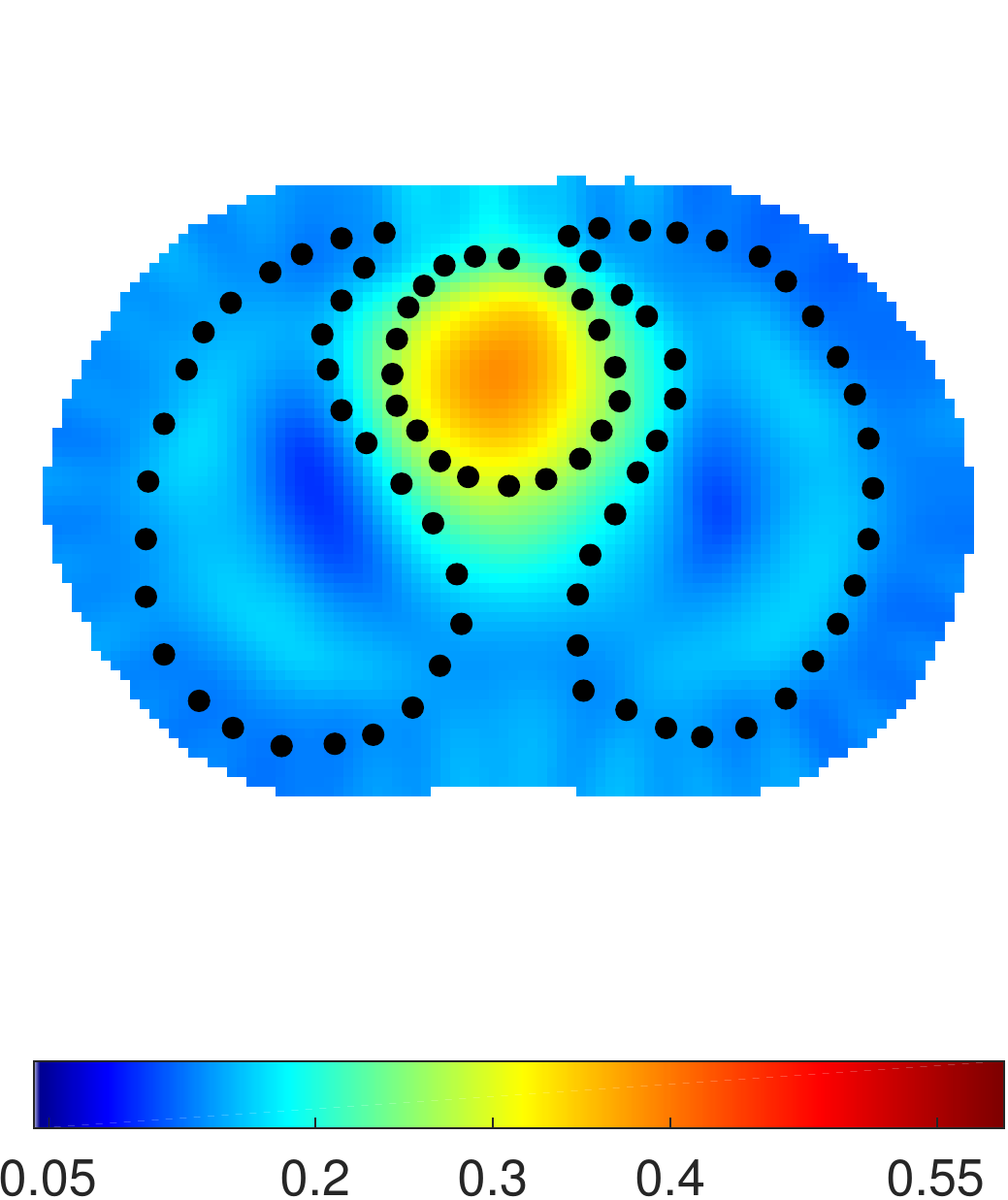}}
\put(90,105){\includegraphics[width=75pt]{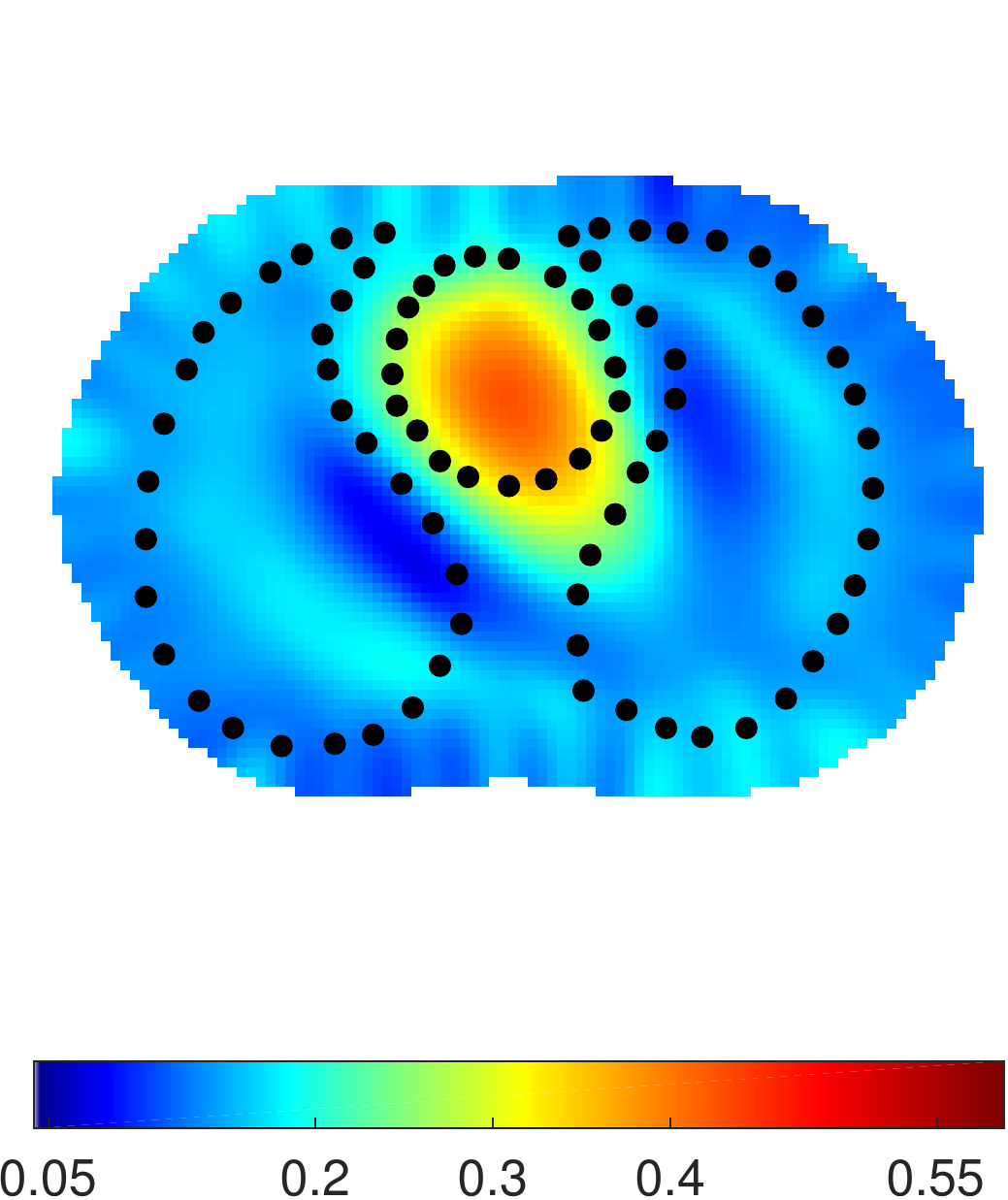}}
\put(180,105){\includegraphics[width=75pt]{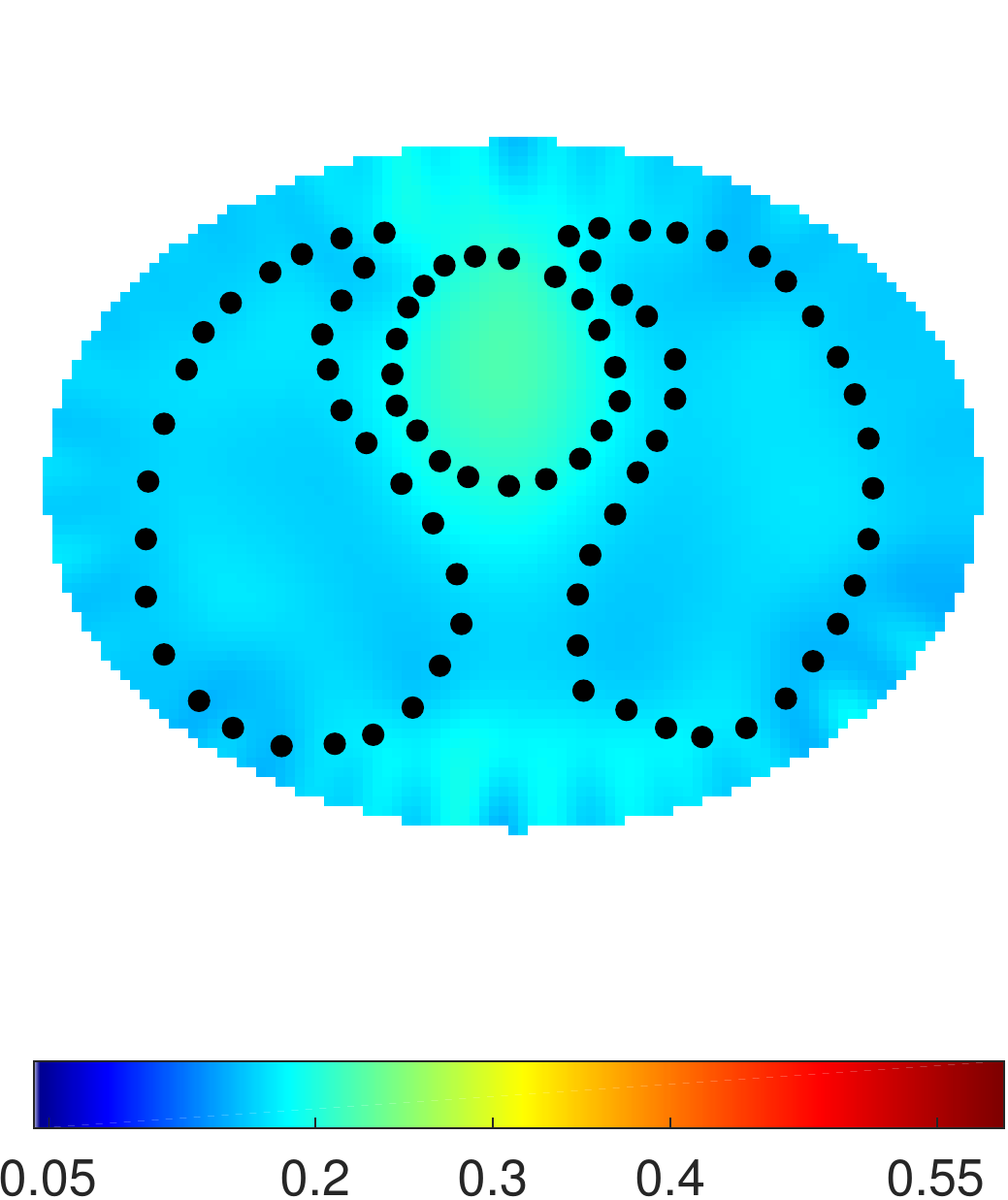}}
\put(270,105){\includegraphics[width=75pt]{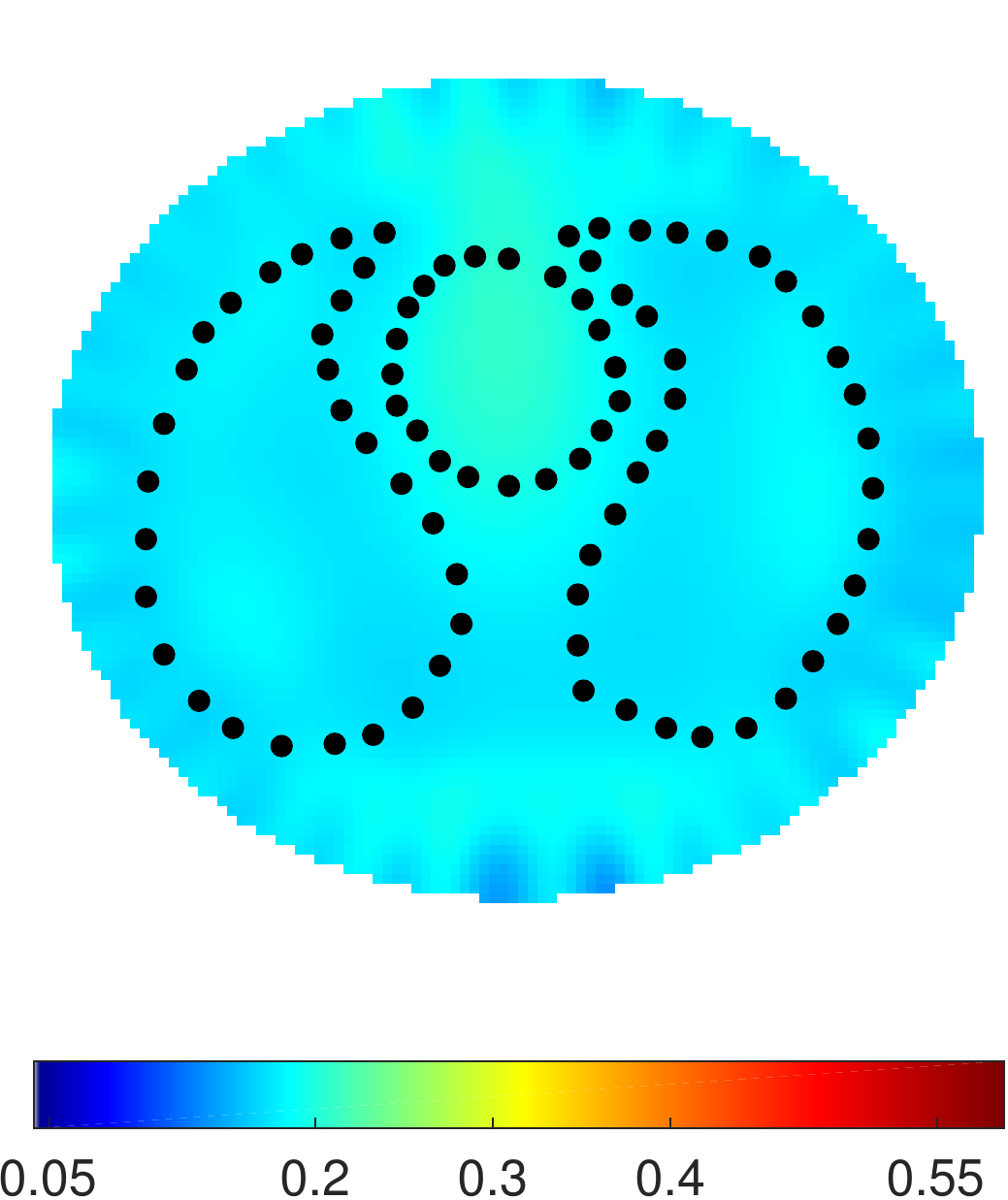}}
\put(360,105){\includegraphics[width=75pt]{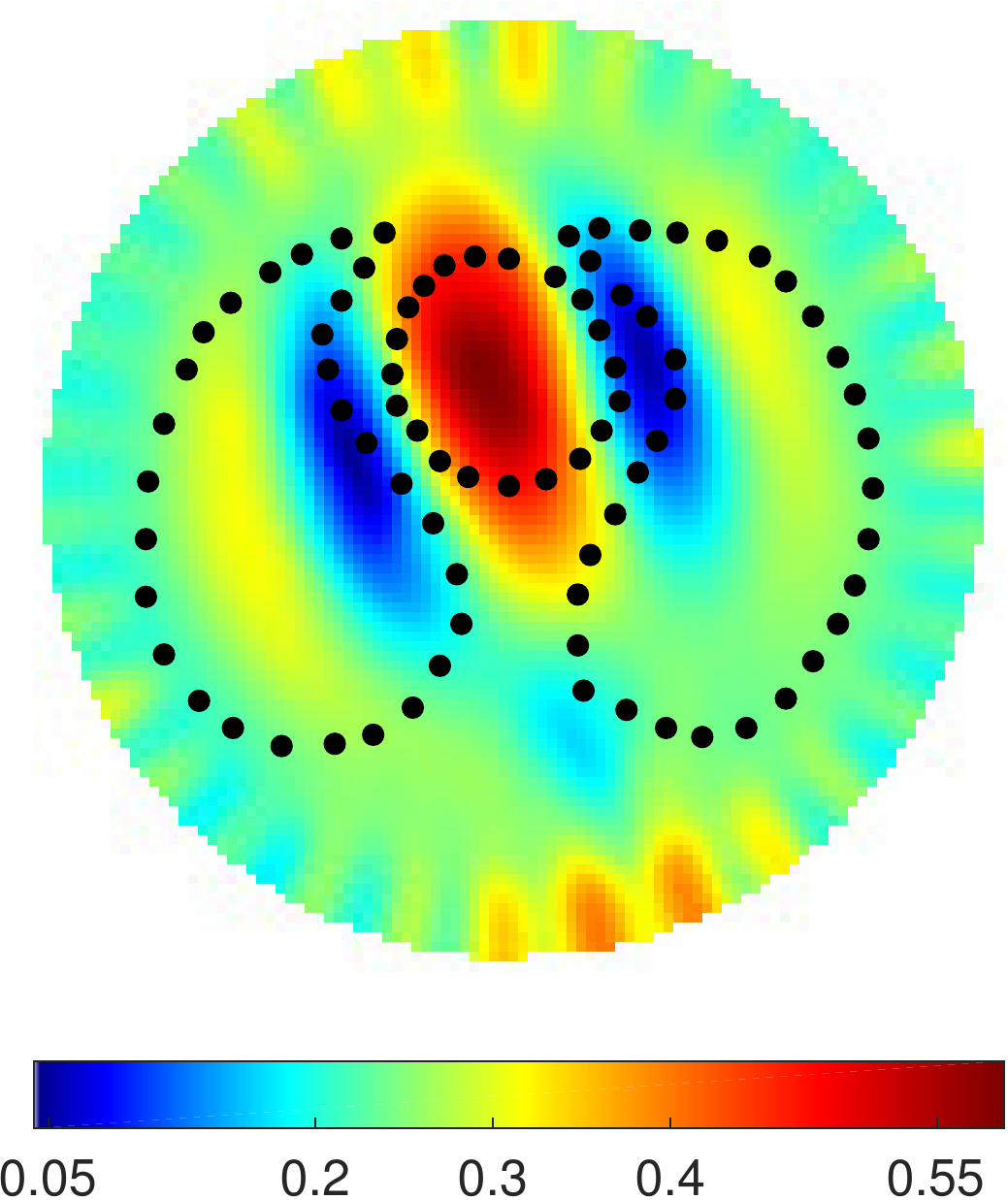}}

% side labels 
\put(-20,35){\rotatebox{90}{\scriptsize{\sc Difference}}}
\put(-12,42){\rotatebox{90}{\scriptsize{\sc Images}}}

\put(-20,140){\rotatebox{90}{\scriptsize{\sc Absolute}}}
\put(-12,145){\rotatebox{90}{\scriptsize{\sc Images}}}

% top labels:
\put(10,215){\footnotesize {\sc True Angles}}
\put(5,205){\footnotesize {\sc True Boundary}}

\put(95,215){\footnotesize {\sc Shifted Angles}}
\put(95,205){\footnotesize {\sc True Boundary}}

\put(195,215){\footnotesize {\sc True Angles}}
\put(180,205){\footnotesize {\sc Ellipse 1 Boundary}}

\put(285,215){\footnotesize {\sc True Angles}}
\put(270,205){\footnotesize {\sc Ellipse 2 Boundary}}

\put(375,215){\footnotesize {\sc True Angles}}
\put(365,205){\footnotesize {\sc Circle Boundary}}

\end{picture}
\caption{\label{fig:ACE1_HL_Talles}  Comparison of conductivity  \trev{(S/m)} reconstructions for the agar heart and lungs ACE1 phantom (see, Figure~\ref{fig:phantoms}, third) using the Gauss-Newton reconstruction algorithm with adjacent current patterns.  Results are compared for knowledge of true vs. incorrect electrode angles as well as boundary shape.  \tRev{Absolute images are in row 1, plotted on the same color scale.  Difference images are in row 2, plotted on the same color scale.}}
\end{figure}
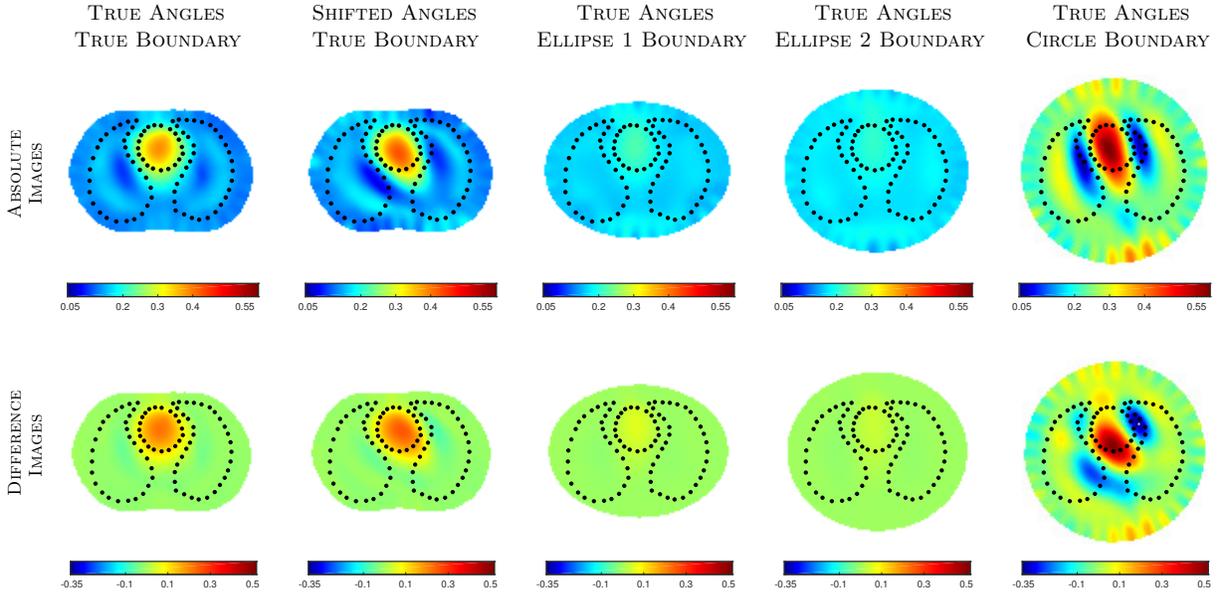

Figure~\ref{fig:ACE1_HL_Talles} shows the Gauss-Newton method reconstructions under the various incorrect boundary and electrode configurations.  Since the FEM for the solution of the forward problem at each iteration requires that the electrodes do not overlap, a different electrode perturbation was used from that of the D-bar reconstructions shown in Figure~\ref{fig:ACE1_HL_FranNata_pula0} (second).  Instead, all electrode centers were first perturbed by the angle $\pi/32$, and then a random vector drawn from the uniform distribution on $[0, Egap/2]$, where $Egap$ is the (uniform) gap between the actual electrodes, was added to the perturbed locations. The Gauss-Newton reconstructions are shown to demonstrate common artefacts that arise in minimization methods without precise knowledge of the boundary or electrode locations.

%%%%%%%%%%%%%%%%%%%%%%%%%%%%%%%%%%%%%%%%%%%%%%%%%
% ACE1 HL -  Conductivity only - Gauss Newton - pula 0
%%%%%%%%%%%%%%%%%%%%%%%%%%%%%%%%%%%%%%%%%%%%%%%%%
%
\trev{T}he Gauss-Newton algorithm \trev{(Figure~\ref{fig:ACE1_HL_Talles})} produced an absolute image with good spatial resolution of the heart, but very small lung regions, which did improve significantly in the difference image.  The small size of the reconstructed lungs was a feature that persisted through the reconstructions on the incorrect boundaries and the shifted electrodes.   As errors in the boundary shape and electrode positions were introduced, conductive regions along the domain boundary appear in the absolute images due to the presence of the conductive electrodes, and were also present in the difference image on the circle.  The shifts in the electrode positions result in a rotated reconstruction, both in the absolute and difference images.  The dynamic ranges for the reconstructed conductivity values using the true boundary, true boundary with incorrect electrode locations, \trev{Ellipse 1} boundary, \trev{Ellipse 2} boundary, and circular boundary were  28.5\%, 35.1\%, 7.2\%, 7\%, and 53.8\%, respectively.  Although the Gauss-Newton method does not include any linearizing assumptions, the regularization term used here serves to damp the maximum value of the conductivity distribution, resulting in a low dynamic range.  Other regularization terms such as total variation regularization or sparsity-promoting strategies may result in sharper images and a higher dynamic range.  However, a thorough comparison with other algorithms is not \trev{in} the scope of this paper; rather we include these images as an example of the results of a standard reconstruction approach.  

%--------------------------------------------------------
\section{Conclusions}\label{sec:conclusions}
%--------------------------------------------------------
The results presented in this paper demonstrate, on experimental data collected on several differing EIT systems, that the D-bar methods considered here are robust to modeling errors of domain shape and electrode placement.  They produce absolute images of very similar quality to difference images computed from the same nonhomogeneous data sets making use of basal saline data.  In particular, it was shown that even when the electrode centers are known so poorly that the resulting electrodes overlap, the D-bar method still produces distinguishable \tRev{absolute} images absent of the artefacts common in traditional reconstruction methods that rely on repeated solutions to the forward problem.  The results hold promise for clinical imaging, where certain applications may benefit from absolute images.  \tRev{For the clinical setting, further
investigation is needed to evaluate the effectiveness of this method in 3D settings.}

%--------------------------------------------------------
\section*{Acknowledgements}
%--------------------------------------------------------
The ACT3 and ACT4 data were provided by the EIT group at RPI [\url{https://www.ecse.rpi.edu/homepages/saulnier/eit/eit.html}], for which we express our thanks.  \tRev{We also thank Matteo Santacessaria for his helpful discussions.}

%%%%%%%%%%%%%%%%%
%--------------------------------------------------------
% BIBLIOGRAPHY
%--------------------------------------------------------
%\bibliographystyle{amsalpha}
%\bibliographystyle{plain}
\small
\bibliography{bibliographyRefs_use.bib}
%--------------------------------------------------------
\end{document}